\documentclass[opre,nonblindrev]{informs3_modified} 
\usepackage{amsmath}
\usepackage{amsfonts}
\usepackage{amssymb}
\usepackage{mathrsfs}
\usepackage{multirow}
\usepackage{xcolor}
\usepackage{enumitem}
\usepackage{algorithm}
\usepackage{algorithmic}
\usepackage{fancyhdr}
\OneAndAHalfSpacedXI 

\usepackage{endnotes}
\let\footnote=\endnote

\usepackage{natbib}
 \bibpunct[, ]{(}{)}{,}{a}{}{,}%
 \def\bibfont{\small}%
 \def\bibsep{\smallskipamount}%
 \def\bibhang{24pt}%
\TheoremsNumberedThrough     
\ECRepeatTheorems

\EquationsNumberedThrough    

\begin{document}


\RUNAUTHOR{Hu, Song, and Fu}

\RUNTITLE{Black-box Quantile Optimization}

\TITLE{Quantile Optimization via Multiple Timescale Local Search for Black-box Functions
}

\ARTICLEAUTHORS{%
\AUTHOR{Jiaqiao Hu}
\AFF{Department of Applied Mathematics and Statistics \\
State University of New York at Stony Brook, Stony Brook, NY 11794, \EMAIL{jqhu@ams.sunysb.edu}}

\AUTHOR{Meichen Song}
\AFF{Department of Applied Mathematics and Statistics \\
State University of New York at Stony Brook, Stony Brook, NY 11794, \EMAIL{meichen.song@stonybrook.edu}}

\AUTHOR{Michael C. Fu}
\AFF{Department of Decision, Operations \& Information Technologies\\
University of Maryland, College Park, MD
20742, \EMAIL{mfu@umd.edu}}
} 

\ABSTRACT{%
We consider quantile optimization of black-box functions that are estimated with noise.
We propose two new iterative three-timescale local search algorithms. 
The first algorithm uses an appropriately modified finite-difference-based gradient estimator that requires $2d$ + 1 samples of the black-box function per iteration of the algorithm, where $d$ is the number of decision variables (dimension of the input vector).
For higher-dimensional problems, this algorithm may not be practical if the black-box function estimates are expensive.
The second algorithm employs a simultaneous-perturbation-based gradient estimator that uses only three samples for each iteration regardless of problem dimension.
Under appropriate conditions, we show the almost sure convergence of both algorithms. {In addition, for the class of strongly convex functions, we further establish their (finite-time) convergence rate through a novel fixed-point argument.}
Simulation experiments indicate that the algorithms work well on a variety of test problems and compare well with recently proposed alternative methods.}



\maketitle

%

\section{Introduction}\label{sec1}
In black-box settings, only estimates of an output function are available,
i.e., there is minimal knowledge of the underlying system generating the output.
Furthermore, the output estimates might also contain observation noise.
For such problems, there is an extensive literature of algorithms addressing the case where the performance measure is an expectation, most commonly the mean, e.g., see \cite{fu2015} and references therein in the context of simulation optimization.
However, in many situations such as many risk management problems, one is interested in tail behavior of the output function or the median rather than the mean,
in which case the performance measure of interest is a quantile,
and the objective is quantile optimization.

Black-box optimization (BBO), defined by \cite{AuHa17} as ``the study of design and analysis of algorithms that assume the objective and/or constraint functions are given by blackboxes"
is a well-developed field in the deterministic (noiseless) setting.
Although the focus of the algorithms described and analyzed in their book are derivative-free approaches,
\cite{AuHa17} begin by strongly recommending that ``if gradient information is available, reliable, and obtainable at reasonable cost, then gradient-based methods should be used."

{In this paper, we consider the stochastic BBO setting where the goal is to optimize the quantile of a black-box output random variable. Our main assumption is that the quantile function is smooth enough so that gradient-based search will yield locally optimal solutions. Such a smoothness condition is guaranteed when the output distribution is differentiable (see equation~(\ref{gd0}) below), which is common in many engineering applications, ranging from queueing network optimization \citep[e.g.,][]{fuhill} to traffic simulation \citep{spallchin, lietal2017} to neural network (NN) training \citep[e.g.,][]{spall97,hong10}. 
For instance, the steady-state waiting time distribution in a queueing network is usually  differentiable with respect to the service rates of the nodes, whereas in traffic simulation, the distribution of vehicle travel time on a road network is typically a smooth function of traffic signal timings \citep[e.g.,][]{cao2014}. Nevertheless, due to the complexity of such problems and/or the lack of model details (e.g., when a simulation program is coded using an off-the-shelf commercial package, or knowledge transfer from the builder to the user of a model is not properly assured), directly estimating gradients based on exploiting model structure is sometimes difficult or practically infeasible. Moreover, in certain applications such as NN controller-design for systems with unknown dynamics \citep{spallchin,spall97}, it is not even possible to determine the gradient of a loss function through direct gradient techniques. In these circumstances, the underlying system model is essentially treated as a black-box, and}
finite-difference(FD)-based estimates of the gradient are often used to carry out local search, dating back to the Kiefer-Wolfowitz stochastic approximation (SA) algorithm \citep{KiWo52, kucl78, kuyin97}.
However, there are very few algorithms in general for what could be called quantile BBO,
in contrast to the abundance of local search algorithms for mean performance.

We will be relying on (approximate) quantile gradient estimates based only on output function samples.
To be specific, let $Y(\theta)$ denote the output random variable and $\theta \in \Theta\subseteq \Re^d$ the set of input decision variables, which we will refer to henceforth as the (input) parameter vector,
which can include both distributional and structural parameters, meaning that the elements of $\theta$ may affect the black-box function both directly and via the input distributions.
The usual optimization problem is of the form
$
\min_{\theta\in \Theta} E[Y(\theta)],
$
whereas we consider the optimization problem:
\begin{equation}\label{obj}
\min_{\theta\in \Theta}q_\varphi(\theta),
\end{equation}
where the quantile function $q_\varphi(\theta)$ is defined by
$$
P(Y(\theta) \leq q_\varphi(\theta)) = \varphi,~~ \varphi \in (0,1).
$$
Since $\varphi$ is fixed throughout this paper, its dependence will be dropped henceforth to simplify notation, i.e., the quantile will simply be denoted by $q(\theta)$ or sometimes just $q$.

Under appropriate smoothness conditions on $q(\theta)$, solving (\ref{obj}) {essentially} becomes equivalent to finding the zero of the gradient $\nabla q(\theta)$, so
a gradient-based iterative local search algorithm would take the general form:
\begin{equation}\label{eqSA}
\theta_{k+1}= \theta_k-\alpha_k \widehat{\nabla} q_k(\theta_k) ,
\end{equation}
where $\widehat{\nabla}q(\cdot)$ denotes an estimator of the quantile gradient,
which is the key element in defining the algorithm.
A straightforward symmetric finite-difference (SD) estimator for the quantile gradient would take the following form:
\begin{equation}\label{qFD}
\widehat{\nabla}_i q(\theta) = \frac{\hat{q}(\theta+c e_i)-\hat{q}(\theta-c e_i)}{2c},
~i=1,\ldots,d,
\end{equation}
where $\widehat{\nabla}_i$ denotes the $i$th component of the gradient estimator,
$e_i$ denotes the unit vector in the $i$th direction, and $c>0$.
This is essentially the approach taken by \cite{kibzun12}.
Note that one could also consider one-sided (forward or backward) finite-difference estimators.
However, one challenge that is apparent in the estimator (\ref{qFD}) is that unlike in the mean case,
the two difference terms in the numerator of the quantile finite-difference gradient are not themselves unbiased, only consistent, which means that the iteration sample size would eventually have to increase to infinity to guarantee convergence. Furthermore, each iteration would require calculation of $\hat{q}$ using order statistics, which may be computationally impractical.
Our alternative approach is to use two additional iterative updates for $q(\theta)$
and ${\nabla} q(\theta)$ based on the following result.

Assume that the output function $Y(\theta)$ is a continuous random variable with
cumulative distribution function (c.d.f.) $F$ and probability density function (p.d.f.) $f$,
we make use of the following relationship \citep{fu2009conditional}:
\begin{equation}\label{gd0}
\nabla_{\theta}q(\theta) = -\frac{\nabla_{\theta}F(y;\theta)|_{y=q}} {f(q;\theta)}.
\end{equation}
The simplest way to use this would be to solve for the zero of the numerator,
and use the SD estimator analogous to (\ref{qFD}) for the c.d.f. gradient:
\begin{equation}\label{F_FD}
\widehat{\nabla}_i F(q; \theta) = \frac{I\{Y(\theta +c e_i)\leq  q\} - I\{Y(\theta -c e_i)\leq  q\}}{2c},
~i=1,\ldots,d,
\end{equation}
where $I\{\cdot\}$ denotes the indicator function,
in which case (\ref{eqSA}) becomes a two-timescale SA algorithm:
\begin{eqnarray}\label{eqSA1}
\label{eqSA1a}
\theta_{k+1} &=& \theta_k + \alpha_k \widehat{\nabla}_\theta F(q_k; \theta_k) ,  \\
q_{k+1} &=& q_k{+}\gamma_k (\varphi - I\{Y_k \leq  q_k\}) ,  \label{eqSA1b}
\end{eqnarray}
%
where (\ref{eqSA1b}) is a recursive quantile estimator replacing the classical sample quantile based on order statistics (see Section~\ref{sec2}), and to make the algorithmic convergent, the perturbation $c$ in (\ref{F_FD}) would also need to go to zero as $k \rightarrow \infty$.
However, empirical results indicate that this approach does not work well in practice, especially when the quantile level $\varphi$ is close to 1 (or 0), in which case the two indicator terms in (\ref{F_FD}) will simultaneously take the value 1 (or 0) with high probability. As a result, a large number of iterations need to be performed in order to obtain a meaningful (non-zero) estimate of the c.d.f. gradient. Moreover, since (\ref{F_FD}) only provides an estimate for the direction of the quantile gradient, not its magnitude, the approach may not be useful in other related applications such as robustness assessment and quantile sensitivity analysis.


The algorithms we propose in this work follow the general structure of (\ref{eqSA1}) and (\ref{eqSA1b}); however, our goal is to estimate the ``real'' quantile gradient (\ref{gd0}) rather than just $\nabla_{\theta}F(y;\theta)|_{y=q}$ in the numerator. As in (\ref{F_FD}), each step of our proposed SD-based estimator requires $2d$+1 function evaluations. When the number of decision variables is large, the SD estimator may become computationally impractical, so we introduce a second algorithm based on using simultaneous perturbation (SP) gradient estimator along the lines of \cite{spall92}, which uses only three function evaluations at each iteration, independent of input parameter dimension. Both the SD and SP estimators require significant adjustments to handle the quantile setting.





As alluded to earlier, the literature on quantile optimization in the stochastic BBO setting is very sparse,
and we now review the most closely related work.
The most closely related work to ours is \cite{kibzun12}, cited earlier, which proposes a stochastic quasi-gradient (QG) algorithm for convex quantile objectives by estimating quantile gradients via a ``traditional" symmetric difference approximation.
Also relevant to the BBO setting are the derivative-free methods using the Bayesian optimization approach, e.g., \cite{ng2021,saba21}, which employ a surrogate model to approximate the response surface of the unknown quantile function.
Lastly, the multi-timescale stochastic approximation (SA) procedure
developed in \cite{huetal2021} has the same structure proposed in our work, but the algorithm cannot be applied in the BBO setting, because it relies on the availability of direct gradients (e.g., through techniques such as perturbation analysis or the likelihood ratio method) that are not available in a black-box setting, as knowledge of the underlying system is needed to derive the gradient estimators, whereas our algorithms use only the black-box function outputs.

Some other related work, albeit much less relevant to our BBO setting, are algorithms that rely on knowledge of the output distribution.
These include the mathematical programming approaches presented in \cite{kibzun91,kibzun13} and \cite{vas15}, and the scenario optimization method of \cite{zamar}.
For differentiable problems, there are also approaches that use gradient information, such as \cite{kim11}, who propose a recursive gradient algorithm for a special class of heavy-tailed distributions that admits the interchange of the derivative and quantile function.

Under appropriate conditions, we analyze the bias and variance of the proposed quantile gradient estimators and establish the almost sure local convergence of
the two FD-based algorithms---SD quantile optimization (SDQO) and
SP quantile optimization (SPQO)---for general multi-modal problems.
Most importantly, {for the class of problems with strongly convex
objective functions,} we are able to analyze the (finite-time) convergence {\it rate} of the algorithms by introducing a novel fixed-point argument. The key idea is to bound the algorithm's estimation errors through the composition of a sequence of suitably constructed contraction mappings, so that the convergence rates of quantile/gradient estimates can be characterized in detail by inspecting the solutions to a collection of fixed-point equations.
As far as we are aware, these are the first quantile BBO algorithms with both guaranteed convergence and known rate of convergence.
While the convergence rate of single-timescale SA is well understood in the literature \citep[cf., e.g.,][]{fabian,spall92,kuyin97,borkar08}, the rate analysis for multi-timescale SA algorithms has been a long-standing open research challenge.
The only existing results seem to be \cite{konda04} and \cite{mokk06} for two-timescale SA algorithms. Our algorithms operate on three timescales, and the convergence rate study of such SA algorithms has not been addressed.
Moreover, the fixed-point argument presented in this work is by no means limited to the analysis of these algorithms, but provides a new general approach that can potentially be applied to address the convergence rate issues of other multi-timescale SA algorithms.

In sum, we view our work as making the following research contributions:
\begin{itemize}
\item
We introduce new FD-based local search algorithms, SDQO and SPQO,
for optimizing a black-box quantile function, prove their convergence, and characterize their convergence rate.
\item
In terms of theory, the convergence rate analysis is the first such results for a three-timescale SA algorithm, and the fixed-point argument used in the analysis is a new general approach that can be applied to other multi-timescale SA algorithms.
\item
In terms of practice, SPQO is particularly well suited to high-dimensional problems, because the number of black-box evaluations per iteration is independent of the number of decision variables (optimization input parameters).
\item
Lastly, the new algorithms provide a practical complement to existing global optimization algorithms that primarily use metamodeling/surrogate functions for BBO.
\end{itemize}

The rest of this paper is organized as follows.
Section~\ref{sec2} begins with an intuitive motivation for the two FD-based black-box quantile gradient estimators and then presents the SD/SP estimators, along with their corresponding optimization algorithms SDPO and SPQO, with a detailed discussion of the proposed simultaneous perturbation estimator.
The convergence and convergence rate analyses of the algorithms are provided in Sections~\ref{sec3} and \ref{sec4}, respectively.
In Section~\ref{sec5}, simulation experiments are used to illustrate and test the performance of the algorithms,
and Section~\ref{sec6} provides some conclusions and future research.

\section{New FD-based Quantile BBO Algorithms}\label{sec2}

We begin with an intuitive informal derivation of the general form of a FD gradient estimator, which will then be specialized to the SD and SP gradient estimators, to be analyzed more rigorously.
To simplify the discussion (and notation) here, we consider the case where $\theta$ is a scalar ($d=1$), so we seek to estimate the derivative $q'(\theta)$, in which (\ref{gd0}) can be viewed as the ratio of two derivatives:
$$
q'(\theta) = -\frac{\partial_2 F(q;\theta)}{\partial_1 F(q;\theta)}, \mbox{~defining~} \partial_i F(q;\theta) \equiv \partial_i F(x;\theta) |_{x=q},
$$
where $\partial_i$ denotes the derivative w.r.t. the $i$th argument and the latter definition is for notational convenience.
When enough is known about the system to develop direct derivative estimators for $\partial_1$ and $\partial_2$, {i.e., the setting considered in \cite{huetal2021},} a natural approach to estimate $q'(\theta)$ would be to solve the equivalent root-finding problem:
\begin{equation}
\label{root0}
q'(\theta) \partial_1 F(q;\theta) + \partial_2 F(q;\theta) = 0.
\end{equation}
Assuming that direct derivative estimators are not available,
{it turns out that a straightforward extension to using FD estimates of
each of these derivatives would not only be computationally burdensome but could also lead to numerical difficulties. In particular, as mentioned in Section~\ref{sec1}, an FD estimator such as (\ref{F_FD}) for either $\partial_1$ or $\partial_2$ would frequently yield a value of 0, in which case the root-finding equation (\ref{root0}) is not even well-posed. Therefore,
we instead motivate an alternative estimator for approximating the entire left-hand side of (\ref{root0})} by considering a simple first-order Taylor series expansion of $F$ in the two arguments:
$$
F^{\pm} \equiv F(q \pm \Delta q; \theta \pm \Delta \theta)
 = F(q;\theta) \pm \partial_1 F(q;\theta) \Delta q  \pm \partial_2 F(q;\theta) \Delta \theta,
$$
where we are ignoring higher-order terms for now, but these arguments will be made more formal shortly.
Taking the difference,
$$
F^+ - F^-
= 2 \partial_1 F(q;\theta) \Delta q + 2 \partial_2 F(q;\theta) \Delta \theta
= 2 \Delta \theta [ \partial_1 F(q;\theta) q'(\theta) + \partial_2 F(q;\theta) ],
$$
where we have taken $\Delta q = q'(\theta) \Delta \theta$,
so solving $(F^+ - F^-)/2\Delta \theta = 0$ is equivalent to the root-finding problem using direct gradients given above by (\ref{root0}).
Noting that  $F^{\pm} = E[I\{Y^{\pm} \leq q + q'(\theta) \Delta \theta]$,
where $Y^{\pm} \sim F(\cdot; \theta \pm \Delta \theta)$,
it thus motivates the two coupled root-finding equations that must be solved:
\begin{eqnarray}
E \left[ \frac{-I\{Y^+ \leq q + q'(\theta) \Delta \theta \} + I\{Y^- \leq q - q'(\theta) \Delta \theta \}}{2 \Delta \theta} \right] \label{eq:rf1}
&=& 0, \\
E \left[ I\{Y \leq q  \}  \right] &=& \varphi, \label{eq:rf2}
\end{eqnarray}
where the second equation is solved via the SA iteration (\ref{eqSA1b}),
and the first equation will be incorporated into a new SA iteration to serve as the gradient estimator in the SA iteration (\ref{eqSA1a}). 

\subsection{SD/SP Quantile Gradient Estimators}

We now return to the multi-dimensional ($d>$1) setting and provide two versions of the FD scheme just described. Both versions can be viewed as different implementations of the SA method for numerically solving the two coupled stochastic root-finding equations (\ref{eq:rf1}) and (\ref{eq:rf2}).

Denote by $\theta^*$ an optimal solution to (\ref{obj}) and let $\hat\theta_k$ be an estimate of $\theta^*$.
Let $\hat\theta_k$ be fixed, $\hat q_k$ and $\hat D_k$ be the current estimates of $q(\hat \theta_k)$ and $\nabla_{\theta}q(\theta)|_{\theta=\hat \theta_k}$.
The SD estimator we propose simultaneously computes new estimates of the quantile and its gradient as follows:
\begin{align}
\hat D_{k+1}&=\hat D_k+\frac{\beta_k}{2c_k}
\begin{pmatrix}
{-I\{Y(\hat\theta_k+c_ke_1)\leq \hat q_k+c_k \hat D_k^Te_1\}+I\{Y(\hat\theta_k-c_ke_1)\leq \hat q_k-c_k \hat D_k^Te_1 \}}\\
 \vdots\\
 {-I\{Y(\hat\theta_k+c_ke_d)\leq \hat q_k+c_k \hat D_k^Te_d\}+I\{Y(\hat\theta_k-c_ke_d)\leq \hat q_k-c_k \hat D_k^Te_d \}}
 \end{pmatrix} \label{hatdkt} \\
\hat q_{k+1}&=\hat q_k+\gamma_k \big(\varphi-I\{Y(\hat\theta_k)\leq \hat q_k \}\big), \label{hatqk}
\end{align}
where $\beta_k,\,\gamma_k >0$ are step-sizes, $c_k>0$ is the perturbation size, and $Y(\hat\theta_k\pm c_ke_i)$ (for $i=1,\ldots,d$) are output random variables obtained by perturbing the $i$th element of $\hat\theta_k$ while holding all other components unchanged. Clearly, each step of (\ref{hatdkt}) requires $2d$ function evaluations. This, together with $Y(\hat\theta_k)$ needed in (\ref{hatqk}) for quantile estimation, results in a total of $2d+1$ function evaluations per iteration of the procedure.

The SP estimator, on the other hand, simultaneously varies all components of the underlying parameter vector in random directions, so that the same effect of the SD scheme can be achieved with only three function evaluations. Compared with the $2d+1$ per-iteration complexity of the SD estimator, this has the potential to lead to significant savings in computational cost, especially when the problem dimension is high and/or black-box function evaluations are expensive.
Let $\theta_k$, $q_k$, and $D_k$ denote the respective SP estimates for $\theta^*$, $q(\theta_k)$, and $\nabla_{\theta}q(\theta)|_{\theta=\theta_k}$. The estimator can be compactly expressed in the following recursive form:
\begin{align}
D_{k+1}&=D_k+\beta_k\bigg[\frac{-I\{Y(\theta_k+c_k\Delta_k)\leq q_k+c_kD_k^T\Delta_k \}+I\{Y(\theta_k-c_k\Delta_k)\leq q_k-c_kD_k^T\Delta_k \}}{2c_k\Delta_k} \bigg] \label{dkt} \\
q_{k+1}&=q_k+\gamma_k \big(\varphi-I\{Y(\theta_k)\leq q_k \}\big), \label{qktmp}
\end{align}
where $\Delta_k=(\Delta_{k,1},\ldots,\Delta_{k,d} )^T$ is a zero-mean random direction with i.i.d. components, and the division by the vector $\Delta_k$ is understood to be element-wise.

The above estimators depart significantly from the usual SD/SP formulations in several  different aspects: $i$) unlike (\ref{qFD}), they involve the difference quotient of an indicator function rather than that of the quantile function whose gradient is sought; $ii$) both are iterated, rather than one-shot (as in (\ref{qFD}) and (\ref{F_FD})), procedures in which the gradient estimation is coupled with another iterative process for estimating the quantile; $iii$) in contrast to conventional SD/SP, where only the parameter vector $\theta$ is varied (see, e.g., (\ref{qFD}) and (\ref{F_FD})), the quantile estimate $\hat q_k$ (resp. $q_k$) is randomly perturbed in (\ref{hatdkt}) (resp. (\ref{dkt})) at the same time, with the magnitude of the perturbation being directly affected by the gradient estimate $\hat D_k$ (resp. $D_k$) itself. As we will see shortly, this last difference further leads to other differences in algorithm design and analysis.

We now provide additional validation for why these estimators work, formalizing the intuitive derivation outlined in the beginning of the section.
We focus on the SP estimator and consider (\ref{dkt}) and (\ref{qktmp}) in their deterministic forms. The SD estimator works in a completely analogous way, so most of the arguments for the SP case also carry over to the SD estimator.
Note that conditional on $\theta_k$, $q_k$, $D_k$, and $\Delta_k$, the expectations of the two indicator terms in (\ref{dkt}) are given by $F(q_k+c_kD_k^T\Delta_k;\theta_k+c_k\Delta_k)$ and $F(q_k-c_kD_k^T\Delta_k;\theta_k-c_k\Delta_k)$.
A two-variable third-order Taylor series expansion of these two functions around $(q_k,\theta_k)$ then yields
\begin{align}\label{bias} \nonumber
\frac{-F(q_k+c_kD_k^T\Delta_k;\theta_k+c_k\Delta_k)+F(q_k-c_kD_k^T\Delta_k;\theta_k-c_k\Delta_k)}{2c_k\Delta_k}\\
&\hspace{-7cm}= \frac{-2f(q_k,\theta_k)c_k\Delta_k^TD_k-2\nabla_{\theta}^TF(q_k;\theta)_{\theta=\theta_k}c_k\Delta_k}{2c_k\Delta_k}+O(c_k^2), \end{align}
where the big-O notation signifies the order of a term, which is formally defined in Section~\ref{sec3}.
Thus, by applying the key argument of SP theory (i.e., $E[{\Delta_{k,i}}/{\Delta_{k,j}}]=0$ for all $i\neq j$) and ignoring the higher-order bias term $O(c_k^2)$ in (\ref{bias}), it is not difficult to observe that the expected-value version of (\ref{dkt}) (with the difference quotient there replaced by its conditional expectation given  $\theta_k$, $q_k$, and $D_k$, but excluding $\Delta_k$) can be written as
\begin{equation}\label{meanflow}
D_{k+1}=D_k+\beta_k\big(-f(q_k;\theta_k)D_k-\nabla_{\theta}F(q_k;\theta)|_{\theta=\theta_k}\big).
\end{equation}
Equation (\ref{meanflow}) is a fixed-point iteration for solving  $-f(q_k;\theta_k)D-\nabla_{\theta}F(q_k;\theta)|_{\theta=\theta_k}=0$ for $D$,
which has solution
$-\nabla_{\theta}F(q_k;\theta)|_{\theta=\theta_k}/f(q_k;\theta_k)$
in exactly the same form as (\ref{gd0}) with the true quantile $q(\theta_k)$ replaced by its estimate $q_k$.
A similar interpretation also applies to (\ref{qktmp}), and by noting that $E[I\{Y(\theta_k)\leq q_k \}|\theta_k,q_k ]=F(q_k;\theta_k)$, the sequence $\{q_k\}$ can be seen to track the unique solution $q(\theta_k)$ to the root-finding problem $\varphi-F(q;\theta_k)=0$. Consequently, as $q_k$ tends to $q(\theta_k)$, it is reasonable to expect that $D_k$ will provide a close approximation to the true gradient $\nabla_{\theta}q(\theta)|_{\theta=\theta_k}$.

The preceding developments ignored the fact that due to the extra perturbations $\pm c_kD_k^T\Delta_k$ introduced in (\ref{dkt}), the iterate $D_k$ itself is contained in the higher-order term in (\ref{bias}) (see the proof of Lemma~\ref{lem:bias}); thus, the sequence $\{D_k\}$ could in fact increase in magnitude to negate the claimed $O(c_k^2)$ order of the term.
This is a technical issue that does not occur in the usual mean-based setting, where the perturbation size is solely determined by $c_k$, so that the order of the bias can be bounded uniformly even without explicitly requiring the boundedness of the iterates \citep[cf. proof of Lemma 1 in][]{spall92}. 
To address this issue, we instead consider
a slight variant of (\ref{dkt}) that replaces $c_k$ in the difference quotient by a perturbation size that adapts to the magnitude of $D_k$.
Specifically, let $M_k=\max\{1, \|D_k\|/\sqrt{d} \}$ and define $\bar c_k=c_k/M_k$.
We suggest the following modification of (\ref{dkt}):
\begin{align}\label{dktmp}
D_{k+1}&=D_k+\beta_k\bigg[\frac{-I\{Y(\theta_k+ \bar c_k{\Delta}_k)\leq q_k+
\bar c_k D_k^T\Delta_k \}+I\{Y(\theta_k- \bar c_k{\Delta}_k)\leq q_k-\bar c_k D_k^T\Delta_k \}}{2\bar c_k{\Delta}_k} \bigg].
\end{align}
It can be easily seen that (\ref{dktmp}) serves the same estimation purpose as (\ref{dkt}) in the sense that its ``mean flow'' (i.e., deterministic counterpart), modulo the higher-order error terms, is identical to (\ref{meanflow}). Nevertheless, because $|\bar c_k D^T_k\Delta_k|\leq c_k\sqrt{d}\|\Delta_k\|$, the use of $\bar c_k$ in (\ref{dktmp}) prevents the perturbations in $q_k$ from becoming excessively large, and thus reduces the influence of $D_k$ on the resulting estimation bias.
We show in Section~\ref{sec:convergence} that under reasonable conditions, the sequence $\{D_k\}$ generated by (\ref{dktmp}) remains bounded, both almost surely and in second-order moments, which in effect justifies the $O(c_k^2)$ bias of the proposed estimator.

The constructions of our SD/SP estimators are based on a symmetric difference scheme. It is possible to consider alternative estimators relying on one-sided difference that use $d+1$ (resp. two) simulation evaluations per iteration. For example, in the SP case, the difference quotient in (\ref{dktmp}) could instead be replaced by either
\begin{align*}
&\frac{-I\{Y(\theta_k)\leq q_k+\bar c_k D_k^T\Delta_k \}+I\{Y(\theta_k-\bar c_k\Delta_k)\leq q_k \}}{\bar c_k{\Delta}_k}~\mbox{or~}  \\
& \frac{-I\{Y(\theta_k+\bar c_k\Delta_k)\leq q_k+\bar c_k D_k^T\Delta_k \}+I\{Y(\theta_k)\leq q_k \}}{\bar c_k{\Delta}_k}.
\end{align*}
Both approaches would lead to the same desired effect for estimating the quantile gradient but come at the cost of introducing large biases of order $O(c_k)$ compared to the $O(c_k^2)$ bias in (\ref{dktmp}).


\subsection{SA Local Search Algorithms}

The quantile optimization algorithms we propose integrate $\hat D_k$ and $D_k$ into the standard gradient descent method, as presented below,
where $\Pi_{\Theta}(\cdot)$ denotes a projection operator that brings an iterate back onto the parameter space $\Theta$ whenever it becomes infeasible. \\[-15pt]

\begin{algorithm}[H]
\renewcommand{\thealgorithm}{}
\caption{\bf SDQO}
\begin{description}
\item[\bf~~Input:] initial estimates $\hat q_0$, $\hat D_0$, $\hat \theta_0$;  sequences $\{\alpha_k\}$, $\{\beta_k\}$, $\{\gamma_k\}$, $\{c_k\}$;
\item[\bf~~Initialize:] iteration counter $k\gets 0$;
\item[\bf~~Loop] until a stopping criterion is met:
\begin{align}
&\hspace{-5mm}\tilde c_k =c_k/\max\{1, \|\hat D_k\|/\sqrt{d} \},  \label{c_tilde} \\
&\hspace{-5mm}\hat q_{k+1}=\hat q_k+\gamma_k \big(\varphi-I\{Y(\hat \theta_k)\leq \hat q_k \}\big), \nonumber\\
&\hspace{-5mm}\hat D_{k+1}=\hat D_k+\frac{\beta_k}{2\tilde c_k}
\begin{pmatrix}
{-I\{Y(\hat\theta_k+\tilde c_k e_1)\leq \hat q_k+\tilde c_k D_k^Te_1\}+I\{Y(\hat\theta_k-\tilde c_k e_1)\leq \hat q_k-\tilde c_k D_k^Te_1 \}}\\
 \vdots\\
 {-I\{Y(\hat\theta_k+\tilde c_k e_d)\leq \hat q_k+\tilde c_k D_k^Te_d\}+I\{Y(\hat\theta_k-\tilde c_k e_d)\leq \hat q_k-\tilde c_k D_k^Te_d \}}
 \end{pmatrix}, \label{vvhatdk}\\
&\hspace{-5mm}\hat \theta_{k+1}=\Pi_{\Theta}\big(\hat \theta_k-\alpha_k \hat D_k \big), \nonumber \\
&\hspace{-5mm} k \gets k+1; \nonumber
\end{align}
\end{description}
\vspace*{-8pt}
\end{algorithm}

Note that (\ref{vvhatdk}) is a variant of (\ref{hatdkt}), in which the same substitution as in (\ref{dktmp}) has been made, i.e., with $\tilde c_k$ defined by (\ref{c_tilde}) replacing $c_k$ in (\ref{hatdkt}).\\[-2mm]

\begin{algorithm}
\renewcommand{\thealgorithm}{}
\caption{\bf SPQO}
\begin{description}
\item[\bf ~~Input:] {initial estimates $q_0$, $D_0$, $\theta_0$;  sequences $\{\alpha_k\}$, $\{\beta_k\}$, $\{\gamma_k\}$, $\{c_k\}$.}
\item[\bf ~~Initialize:] iteration counter $k \gets 0$;
\item[\bf ~~Loop] until a stopping criterion is met:
\begin{align}
&\hspace{-5mm}\bar c_k=c_k/\max\{1, \| D_k\|/\sqrt{d} \}, \nonumber \\
&\hspace{-5mm}q_{k+1}=q_k+\gamma_k \big(\varphi-I\{Y(\theta_k)\leq q_k \}\big), \label{qk}\\
&\hspace{-5mm}D_{k+1}=D_k+\beta_k\bigg[\frac{-I\{Y(\theta_k+\bar c_k\Delta_k)\leq q_k+\bar c_k D_k^T\Delta_k \}+I\{Y(\theta_k-\bar c_k\Delta_k)\leq q_k-\bar c_k D_k^T\Delta_k \}}{2\bar c_k{\Delta}_k} \bigg], \label{dk} \\
&\hspace{-5mm}\theta_{k+1}=\Pi_{\Theta}\big(\theta_k-\alpha_kD_k \big), \label{thetak}\\
&\hspace{-5mm}k\gets k+1; \nonumber
\end{align}
\end{description}
\vspace*{-8pt}
\end{algorithm}

Both SDQO and SPQO have a multi-timescale structure, as reflected by the use of distinct step-sizes $\gamma_k$, $\beta_k$, and $\alpha_k$ in the recursions. Intuitively speaking, since our discussion on the convergence behavior of the SP estimator has assumed a fixed value for $\theta_k$, the step-size $\alpha_k$ in SPQO should be chosen very small relative to $\beta_k$ and $\gamma_k$. As a result, when viewed from (\ref{dk}) and (\ref{qk}), the increment in $\theta_k$ at each step of (\ref{thetak}) is almost negligible as if the parameter vector $\theta_k$ were held at a constant value. On the other hand, because (\ref{dk}) is designed to iteratively approximate
$-\nabla_{\theta}F(q_k;\theta)|_{\theta=\theta_k}/f(q_k;\theta_k)$, the step-size $\beta_k$ should be taken to be the largest to warrant proper tracking of the ratio as both $q_k$ and $\theta_k$ vary over time.
In the same manner, the three recursions in SDQO should also be carried out at different speeds,  with their step-sizes satisfying $\alpha_k=o(\gamma_k)$ and $\gamma_k=o(\beta_k)$.

When the black-box function is given by a computer simulation program,
it is natural to exploit the use of common random numbers (CRN) \citep[e.g.,][]{law}---in the same spirit as the use of CRN in e.g. simultaneous perturbation stochastic approximation (SPSA) \citep{spallcrn}---for reducing the variance in the difference estimates.
With a slight abuse of notation, let $Y(U_k; \hat\theta_k\pm\tilde c_k e_i)$ be the output random variables simulated using the same input random number stream $U_k$ under the perturbed vectors
$\hat \theta_k\pm \tilde c_k e_i$ ($i=1,\ldots,d$).
The CRN version of SDQO simply works by replacing $Y(\hat \theta_k\pm \tilde c_k e_i)$ in (\ref{vvhatdk}) with $Y(U_k;\hat \theta_k\pm \tilde c_k e_i)$.
Likewise, a CRN version of SPQO can be easily implemented by substituting $Y(U_k;\theta_k\pm \bar c_k{\Delta}_k)$  for $Y(\theta_k\pm \bar c_k \Delta_k)$ in (\ref{dk}).
In Section~\ref{sec:convergence}, we give conditions under which we show that this CRN approach helps to induce a positive correlation between the indicator terms and reduces the estimation variance at each step.

\section{Convergence Results}\label{sec3}
To fix ideas and to avoid unnecessary repetition, we perform detailed analysis and present results mainly for the SPQO algorithm. All results obtained for SPQO can be shown (with appropriate/slight modifications) to hold for SDQO. We begin by defining $(\Omega, \mathcal{F},P)$ as the probability space induced by SPQO, where $\Omega$ is the set of all sample trajectories that could possibly be observed by executing the algorithm, $\mathcal{F}$ is the $\sigma$-field of subsets of $\Omega$, and $P$ is a probability measure on $\mathcal{F}$. We also define $\mathcal{F}_k=\sigma\{D_0,q_0,\theta_0,\ldots,D_k,q_k,\theta_k \}$ as an increasing $\sigma$-field representing the information available at iteration $k=0,1,\ldots.$ For a given vector $v$, let $\|v\|$ be the Euclidean norm of $v$; whereas for a matrix $A$, let $\|A\|$ be the matrix norm induced by the Euclidean norm. For any two real-valued functions $u(k)$ and $v(k)$, we write $u(k)=O(v(k))$ if $\limsup_{k\rightarrow \infty}{u(k)}/{v(k)}<\infty$ and $u(k)=o(v(k))$ if $\lim_{k\rightarrow \infty}{u(k)}/{v(k)}=0$.
We assume that the parameter space $\Theta$ is a compact convex set described by functional constraints and takes the form $\Theta=\{\theta\in \Re^d: h_j(\theta)\leq 0,~j=1,\ldots,m\},$ where $h_j(\cdot),~j=1,\ldots,m$ are continuously differentiable functions with their gradients satisfying $\nabla_{\theta}h_j(\theta)\neq 0$ whenever $h_j(\theta)=0$ \citep[cf., e.g.,][]{kuyin97}. Such a characterization is satisfied by many common choices for $\Theta$, including hyper-balls, hyper-rectangles and more general convex polytopes. For notational convenience, throughout our analysis,
we denote $I_k^+:=I\{Y(\theta_k+\bar c_k\Delta_k)\leq q_k+\bar c_k{D}_k^T\Delta_k \}$, $I_k^-:=I\{Y(\theta_k-\bar c_k\Delta_k)\leq q_k-\bar c_k{D}_k^T\Delta_k \}$, $F^+_k:=F(q_k+\bar c_k{D}_k^T\Delta_k;\theta_k+\bar c_k{\Delta}_k)$, and $F^-_k:=F(q_k-\bar c_k{D}_k^T\Delta_k;\theta_k-\bar c_k{\Delta}_k)$.

\subsection{Strong Convergence}\label{sec:convergence}

The projection $\Pi_{\Theta}(\cdot)$ in (\ref{thetak}) ensures the boundedness of $\theta_k$ by projecting an iterate onto the compact region $\Theta$. This operation can be defined through adding an extra correction term $Z_k$ to the recursion \cite[see, e.g.,][Chapter 5]{kuyin97}, leading to
\begin{align}
\theta_{k+1} = \theta_k - \alpha_k D_k + \alpha_k Z_k, \label{projection}
\end{align}
where $\alpha_k Z_k := \theta_{k+1} - \theta_k + \alpha_k D_k$ is the vector with the shortest Euclidean length needed to bring $\theta_k - \alpha_k D_k$ back onto $\Theta$. In our setting, since $\Theta$ is a convex set, $Z_k$ takes values in the convex cone generated by the inward normals to the surface of $\Theta$ at the point $\theta_{k+1}$, i.e.,
$Z_k\in -C(\theta_{k+1})$, where
\begin{equation}\label{cone}
C(\theta):=\{\nu\in \Re^d: \nu^T(\tilde{\theta}-\theta)\leq 0,~\forall\tilde{\theta}\in\Theta \}
\end{equation}
is the normal cone to $\Theta$ at $\theta$. Note that $C(\theta)=\{0\}$ whenever $\theta$ lies in the interior of $\Theta$.

The convergence of SPQO is investigated by following an ordinary differential equation (ODE)
argument \citep[e.g.,][]{kuyin97,borkar08,huetal2021}. The general idea is to construct interpolations of the iterates $\{D_k,q_k,\theta_k\}_{k=0}^{\infty}$ by ``stretching'' them continuously in time and then capture the long run behavior of these interpolations using a set of coupled ODEs. In particular, our main result is to show that the sequence $\{\theta_k\}$ generated by (\ref{projection}) asymptotically approaches the limiting solution to a projected ODE of the form
\begin{equation}\label{limitode}
\dot{\theta}(t)=-\nabla_{\theta}q_{\varphi}(\theta)|_{\theta=\theta(t)}+z(t),~~t\geq 0,
\end{equation}
where $z(t)\in -C(\theta(t))$ is the minimum force (the real vector with the smallest Euclidean norm) needed to keep the trajectory $\theta(t)$ within the constraint set $\Theta$.

We introduce the list of assumptions that will be used in our analysis.\\[-5mm]

\noindent \textbf{Assumptions:} \\
\textbf{A1:} {\it For almost all $(q_k,\theta_k)$ pairs, there exists an open neighborhood of $(q_k,\theta_k)$, independent of $k$ and $\omega\in\Omega$, such that

(a) ${\partial^2{f(y;\theta)}}/{\partial{y^2}}$, ${\partial^2 f(y;\theta) }/\partial y\partial\theta^T$, ${\partial^2{f(y;\theta)}}/{\partial{\theta^T}\partial{\theta^T}}$, and $\partial^3{F(y;\theta)}/\partial{\theta^T}\partial{\theta^T}\partial{\theta^T}$
all exist and are continuous on the neighborhood with their elements uniformly bounded in $k$ and $\omega$.

(b) $f(y,\theta)\geq\epsilon$ for all $(y,\theta)$ pairs in the neighborhood for some constant  $\epsilon>0$.}

\noindent\textbf{A2:} {\it  The random directions $\{\Delta_k\}$ are i.i.d., independent of $\mathcal{F}_k$. Each $\Delta_k$ has mutually independent components with the Bernoulli distribution $P(\Delta_{k,i}=1)=P(\Delta_{k,i}=-1)={1}/{2}$ for all $i=1,\ldots,d$.}

\noindent\textbf{A3:} {\it The sequences $\{\alpha_k\}$, $\{\beta_k\}$, $\{\gamma_k\}$, and $\{c_k\}$ satisfy the following conditions:

(a) $\beta_k,\,c_k >0$, $c_k\rightarrow 0$, $\sum_{k=0}^{\infty} \beta_k = \infty$, $\sum_{k=0}^{\infty} {\beta_k^2}/{c_k^2} < \infty$;

(b) $\gamma_k >0$, $\sum_{k=0}^{\infty} \gamma_k = \infty$;

(c) $\alpha_k >0$, $\sum_{k=0}^{\infty} \alpha_k = \infty$;

(d) $\alpha_k=o(\gamma_k),~\gamma_k=o(\beta_k).$}\\[-5mm]


Assumption A1(a) is consistent with the condition used in Lemma 1 of \cite{spall92} but is stated within a quantile optimization context. It requires the output distribution to be sufficiently smooth, and is satisfied when $F$ is three times continuously differentiable (in both arguments) with bounded derivatives. {Note that the condition can be weakened to twice differentiability when one-sided difference estimators are used.} From the discussion at the end of Section~\ref{sec2}, since (\ref{dk}) iteratively approximates a gradient of the form $-\nabla_{\theta}F(q_k;\theta)|_{\theta=\theta_k}/f(q_k;\theta_k)$, A1(b) ensures that the denominator of the ratio is bounded away from zero, so that the limit of the $\{D_k\}$ sequence (assuming its existence) does not get arbitrarily large; see \cite{ng2021} for a similar assumption. The suitability of A1(b) has been discussed in \cite{huetal2021}. Specifically, because $f(\cdot;\cdot)$ is continuous and $\Theta$ is compact, the condition holds trivially when $\{q_k\}$ lies in a compact set. In practice, this can be guaranteed by truncating the sequence to a large closed interval containing the true quantiles $q(\theta)$ for all $\theta\in\Theta$. It has been shown in \cite{huetal2021} that such a truncation will not have an influence on the convergence behavior of $\{q_k\}$.
Both A2 and A3 are conditions on the algorithm input parameters. The Bernoulli random direction is perhaps the most commonly used choice when implementing SP estimators.   A3 is also standard in the SA literature \citep[e.g.,][]{kucl78, kuyin97,spall92}. A3(d) is needed in multi-timescale SA methods \citep[cf., e.g.,][]{bhat05,borkar08,zhanghu,huetal2021}; it guarantees the three recursions to be performed at timescales that are noticeably distinct from each another (see the discussion at the end of Section~\ref{sec2}). The condition, when combined with $\sum_{k=0}^{\infty} {\beta_k^2}/{c_k^2} < \infty$ in A3(a), implies that $\sum_{k=0}^{\infty}\gamma_k^2<\infty$ and $\sum_{k=0}^{\infty}\alpha_k^2<\infty$.

We begin by stating a result that is essential for characterizing the convergence behavior of the algorithm. It shows that the gradient estimators constructed through (\ref{dk}) have finite second-order moments and that the sequence $\{D_k\}$ itself remains bounded almost surely for all $k$.


\begin{lemma}\label{lem:dk}
Assume that A1, A2, and A3(a) hold, then we have ($i$) $\sup_k E[\|D_k \|^2]<\infty$; ($ii$) $\sup_k \|D_k\|<\infty$ w.p.1.
\end{lemma}
\proof{Proof.}
See Section~\ref{appendix_dk} of the Appendix. \Halmos
\endproof

Lemma~\ref{lem:bias} below gives an explicit bound on the (conditional) bias introduced by the symmetric SP scheme used in (\ref{dk}). As a result of Lemma~\ref{lem:dk}, the estimation bias goes to zero at the rate  $O(c^2_k)$, both almost surely and in expectation.

\begin{lemma}\label{lem:bias}
Let Assumptions A1, A2, and A3(a) hold, and define the bias $$b_k(q_k,D_k,\theta_k):=E\Big[\frac{-I^+_k+I^-_k}{2\bar c_k{\Delta}_k}\Big|\mathcal{F}_k \Big]+f(q_k,\theta_k)D_k+\nabla_{\theta}F(q_k;\theta_k).$$
Then we have that
($i$) $b_k(q_k,D_k,\theta_k) =O(c^2_k)$ w.p.1; ($ii$) $E[\|b_k(q_k,D_k,\theta_k)\|]=O(c^2_k).$
\end{lemma}
\proof{Proof.}
From A2, $1/{\Delta}_k=\Delta_k$. It follows that
\begin{equation*}
E\Big[\frac{-I^+_k+I^-_k}{2\bar c_k{\Delta}_k}\Big|\mathcal{F}_k \Big]=\frac{M_k}{2c_k}E\big[(-I^+_k+I^-_k)\Delta_k\big|\mathcal{F}_k \big]=\frac{M_k}{2c_k}E\big[(-F^+_k + F^-_k )\Delta_k\big|\mathcal{F}_k \big].
\end{equation*}
Let $\nabla^3_{\theta}F(y;\theta):=\partial^3{F(y;\theta)}/\partial{\theta^T}\partial{\theta^T}\partial{\theta^T}$ be the tensor of $F$. Note that the tensor when evaluated at a vector $v$ of appropriate dimension, denoted by $\nabla^3_{\theta}F(y;\theta)[v]$, gives a matrix. Thus,
by a third-order Taylor series expansion of $F^+_k$ and $F^-_k$ around $(q_k,\theta_k)$ and then using
$\Delta^2_{k,i}=1$, $E[\Delta_{k,i}\Delta_{k,j}|\mathcal{F}_k]=0$ for all $i\neq j$, we obtain
\begin{align*}
E\Big[\frac{-I^+_k+I^-_k}{2\bar c_k{\Delta}_k}\Big|\mathcal{F}_k \Big]&=E\Big[\big(-f(q_k;\theta_k)\Delta^T_kD_k-\nabla^T_{\theta}F(q_k;\theta_k)\Delta_k\big)\Delta_k\big|\mathcal{F}_k \Big]+E\big[R_3(\bar q_k^+,\bar q_k^-,\bar{\theta}_k^+,\bar{\theta}_k^-)\Delta_k \big|\mathcal{F}_k\big]\\
&=-f(q_k;\theta_k)D_k-\nabla_{\theta}F(q_k;\theta_k)+E\big[R_3(\bar q_k^+,\bar q_k^-,\bar{\theta}_k^+,\bar{\theta}_k^-)\Delta_k\big|\mathcal{F}_k\big],
\end{align*}
where $\bar q_k^+$, $\bar q_k^-$ are on the line segments between $q_k$ and $q_k\pm \bar c_k{D}_k^T\Delta_k$, $\bar{\theta}_k^+$, $\bar{\theta}_k^-$ are on the line segments connecting $\theta_k$ and $\theta_k\pm \bar c_k{\Delta}_k$, and $R_3$ is a remainder term whose absolute value is bounded by
\begin{align*}
|R_3|&\leq \frac{c_k^2}{12M_k^2}\big(|f_{yy}(\bar q^+_k;\bar{\theta}_k^+)|+|f_{yy}(\bar q^-_k;\bar{\theta}_k^-)|\big)|D_k^T\Delta_k |^3 \\
&\hspace{1cm}+\frac{c_k^2}{12M_k^2}\Big(\big|\Delta_k^T\nabla^3_{\theta}F(\bar q^+_k;\bar{\theta}_k^+)[\Delta_k]\Delta_k \big|+\big|\Delta_k^T\nabla^3_{\theta}F(\bar q^-_k;\bar{\theta}_k^-)[\Delta_k]\Delta_k \big| \Big)\\
&\hspace{1.5cm}+\frac{c^2_k}{4M_k^2}\Big(\big|\Delta_k^T\nabla^2_{\theta}f(\bar q_k^+;\bar{\theta}_k^+)D_k^T\Delta_k\Delta_k \big|+\big|\Delta_k^T\nabla^2_{\theta}f(\bar q_k^-;\bar{\theta}_k^-)D_k^T\Delta_k\Delta_k \big|  \Big)\\
&\hspace{2cm}+\frac{c_k^2}{4M_k^2}\Big(\big|\nabla^T_{\theta}f_y(\bar q^+_k;\bar{\theta}_k^+)(D_k^T\Delta_k)^2\Delta_k \big|+\big|\nabla^T_{\theta}f_y(\bar q^-_k;\bar{\theta}_k^-)(D_k^T\Delta_k)^2\Delta_k \big|  \Big).
\end{align*}
Note that the facts $M_k\geq 1$, $\|\bar c_k D_k\|=c_k{\|D_k\|}/{M_k}\leq c_k\sqrt{d}$, and $\|\bar c_k {\Delta}_k\|\leq c_k \|\Delta_k\|=c_k\sqrt{d}$ ensure that the pairs $(\bar q_k^+,\bar{\theta}_k^+)$ and $(\bar q_k^-,\bar{\theta}_k^-)$ will all lie within the neighborhood of $(q_k,\theta_k)$ stated in A1 as $c_k\rightarrow 0$. Thus, invoking A1 and part ($ii$) of Lemma~\ref{lem:dk}, we arrive at the conclusion that $E[\|R_3(\bar q_k^+,\bar q_k^-,\bar{\theta}_k^+,\bar{\theta}_k^-)\Delta_k\||\mathcal{F}_k]=O(c_k^2)$ w.p.1. This completes the proof of part ($i$) of the lemma. Part ($ii$) follows from
a simple application of H\"{o}lder's inequality and the fact that $\sup_kE[\|D_k\|^2]<\infty$ (Lemma~\ref{lem:dk} part ($i$)).
\Halmos
\endproof

We now present the main convergence theorem. In its most general form, the result
implies that the local/global convergence of the algorithm can be determined by examining the asymptotic behavior of the ODE (\ref{limitode}). For example, if the ODE has multiple isolated stable equilibrium  points, then the sequence $\{\theta_k\}$ will converge to one of them.
A strengthened version of the result is obtained when the quantile function is strictly convex, in which case the unique optimal solution $\theta^*$ to (\ref{obj}) turns out to be a globally asymptotically stable equilibrium of (\ref{limitode}) \citep[see, Corollary 1 of][]{huetal2021}, so the sequence $\{\theta_k\}$ will converge to $\theta^*$ w.p.1.
Because the proof is similar to the convergence analysis in \cite{huetal2021}, it is included in the appendix.

\begin{theorem}\label{thm:con}
Assume that conditions A1, A2, and A3 hold. Then the sequence $\{\theta_k\}$ generated by SPQO converges to some limit set of the ODE (\ref{limitode}) w.p.1. In addition, if the objective function $q_{\varphi}(\theta)$ is strictly convex on $\Theta$, then the sequence $\{\theta_k\}$ converges to the unique optimal solution $\theta^*$ to (\ref{obj}) w.p.1.
\end{theorem}

Note that as with Lemma~\ref{lem:bias}, it is readily seen that the SD estimator in (\ref{vvhatdk}) also has $O(c_k^2)$ bias. Using this observation it can be shown, along the same lines as in the proof of Theorem~\ref{thm:con}, that the following (same) convergence result holds for SDQO.

\begin{theorem}\label{thm:sdqo}
Assume that condition A1 holds with $(q_k,\theta_k)$ replaced by $(\hat q_k,\hat \theta_k)$. Then under A3, the sequence $\{\hat \theta_k\}$ generated by SDQO converges to some limit set of the ODE (\ref{limitode}) w.p.1. In addition, if the objective function $q_{\varphi}(\theta)$ is strictly convex on $\Theta$, then the sequence $\{\hat \theta_k\}$ converges to the unique optimal solution $\theta^*$ to (\ref{obj}) w.p.1.
\end{theorem}

As mentioned in Section~\ref{sec3}, under the special setting of simulation optimization, the CRN versions of (\ref{vvhatdk}) and (\ref{dk}) could be used in some cases to improve the algorithm efficiency. We close this section by providing conditions that guarantee the variance reduction property of this approach. For simplicity, we only state and prove the result for (\ref{dk}). The variance reduction property of (\ref{vvhatdk}) is ensured under the same set of conditions given in Proposition~\ref{prop:crn} below.

\begin{proposition}\label{prop:crn}
Denote by $\hat I^\pm_k:=I\{Y(U_k;\theta_k\pm \bar c_k {\Delta}_k)\leq q_k\pm c_k\bar D_k^T\Delta_k \}$. Let conditions A1, A2, and A3(a) hold. Suppose that for any given parameter vector $\theta$, the output random variable $Y(U_k;\theta)$, when viewed as a function of input random numbers, is monotone in each argument, that is, $U_{k,i}\leq U'_{k,i}$ implies $Y(\ldots,U_{k,i-1},U_{k,i},U_{k,i+1},\ldots;\theta)\leq Y(\ldots,U_{k,i-1},U'_{k,i},U_{k,i+1},\ldots;\theta)$ for all $i$ or $Y(\ldots,U_{k,i-1},U_{k,i},U_{k,i+1},\ldots;\theta)\geq Y(\ldots,U_{k,i-1},U'_{k,i},U_{k,i+1},\ldots;\theta)$ for all $i$. Then we have
$$Var \Big(\frac{-\hat I^+_k +\hat I^-_k}{2\bar c_k{\Delta}_{k,i}}\Big|\mathcal{F}_k \Big)\leq Var \Big(\frac{-I^+_k + I^-_k}{2\bar c_k{\Delta}_{k,i}}\Big|\mathcal{F}_k \Big),~~\mbox{$i=1,\ldots,d,$}$$
for all $k$ w.p.1.
\end{proposition}
\proof{Proof.}

The proof relies on an inequality given on pp. 187-188 of \cite{comstat}, which can be stated as follows:
Let $X_1,\ldots,X_k$ be a sequence of i.i.d. random variables, $g_1(x_1,\ldots,x_k)$ and $g_2(x_1,\ldots,x_k)$ be two functions that are both monotonically non-decreasing (or non-increasing) in each argument. Then
\begin{equation}\label{pc}
E[g_1(X_1,\ldots,X_k)g_2(X_1,\ldots,X_k)]\geq E[g_1(X_1,\ldots,X_k) ]E[g_2(X_1,\ldots,X_k) ].
\end{equation}

Now fix an $i=1,\ldots,d$, it is straightforward to show that
\begin{align}\label{vdiff}
&Var\Big(\frac{-\hat I^+_k +\hat I^-_k}{2\bar c_k{\Delta}_{k,i}}\Big|\Delta_k,\mathcal{F}_k \Big)-
Var\Big(\frac{-I^+_k+ I^-_k}{2\bar c_k{\Delta}_{k,i}}\Big|\Delta_k,\mathcal{F}_k \Big) \nonumber\\
&=\frac{1}{2\bar c_k^2{\Delta}_{k,i}^2}\Big[E\big[I_k^+\big|\Delta_k,\mathcal{F}_k\big]E\big[ I_k^- \big|\Delta_k,\mathcal{F}_k \big]-E\big[\hat I_k^+\hat I_k^-  \big|\Delta_k,\mathcal{F}_k \big]\Big]. \nonumber
\end{align}
Since for fixed $q_k$, $D_k$, $\Delta_k$, and $\theta_k$, the indicator functions $I\{\cdot\leq q_k\pm \bar c_k D_k^T\Delta_k \}$ are non-increasing and by our assumption $Y(\cdot;\theta_k\pm \bar c_k{\Delta}_k)$ are monotone in each argument, the compositions $\hat I_k^{\pm}=I\{Y(\cdot;\theta_k\pm \bar c_k{\Delta}_k)\leq q_k\pm \bar c_k D_k^T\Delta_k  \}$ are also monotone in each argument. Hence, we have from (\ref{pc}) that $E[\hat I_k^+\hat I_k^-  |\Delta_k,\mathcal{F}_k ]\geq E[\hat I_k^+|\Delta_k,\mathcal{F}_k]E[\hat I_k^-|\Delta_k,\mathcal{F}_k]=E[ I_k^+|\Delta_k,\mathcal{F}_k]E[I_k^-|\Delta_k,\mathcal{F}_k]$. This shows that
$Var\big(\frac{-\hat I^+_k +\hat I^-_k}{2\bar c_k{\Delta}_{k,i}}\big|\Delta_k,\mathcal{F}_k \big)\leq
Var\big(\frac{-I^+_k+ I^-_k}{2\bar c_k{\Delta}_{k,i}}\big|\Delta_k,\mathcal{F}_k \big).$
Finally, by noticing that
$$Var\Big(E\big[\frac{-\hat I^+_k +\hat I^-_k}{2\bar c_k{\Delta}_{k,i}}\big|\Delta_k,\mathcal{F}_k \big]\Big|\mathcal{F}_k \Big)=
Var\Big(E\big[\frac{-I^+_k+ I^-_k}{2\bar c_k{\Delta}_{k,i}}\big|\Delta_k,\mathcal{F}_k \big]\Big|\mathcal{F}_k\Big),$$
the proof is completed by unconditioning on $\Delta_k$ using the law of total variance.
\Halmos
\endproof

Proposition~\ref{prop:crn} shows that when the simulation output random variables react monotonically with respect to the input random numbers, the conditional variance of the gradient estimator $D_{k+1}$ constructed using CRN is always no greater than that obtained under independent sampling. This monotonicity requirement can be expected in certain applications such as in the simulation of regenerative processes and queueing systems \cite[see, e.g.,][]{law}. Note that Proposition~\ref{prop:crn} is a finite time result that holds almost surely for every $k$. This is different from the work of \cite{spallcrn}, in which the same CRN approach has been used in SPSA
and shown to lead to a faster asymptotic convergence rate than the original SPSA without CRN. \cite{spallcrn}, however, consider mean-based SO, and a key condition used in deriving their result, when put into our current context, requires the two indictor functions $\hat I_k^{\pm}$ to be differentiable with bounded derivatives. So their result does not directly carry over to the quantile setting.

\section{Convergence Rate Analysis}\label{sec4}
Again, because SPQO differs from SDQO only in the perturbation scheme used in constructing gradient estimates (i.e., simultaneous vs. element-wise perturbation), we use SPQO as a representative algorithm and investigate its rate of convergence in detail. A completely analogous argument, which we omit, yields essentially the same rate result for SDQO (see Theorem~\ref{thm:ratesdqo} at the end of this section).
Our analysis assumes that the ODE (\ref{limitode}) has a unique globally asymptotically stable equilibrium $\theta^*$ that lies in the interior of $\Theta$.
Clearly, because $C(\theta^*)=\{0\}$, $\theta^*$ satisfies $\nabla_{\theta}q_{\varphi}(\theta)|_{\theta=\theta^*}=0$, and
according to Theorem~\ref{thm:con}, the sequence $\{\theta_k\}$ generated by (\ref{thetak}) converges to $\theta^*$ w.p.1.
The main result for SPQO is obtained through the repeated application of a fixed-point argument for characterizing the MAEs of SA estimates. 
The novelty of the approach resides in the use of a sequence of suitably designed contraction mappings to quantify the estimation errors accumulated over the iterations. This, in essence, translates a difficult rate analysis problem into the simple task of examining the fixed points of a sequence of contraction mappings. Notice that the three recursions in SPQO each represents a major category of SA algorithms, i.e., a stochastic rooting finding procedure (\ref{qk}), a Kiefer–Wolfowitz/SPSA type algorithm (\ref{dk}), and a Robbins–Monro type gradient iteration (\ref{thetak}). We show that these recursions, even when coupled through different timescales, can all be analyzed using the proposed fixed-point argument.

The analysis proceeds in three steps. First, we consider the two coupled iterations (\ref{qk}) and (\ref{thetak}) and
derive the convergence speed of the quantile recursion (\ref{qk}) while taking into account the variations in the value of $\theta_k$.
Then, we 
characterize the rate at which the gradient recursion (\ref{dk}) converges as both $q_k$ and $\theta_k$ vary over time. Finally, 
we present the main convergence rate result for (\ref{thetak}) and discuss the selection of algorithm parameters that optimizes the asymptotic performance of SPQO.
Throughout this section, we focus on standard step- and perturbation-sizes of the forms: $\alpha_k=a/k^{\alpha}$, $\beta_k=b/k^{\beta}$, $\gamma_k=r/k^{\gamma}$, and $c_k=c/k^\tau$, where $\alpha,\,\beta,\,\gamma,\,\tau\in (0,1)$ and $a,\, b,\, r,\,c>0$. Let $q(\theta)$ be twice continuously differentiable with Hessian matrix $H(\theta):=\nabla^2_{\theta}q(\theta)$. In addition to the assumptions used in Section~\ref{sec:convergence}, we also impose the following regularity conditions:\\[-5mm]

\noindent \textbf{Assumptions:} \\
\textbf{B1:} {\it For almost all $(q_k,\theta_k)$ pairs, there exist constants $\varepsilon,\,C_f>0$ such that $\varepsilon\leq f(y;\theta_k)\leq C_f$ for all $y$ in the interval between $q_k$ and $q(\theta_k)$.}

\noindent \textbf{B2:} {\it The output density $f(y;\theta)$ is jointly continuous in both $y$ and $\theta$. There are constants $L_f,\,L_F>0$ such that
$|f(y_1;\theta)-f(y_2;\theta)|\leq L_f\|y_1-y_2\|$ and
$\|\nabla_{\theta}F(y_1;\theta)-\nabla_{\theta}F(y_2;\theta) \|\leq L_F\|y_1-y_2 \|$ for all $\theta\in\Theta$.}

\noindent \textbf{B3:} {\it Let $\lambda(\theta)$ be the smallest eigenvalue of $H(\theta)$. There is a constant $\varrho>0$ such that $\lambda(\theta)\geq \varrho$ for all $\theta$ that lies on the line segment between $\theta_k$ and $\theta^*$.}\\[-5mm]

\noindent Note that since no knowledge of the bounding constants $\epsilon$ and $C_f$ is required, condition B1 is acceptable in many practical situations. For example, when the simulation outputs themselves are bounded or truncated to a large interval, it is easy to see from (\ref{qk}) that $\{q_k\}$ will stay bounded. Thus, the assumption holds if $f(\cdot;\cdot)$ is continuous on $\Theta$ and a compact interval that contains $\{q_k\}$ and $\{q(\theta_k)\}$ (see also the remarks on Assumption A1(b) in Section~\ref{sec:convergence}). B2 is roughly a globalized version of A1(a) but without requirements on the higher-order derivatives of $F$. A sufficient condition for B2 to hold is that the output distribution $F$ is twice differentiable in both arguments and has bounded derivatives.
B3, in a sense, is a stronger version of the strict convexity assumption on $q(\theta)$ used in Theorem~\ref{thm:con} and is satisfied when $q(\theta)$ is strongly convex on $\Theta$, {which is a standard condition frequently adopted in the literature for analyzing the convergence rates of gradient descent algorithms \citep[e.g.,][]{GhLan12,BoCuNo18}.}

The following lemma provides a good estimate for the weighted sum of a sequence of decreasing functions of order $O(1/k^s)$, $s\in(0,1)$. The result will be repeatedly used in the subsequent analysis.
\begin{lemma}\label{lem:order}
Let $u(i)={a}/{i^p}$ and $w(i)=O(1/i^s)$, where $a>0$, $p,\,s\in(0,1)$, and $w(i)>0$ for all $i=1,2,\ldots$. Then
$$\sum_{i=1}^k\Big[ \prod_{j=i+1}^k \big(1-u(j) \big)\Big]u(i)w(i)=O(k^{-s}).$$
\end{lemma}
\proof{Proof.}
See Section~\ref{appendix_order} of the Appendix.
\endproof

A characterization of the convergence rate for the mean squared errors of
the quantile estimates is given below; see Appendix~\ref{appendix_qkrate} for a proof.
\begin{lemma}\label{lem:qkrate}
Assume {A1}$-$A3 and B1 hold. Then the sequence $\{q_k\}$ generated by (\ref{qk}) satisfies
$$\sqrt{E\big[\big(q_k-q(\theta_k)\big)^2  \big]}=O\Big(\frac{\alpha_k}{\gamma_k} \Big)+O(\gamma_k^{\frac{1}{2}}).$$
\end{lemma}

From the discussion in Section~\ref{sec2}, the $q_k$ iteration (\ref{qk}) is an SA method for solving a sequence of time-varying root-finding problems. The $O({\alpha_k}/{\gamma_k})$ term above reflects the influence of the slowest component $\theta_k$ on the tracking ability of the $\{q_k\}$ sequence.
In particular, a large $\alpha_k$ value implies that the underlying input parameter $\theta_k$ will change quickly over the iterations, in which case
the step size $\gamma_k$ should decay sufficiently slowly to ensure that the $\{q_k\}$ sequence could properly follow the true quantiles $\{q(\theta_k)\}$.
When $\theta_k$ is fixed, i.e., $\alpha_k=0$, it can been seen from the proof of Lemma~\ref{lem:qkrate} that the rate of convergence of $\{q_k\}$ to the true quantile in MAE is of order $O(\gamma_k^{1/2})$. Therefore, if (\ref{qk}) is used as a stand-alone procedure for estimating distribution quantiles, its best rate of convergence is $O(1/\sqrt{k})$ (e.g., when $\gamma_k=r/k$). This is consistent but stronger than the classical (weak) convergence rate result for root-finding SA algorithms.

Based on Lemma~\ref{lem:qkrate}, we further obtain the following convergence rate result for the gradient estimates $\{D_k\}$.
\begin{lemma}\label{lem:dkrate}
{Assume A1$-$A3, B1, and B2 hold.} Then the $\{D_k\}$ sequence satisfies
\begin{align*}
\sqrt{E[\|D_k-\nabla_{\theta}q(\theta)|_{\theta=\theta_k} \|^2 ]}=O\Big(\frac{\alpha_k}{\gamma_k} \Big)+{O(c_k^2)}+O\Big(\frac{\beta_k^{\frac{1}{2}}}{c_k} \Big).
\end{align*}
\end{lemma}
\proof{Proof.}
See Section~\ref{appendix_dkrate} of the Appendix. \Halmos
\endproof
The result also has an intuitive explanation. A comparison of the results of Lemmas~\ref{lem:qkrate} and \ref{lem:dkrate} indicates that the $O\big({\alpha_k}/{\gamma_k} \big)$ term is attributed to the approximation error of the quantile estimator. Lemma~\ref{lem:dkrate} shows that this error sets a limit on the convergence speed of $\{D_k\}$, suggesting that performance of the two coupled SA recursions is primarily governed by the rate of the slower component.
In the special case when the input parameter vector is fixed, the $O\big({\alpha_k}/{\gamma_k}\big)$ term vanishes. Consequently, if (\ref{qk}) and (\ref{dk}) are jointly used as a means for quantile sensitivity analysis, then the lemma implies that the convergence rate of the gradient estimates is $O(c_k^2)+O\big({\beta_k^{{1}/{2}}}/{c_k}\big)$ provided that $\gamma_k=o(\beta_k)$ (due to A3(d)). If in addition the quantile is also fixed, then (\ref{dk}) alone is in the form of the standard SPSA algorithm, and the rate result is simply given by $O(c_k^2)+O\big({\beta_k^{1/2}}/{c_k}\big)$. When $\beta\leq 6\tau$, this further reduces to $O\big({\beta_k^{1/2}}/{c_k}\big)$, which becomes identical to the rate result previously obtained in \cite{spall92}, except in the mode of convergence.

Finally, we arrive at the following main convergence rate result for SPQO.
\begin{theorem}\label{thm:rate}
{Assume conditions A1-A3, B1-B3 hold,} then the sequence $\{\theta_k \}$ generated by SPQO satisfies
\begin{equation}\label{rate}
E\big[\|\theta_k-\theta^*\| \big]=O\Big(\frac{\alpha_k}{\gamma_k} \Big)+{O(c_k^2)}+O\Big(\frac{\beta_k^{\frac{1}{2}}}{c_k} \Big).
\end{equation}
\end{theorem}
\proof{Proof.}
{Define $\psi_k:=\theta_k-\theta^*$. We have from (\ref{projection}) that
\begin{align*}
\psi_{k+1}&=\psi_k-\alpha_k(D_k-\nabla_{\theta}q(\theta)|_{\theta=\theta_k})-\alpha_k\nabla_{\theta}q(\theta)|_{\theta=\theta_k}+\alpha_kZ_k\\
&=\psi_k-\alpha_k\eta_k-\alpha_k\nabla_{\theta}q(\theta)|_{\theta=\theta_k}+\alpha_kZ_k,
\end{align*}
where $\eta_k=D_k-\nabla_{\theta}q(\theta)|_{\theta=\theta_k}$. Since $\nabla_{\theta}q(\theta)|_{\theta=\theta^*}=0$, a Taylor series expansion of $\nabla_{\theta}q(\theta)|_{\theta=\theta_k}$ around $\theta^*$ shows that
\begin{align*}
\psi_{k+1}&= \psi_k-\alpha_k\eta_k-\alpha_kH(\bar\theta_k)\psi_k+\alpha_kZ_k\\
&=\big(I-\alpha_kH(\bar\theta_k)\big)\psi_k-\alpha_k\eta_k+\alpha_kZ_k,
\end{align*}
where $\bar \theta_k$ is on the line segment between $\theta_k$ and $\theta^*$. Taking norms at both side, using the Rayleigh-Ritz inequality \citep[cf., e.g.,][]{rugh96} and condition B3, we obtain that for all $k$ sufficiently large such that $\alpha_k\varrho<1$,
\begin{align*}
\|\psi_{k+1}\|&\leq \|\big(I-\alpha_kH(\bar\theta_k)\big)\psi_k\|+\alpha_k\|\eta_k \|+\alpha_k\|Z_k\|\\
&\leq (1-\alpha_k\varrho)\|\psi_k \|+\alpha_k\|\eta_k\|+\alpha_k\|Z_k\|.
\end{align*}
It follows that
\begin{equation}\label{tmp8}
E[\|\psi_{k+1}\|]\leq (1-\alpha_k\varrho)E[\|\psi_k \|]+\alpha_kE[\|\eta_k\|]+\alpha_kE[\|Z_k\|].
\end{equation}}

We now derive a bound for $E[\|Z_k \|]$. Since $\theta^*$ is in the interior of $\Theta$, there is a constant $\varsigma>0$ such that the $2\varsigma$-neighborhood of $\theta^*$ is contained $\Theta$. Let $\mathcal{E}_k=\{\|\theta_{k+1}-\theta^* \|\geq 2\varsigma \}$. Note that
by (\ref{cone}), the occurrence of $\mathcal{E}_k^c$ implies that $Z_k=0$.
Using this observation, we obtain that
\begin{align}\label{tmp9}
E[\|Z_k \|]&=E[\|Z_k\| |\mathcal{E}_k]P(\mathcal{E}_k)+E[\|Z_k\| |\mathcal{E}^c_k]P(\mathcal{E}^c_k) \nonumber \\
&\leq E[\|D_k \|]P(\|\theta_{k+1}-\theta^* \|\geq 2\varsigma) \nonumber \\
&\leq E[\|D_k \|]P\big(\|\theta_{k+1}-\theta_k \|\geq \varsigma \cup \|\theta_k-\theta^* \|\geq \varsigma \big) \nonumber \\
&\leq E[\|D_k \|]\frac{E[\|\theta_{k+1}-\theta_k \| ]}{\varsigma}+E[\|D_k \|]\frac{E[ \|\psi_k\|]}{\varsigma} \mbox{~~~~by Markov's inequality}  \nonumber \\
&\leq \frac{2\alpha_kE^2[\|D_k \| ]}{\varsigma}+E[\|D_k \|]\frac{E[ \|\psi_k\|]}{\varsigma},
\end{align}
where the first inequality is due to the fact that $\|Z_k\|\leq \|D_k\|$ \citep[see, the proof of Lemma 5 in][]{huetal2021}, and the last step follows from $\|\theta_{k+1}-\theta_k\|\leq \alpha_k\|D_k-Z_k\|\leq 2\alpha_k\|D_k\|$.

{Next, substitute the bound (\ref{tmp9}) into (\ref{tmp8}) and combine like terms,
\begin{equation*}
E[\|\psi_{k+1} \|]\leq \Big(1-\alpha_k\big(\varrho-\frac{E[\|D_k \| ]}{\varsigma} \big) \Big)E[\|\psi_k \| ] +\alpha_kE[\|\eta_k \| ]+\frac{2\alpha_k^2E^2[\|D_k \|]}{\varsigma}.
\end{equation*}
Note that since $\theta_k\rightarrow \theta^*$ w.p.1 and $\nabla_{\theta}q(\theta)|_{\theta=\theta^*}=0$, the continuity of $\nabla_{\theta}q(\theta)$ indicates that $\|\nabla_{\theta}q(\theta)|_{\theta=\theta_k} \|\rightarrow 0$ w.p.1. This, together with the boundedness of $\nabla_{\theta}q(\theta)$ (due to the compactness of $\Theta$), shows that $E[\|\nabla_{\theta}q(\theta)|_{\theta=\theta_k} \|]\rightarrow 0$ by the dominated convergence theorem. Thus, we have from Lemma~\ref{lem:dkrate} that $E[\|D_k \|]\leq E[\|D_k-\nabla_{\theta}q(\theta)|_{\theta=\theta_k} \|]+E[\|\nabla_{\theta}q(\theta)|_{\theta=\theta_k} \|]\rightarrow 0$. This then implies the existence of an integer $\mathcal{N}>0$ such that $\varrho-E[\|D_k \| ]/\varsigma\geq \varrho/2:=\bar\varrho$ for all $k\geq \mathcal{N}$. Consequently we get that for all $k\geq \mathcal{N}$
\begin{align*}
E[\|\psi_{k+1} \| ]&\leq (1-\alpha_k\bar\varrho)E[\|\psi_k \| ]+\alpha_kE[\|\eta_k \| ]+\frac{2\alpha_k^2E^2[ \|D_k\|]}{\varsigma}.
\end{align*}
Directly expand the above inequality from term $\mathcal{N}$ onwards,
\begin{align*}
E[\|\psi_k \| ]&\leq \prod_{i=\mathcal{N}}^k(1-\alpha_i\bar\varrho)E[\|\psi_{\mathcal{N}} \|]+\sum_{i=\mathcal{N}}^k\Big[\prod_{j=i+1}^k(1-\alpha_j\bar\varrho) \Big]\alpha_iE[\|\eta_i \| ]\\
&\hspace{1cm}+\sum_{i=\mathcal{N}}^k\Big[\prod_{j=i+1}^k(1-\alpha_j\bar\varrho) \Big]\alpha_i\frac{2\alpha_iE^2[\|D_i \| ]}{\varsigma}.
\end{align*}
Finally, because $\prod_{i=\mathcal{N}}^k(1-\alpha_i\bar\varrho)=o(k^{-1})$, $E[\|\eta_k\|]=O(\alpha_k/\gamma_k)+O(c_k^2)+O({\beta_k^{1/2}}/{c_k})$ (Lemma~\ref{lem:dkrate}), and $\alpha_kE^2[\|D_k\|]=O(\alpha_k)$ (Lemma~\ref{lem:dk}(i)), a direct application of  Lemma~\ref{lem:order} leads to the conclusion
$E[\|\psi_k \|]=O(\alpha_k/\gamma_k)+O(c_k^2)+O({\beta_k^{1/2}}/{c_k})$.}
\Halmos
\endproof
In view of Lemma~\ref{lem:dkrate}, regardless of the choice of $\alpha_k$, the MAE of $D_k$ converges at a rate that is always slower than $\alpha_k$ itself, which results in the long run behavior of (\ref{thetak}) being dominated by the errors in gradient estimation. Therefore, in contrast to single-timescale SA, whose rate is determined by its step-sizes, an interesting observation from Theorem~\ref{thm:rate} is that the step size $\alpha_k$ in (\ref{thetak}) does not have a direct effect on the convergence rate of $\{\theta_k\}$, but only does so indirectly through the expression for the convergence rate of the faster component $D_k$ in Lemma~\ref{lem:dkrate}.

{From A3(a), both $\beta_k$ and $c_k$ should be chosen to satisfy $\tau+1/2<\beta\leq 1$. In addition, to improve the rate in (\ref{rate}), it is clear that $\alpha$ should be taken large, which gives the obvious choice $\alpha\approx 1$ (note that $\alpha\in(0,1)$). On the other hand, since  $\beta_k^{1/2}/c_k=O(k^{-(\beta/2-\tau)})$ and $c_k^2=O(k^{-2\tau})$, the $O(c_k^2)+O(\beta_k^{1/2}/c_k)$ term is optimized when $\beta=6\tau$, yielding $O(c_k^2)=O(\beta_k^{1/2}/c_k)=O(k^{-2\tau})$. Consequently, by equating the terms in (\ref{rate}), we find that under the above choice of $\alpha,\,\beta$ and the constraint $\alpha>\gamma>\beta$ (A3(d)), an upper bound on the MAEs of $\{\theta_k\}$ diminishes at an optimal rate that can be made arbitrarily close to $O(k^{-1/4})$, which is approximately attained when $\alpha\approx 1$, $\gamma=3/4$, $\beta\approx3/4$ and $\tau=1/8$.}

Finally, for the sake of completeness, we conclude by stating the following rate result we have for SDQO, implying that the algorithm essentially shares the same $O(k^{-1/4})$ {best convergence rate bound} as SPQO, except possibly in the constant contained in the big-O notation.
\begin{theorem}\label{thm:ratesdqo}
Assume that the conditions of Theorem~\ref{thm:sdqo} and B1-B3 hold but with $(\hat q_k,\hat \theta_k)$ replacing $(q_k,\theta_k)$ in B1. Then the MAEs of the sequence $\{\hat \theta_k \}$ generated by SDQO satisfies
\begin{equation*}
E\big[\|\hat \theta_k-\theta^*\| \big]=O\Big(\frac{\alpha_k}{\gamma_k} \Big)+{O(c_k^2)}+O\Big(\frac{\beta_k^{\frac{1}{2}}}{c_k} \Big).
\end{equation*}
\end{theorem}

\section{Simulation Experiments}\label{sec5}

To illustrate the algorithms, we perform some computational experiments on a set of artificially created black-box test functions and a queueing example. In all cases, the performance of SPQO and SDQO is compared with that of the QG algorithm and the surrogate-based gTSSO-QML algorithm proposed in \cite{ng2021}. We describe the latter two algorithms in more detail now.
As discussed in Section~\ref{sec1}, QG is a single-timescale SA algorithm that employs the conventional SD (\ref{qFD}) to approximate quantile gradients, where the true quantiles are estimated by order statistics. Denote by $\upsilon_k>0$ the perturbation size, and let $\widehat q(\theta)$ stand for the $\lceil n_k\varphi \rceil$th order statistic of an output sample $Y_1,\ldots, Y_{n_k}\sim F(\cdot;\theta)$ of size $n_k$. The QG algorithm uses the following update:
\begin{equation}\label{qg}
\theta_{k+1}=\Pi_{\Theta}\big(\theta_k-\rho_k \tilde D_k \big),
\end{equation}
where $\rho_k>0$ is the step size and $\tilde D_k$ is the gradient estimate whose $i$th element is given by
$$\tilde{D}_{k,i}=\frac{\widehat q\big((\tilde \theta_{k,1},\ldots,\theta_{k,i}+\upsilon_k,\ldots,\tilde \theta_{k,d})^T\big)-\widehat q\big((\tilde \theta_{k,1},\ldots,\theta_{k,i}-\upsilon_k,\ldots,\tilde \theta_{k,d})^T\big)}{2\upsilon_k},~i=1,\ldots,d.$$
Here, $\theta_{k,i}$ is the $i$th element of $\theta_k$ and $\tilde \theta_{k,j}$'s are random variables uniformly distributed over $[\theta_{k,j}-\upsilon_k,\theta_{k,j}+\upsilon_k ]$ for all $j=1,\ldots,d$, $j\neq i$.
Note that this construction differs from our proposed algorithms in that $2dn_k$ function evaluations (as opposed to three in SPQO and $2d+1$ in SDQO) are needed at each iteration. 
Our implementation of QG is based on the following parameter values: $\rho_k=1/k$, $\upsilon_k=1/k^{0.501}$, and $n_k=\lceil k^{2.003}\rceil$, which are the minimum required to satisfy the conditions for the convergence of the algorithm \citep[see Theorem 8 of][]{kibzun12}.
The gTSSO-QML algorithm uses a stochastic cokriging model to approximate the response surfaces of a set of quantile functions with progressively increasing quantile levels and selects new design points by optimizing an expected improvement criterion.
Following \cite{ng2021}, the algorithm parameter values are determined in our implementation based on a cross-validation test, and at each step, simulation samples are adaptively allocated to the selected design points by using the optimal computing budget allocation method.

Similar to many other Bayesian optimization approaches, gTSSO-QML is very computationally demanding on high-dimensional problems, so we have implemented the algorithm on a parallel computing platform with 164 nodes. Each node has two Intel Xeon E5-2683v3 processors (with 14 cores each running at 2.0 GHz) and 128 GB of memory. The computational experiments for all other algorithms are performed on a Window PC with an Intel Core i5 1.8GHz processor and 8 GB of memory.

\subsection{Black-box Test Functions}\label{sec51}
Six noisy black-box functions are tested, with dimensions varying from 2 to 20.
Case 1 contains multiplicative noise but only two decision variables and is relatively easy to solve.
In case 2, the noise also scales the function, but as the problem dimension increases, the distribution of the function may become extremely flat, making its extreme quantiles challenging to estimate.
In case 3, on the other hand, the noise is additive, and under the optimal parameter values, the quantiles are very distant from the origin, so predicting their values could also be challenging, especially when the initial estimates are far from the true values.
In case 4, the noise is both multiplicative and additive. {Cases 5 and 6 are relative low-dimensional multi-modal problems, and each contains a large number of local optima.}
\begin{enumerate}[leftmargin=*,labelindent=16pt,label=\bfseries Step \arabic*.]
\item[Case 1.] $Y(\theta)=\big(2.6(\theta_1^2+\theta_2^2)-4.8\theta_1\theta_2 \big)X+10$, where  $\Theta=[-2,2]^2$.
\item[Case 2.] $Y(\theta)=\big(\sum_{i=1}^d (\theta_i -i)^2+1 \big)X$, where $d=10$ and
$\theta_i\in [i-1,i+1]$ for $i=1,\ldots,d$.
\item[Case 3.] $Y(\theta)=X+\sum_{i=1}^d(\theta_i-i)\theta_i$, where $d=20$ and $\Theta=[-20,20]^d$.
\item[Case 4.] $Y(\theta)=\frac{1}{d}\sum_{i=1}^d (\theta_i-1)^2 X+\frac{1}{d}\sum_{i=1}^d(\theta_i^4-16\theta_i^2+5\theta_i),$ where $d=20$ and $\Theta=[1,4]^d$.
{\item[Case 5.] $Y(\theta)=\Big[-10\exp\Big(-0.2\sqrt{\frac{1}{d}\sum_{i=1}^d \theta_i^2}\Big)-\exp(\frac{1}{d}\sum_{i=1}^d \cos(\pi \theta_i))+11+e\Big]X$, where $d=5$ and $ \Theta=[-5,5]^d$.
\item[Case 6.] $Y(\theta)=\frac{1}{d}\sum_{i=1}^d\left[0.4\sin^2(0.2\pi(\theta_i-0.9))+0.3\sin^2(0.4\pi(\theta_i-0.9))+0.001(\theta_i-0.9)^2\right]+X$, where $d=5$ and $ \Theta=[-10,10]^d$.} 
\end{enumerate}
For each test problem, we consider two quantile levels: $\varphi=0.6$, $\varphi=0.95$ and two choices of the (unknown) noise distribution: $X\sim N(0,1)$ and $X\sim Cauchy(0,1)$, resulting in {24} total test scenarios. {The optimal quantile values in all scenarios are listed in Table~\ref{tab1}.} Note that under the Cauchy noise, neither the mean nor the variance of the output distribution exists.

\begin{table}[!t]
\caption{Optimal quantile values in the 24 test scenarios.}
  \label{tab1}
 \centering
\begin{tabular}{c c c c c c}
    \hline
    \multirow{2}{*}{Case} & \multicolumn{2}{c}{Normal} & & \multicolumn{2}{c}{Cauchy} \\
        & $\varphi=0.6$ & $\varphi=0.95$ & & $\varphi=0.6$ & $\varphi=0.95$\\
    \hline
     1 & 10 & 10 & & 10 &10\\
     2 & 0.25  & 1.64 & & 0.32 & 6.31\\
     3 & -717.25  &-715.86& & -717.18 &-711.19\\
     4 & -49.29 &-45.32& &-49.08&-34.62\\
     5 & 0.25 &1.64& &0.32 & 6.31 \\
     6 & 0.25 &1.64& & 0.32 & 6.31 \\
    \hline
\end{tabular}
\end{table}

In the implementation of SPQO and SDQO, the decay rates of the parameters are determined from the result of Section~\ref{sec4}, i.e., {$\alpha=0.99$, $\gamma=0.75$, $\beta=0.74$, and $\tau=0.125$} (see the discussion at the end of that section).
Our experience indicates that their performance is not very sensitive to the choice of $\alpha_k$, in that the standard step-size $\alpha_k=2/k^{\alpha}$ seems to work well across a variety of test cases. The parameters $\beta_k$ and $c_k$ resemble those of SPSA, and we choose them to be of the forms $\beta_k=b/(k+R)^\beta$ and $c_k=c/(k+R)^\tau$ as suggested in \cite{spall03}, where $R$ is set to $10$ percent of the maximum number of iterations allowed, $b=\kappa_1(2R)^{\beta}$, and $c=\kappa_2(2R)^{\tau}$. This choice maintains the respective values of $\beta_k$ and $c_k$ to be
greater than $\kappa_1$ and $\kappa_2$ during the first $R$ iterations. The constants $\kappa_1$ and $\kappa_2$ are then selected by trial and error, and we find that values
satisfying $\kappa_2\in[0.1,0.9]$ and $0.005\leq \kappa_1/\kappa_2\leq 0.1$ all yield reasonable performance. Note that from (\ref{dk}) (resp. (\ref{vvhatdk})), the increment (if there is any) in each component of $D_k$ (resp. $\hat D_k$) is exactly $\beta_k/2\bar c_k$ (resp. $\beta_k/2\tilde c_k$). So the lower bound $0.005$ on $\kappa_1/\kappa_2$ prevents the updates in gradient estimates from becoming too small. Our numerical results reported here are based
on $\kappa_1=0.05$ and $\kappa_2=0.5$.
The choice of $\gamma_k$, on the other hand, is most critical to the performance. This is mainly because of the recursive procedure used for estimating $q_k$ (resp. $\hat q_k$). As can be observed from (\ref{qk}), each increment in quantile estimate is bounded in magnitude by $\gamma_k$. Thus, if a desired (true) quantile is far from the initial $q_0$, then a reasonable estimate of its value would take an enormously large number of iterations under the standard choice $\gamma_k=1/k^{\gamma}$, leading to excessively slow (finite-time) convergence behavior. One way to address this issue would be to take $\gamma_k$ to be of the same form as $\beta_k$, so that a non-negligible gain could be maintained in tracking the true quantile values. In our study, however, we simply take $\gamma_k=R/k^{\gamma}$. The intuitive reason is that
the large constant $R$ will provide enough impetus in early iterations to help the iterates move quickly towards the ``correct" quantile range, whose values can then be further fine-tuned as $\gamma_k$ decreases rapidly with $k$, due to the large decay rate $\gamma$.

\begin{table}[!t]
\caption{\baselineskip8pt Performance on test functions for normally distributed noise, based on 40 independent runs (standard error in parentheses).}
  \label{tab2}
 \centering
  \setlength{\tabcolsep}{1pt}
 {\small\begin{tabular}{c c c c c c c}
    \hline
    Case & SPQO & SPQO-CRN & SDQO & SDQO-CRN & QG & gTSSO-QML\\
    \hline
     & & & $\varphi=0.6$ & & & \\
     \hline
     1& 10.06 (8.0e-3) &  10.04 (7.8e-3)&10.06 (1.0e-2)  & 10.04 (7.0e-3) & 10.11 (1.3e-2) & {\bf10.01 (7.3e-4)}\\
     2&  0.30 (2.7e-3) & {\bf0.28 (1.0e-3)} & 0.55 (1.3e-2)  & 0.43 (1.6e-2) & 0.39 (1.4e-2) & 0.33 (1.7e-2) \\
     3&  -717.24 (3.3e-4) &  -717.24 (7.2e-4)  & -717.22 (1.4e-3) & {\bf-717.25 (1.9e-4)} & -717.23 (6.7e-4) &  -168.00 (2.71)\\
     4&  -49.22 (1.8e-3) & -49.19  (8.3e-4)  &  -49.16  (5.0e-3) &  {\bf -49.25  (3.7e-3)} &-49.17 (7.1e-3) &  -48.60 (3.9e-2)\\
     5& 1.13 (4.4e-2) & 1.05 (3.2e-2) & 1.28 (3.8e-2) & 1.18 (2.5e-2) &1.51 (3.7e-2) & {\bf0.62 (4.1e-2)}\\
     6& 0.50 (1.4e-2) & 0.52 (1.4e-2) & 0.57 (1.6e-2) & 0.54 (1.4e-2) & 0.54 (1.2e-2) & {\bf0.36 (5.1e-3)} \\
    \hline
    & & & $\varphi=0.95$ & & & \\
     \hline
    1 &  10.07 (6.6e-3) & 10.09 (1.0e-2)& 10.04 (5.3e-3) &  10.02 (4.0e-3) & 10.22 (2.9e-2) &{\bf10.00 (6.0e-6)} \\
    2 & 1.65 (3.3e-4)  &  {\bf1.64 (6.0e-6)} & 1.68 (2.3e-3) & 1.65 (4.6e-5) & 1.66 (9.5e-4) & 1.69 (3.0e-3)\\
    3 & -715.85 (5.6e-5)  &  {\bf-715.86 (5.1e-6)}  &   -715.82 (1.6e-3)&  -715.85 (1.7e-4)& -715.82 (1.8e-3)& -174.96 (3.59)\\
    4 & -45.22 (1.5e-3) &  -45.21 (6.8e-4)  & -45.13 (9.6e-3)  & {\bf -45.31 (9.5e-4)}&-45.05 (1.3e-2)&-44.41 (4.7e-2)\\
    5 & 4.85 (4.0e-1) & 5.31 (3.3e-1) & 6.43 (3.0e-1) & 5.17, (2.9e-1) &8.52 (1.9e-1) & {\bf4.10 (2.6e-1)}\\
    6 &  1.91 (1.4e-2) & 1.91 (1.6e-2) & 1.94 (1.4e-2)& 1.93 (1.3e-2) & 1.93 (1.5e-2) &{\bf1.76 (6.3e-3)}\\
    \hline
\end{tabular}}
\end{table}

\begin{table}[!t]
\caption{\baselineskip8pt Performance on test functions for Cauchy distributed noise, based on 40 independent runs (standard error in parentheses).}
  \label{tab3}
 \centering
  \setlength{\tabcolsep}{1pt}
 {\small
\begin{tabular}{c c c c c c c}

    \hline
    Case & SPQO & SPQO-CRN & SDQO & SDQO-CRN & QG & gTSSO-QML\\
    \hline
     & & & $\varphi=0.6$ & & & \\
    \hline
    1 & 10.06 (9.4e-3) & 10.03 (6.6e-3) & 10.07 (1.3e-2) & 10.03 (6.3e-3) &10.11 (1.3e-2) & {\bf10.00 (8.7e-4)}\\
    2 & 0.37 (3.2e-3)&  {\bf0.33 (4.8e-4)} &  0.55 (1.9e-2) &  0.36 (5.0e-3) & 0.47 (1.5e-2) & 0.37 (4.5e-3)\\
    3 & -717.16 (8.5e-4)&{\bf -717.17 (2.9e-5)}  &  -717.08 (5.5e-3) &   -717.17 (1.3e-4)& -717.15 (1.1e-3) & -165.34 (2.52)\\
    4 & -48.99 (2.9e-3) &  -48.98  (1.3e-3) &  -48.75  (1.8e-2)&  {\bf -49.03  (4.6e-3)} & -48.87 (1.2e-2) & -48.20 (4.4e-2)\\
    5 &1.64 (4.6e-2) &  1.34 (4.4e-2) & 1.76 (5.7e-2)  & 1.31, (5.5e-2)&1.97 (5.0e-2)&{\bf0.83 (5.8e-2)}\\
    6 & 0.60 (1.3e-2) & 0.56 (1.2e-2)  & 0.61 (1.6e-2) & 0.60 (1.4e-2) &0.63 (1.3e-2)& {\bf0.43 (5.6e-3)}\\
    \hline
     & & & $\varphi=0.95$ & & & \\
    \hline
    1 & 10.03 (6.4e-3) &   10.00 (5.3e-4)  &   10.04 (7.9e-3) &  10.01 (1.6e-3) & 120.03 (18.00) & \bf{10.00 (3.6e-5)}\\
    2 &6.48 (8.6e-3) &  {\bf6.31 (9.1e-5)} & 8.15 (3.5e-1) &  6.32 (1.4e-3) & 48.05 (1.25) & 6.67 (2.2e-2)\\
    3 & -711.17 (9.2e-4) &   {\bf-711.19 (2.2e-5)} &  -709.24 (1.9e-1) &  -710.60 (7.0e-2) & -703.94 (3.8e-1) & -160.79 (3.06)\\
    4 & -33.80 (3.9e-2) & {\bf-34.20 (1.0e-2)} &   -30.16 (1.9e-1) & -33.92 (1.0e-1) &  -19.54 (1.06) & -32.10 (1.1e-1)\\
    5 & 23.96 (2.99)& {\bf 8.21 (9.0e-1)}& 32.11 (2.99) & 16.69 (1.30) & 36.74 (1.30) & 14.64 (1.23)\\
    6 & 6.58 (1.6e-2)  & 6.58 (1.7e-2) & 6.64 (1.4e-2) & 6.61 (1.6e-2) &6.72 (1.2e-2)& {\bf6.43 (6.5e-3)}\\
    \hline
\end{tabular}}
\end{table}

\begin{figure}[!htbp]
	\centering
	\includegraphics[trim=1cm 7.5cm 1cm 6.9cm, clip,width=0.45\textwidth]{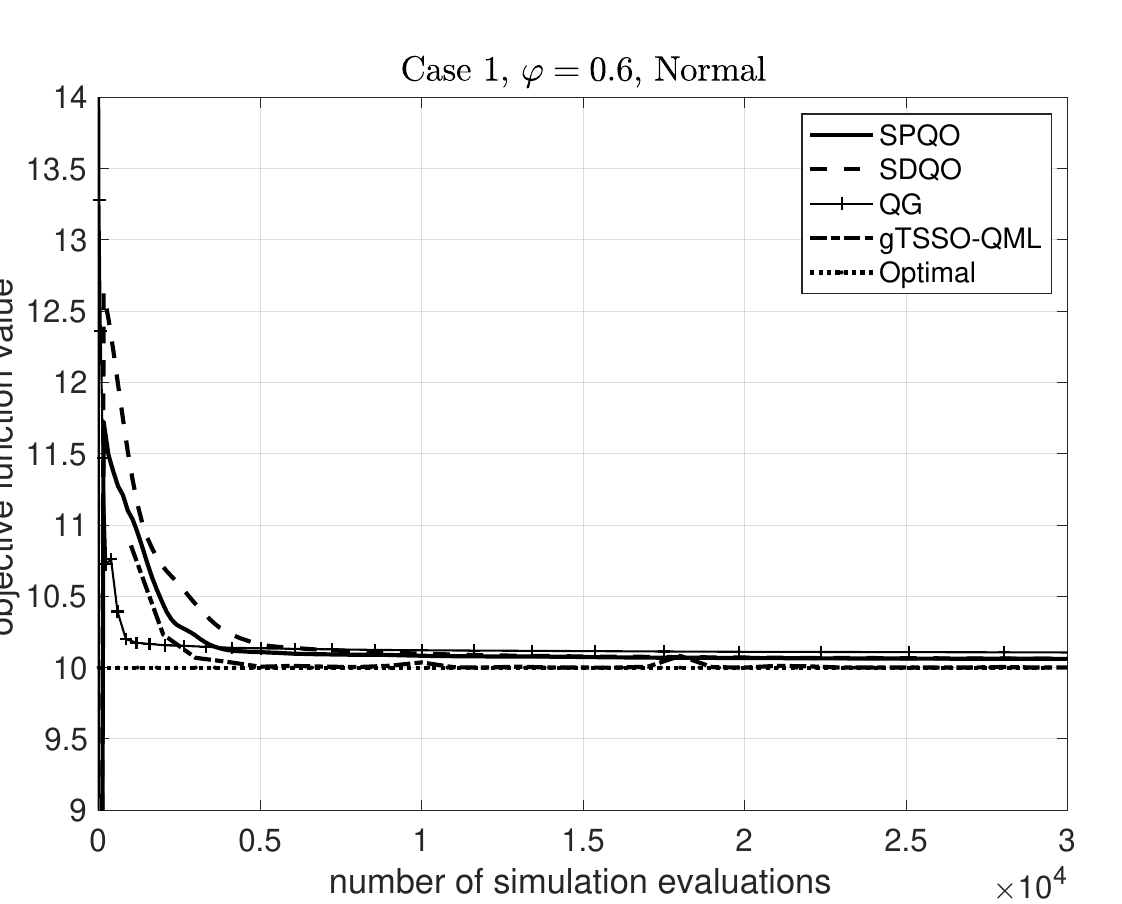}
	\includegraphics[trim=1cm 7.5cm 1cm 6.9cm, clip,width=0.45\textwidth]{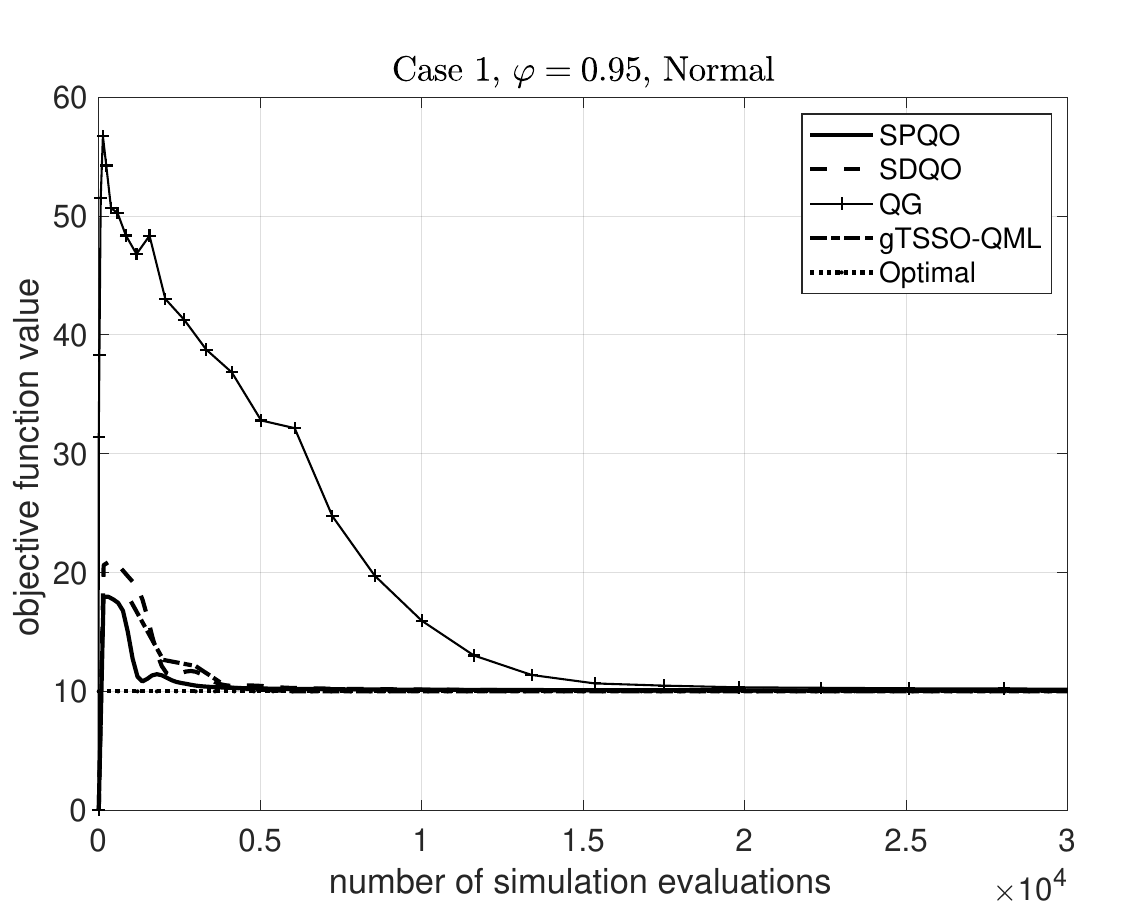}\\
	\vspace{1pt}
	\includegraphics[trim=1cm 7.5cm 1cm 6.9cm, clip,width=0.45\textwidth]{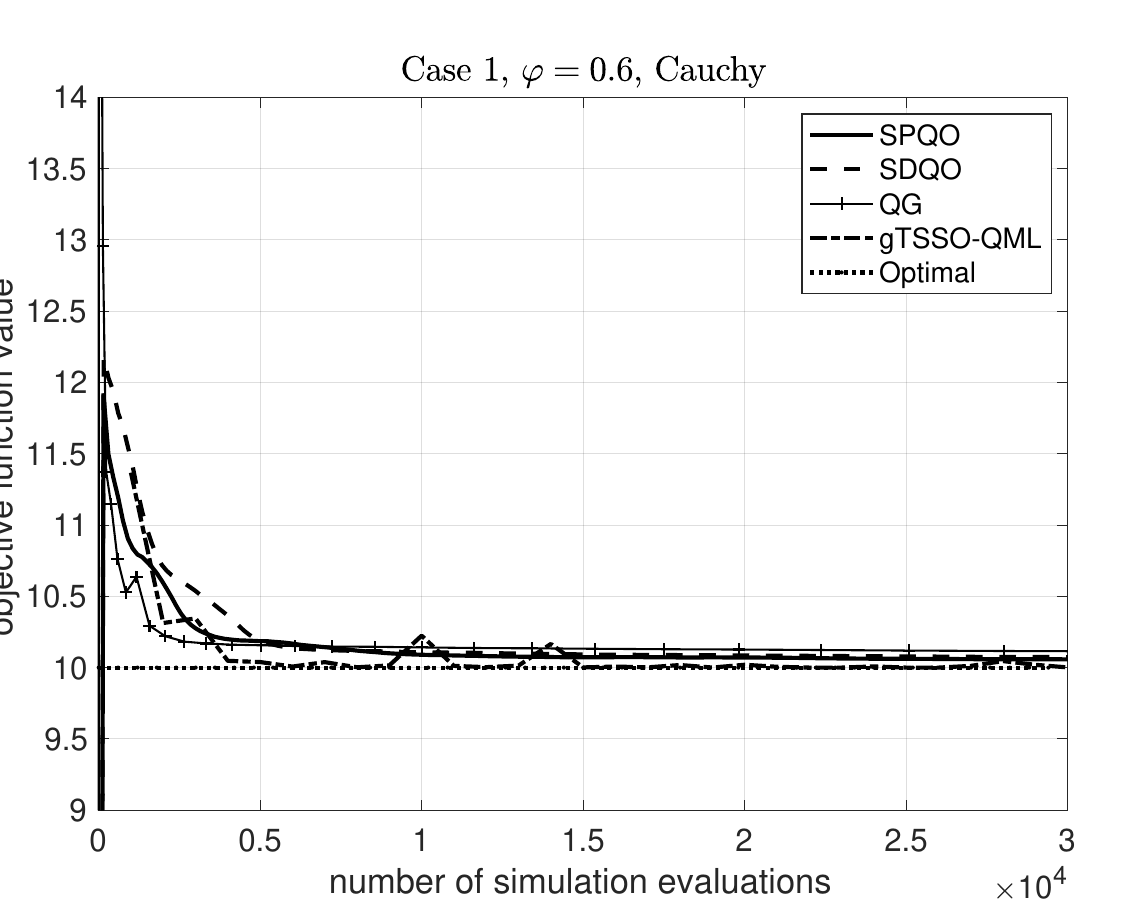}
	\includegraphics[trim=1cm 7.5cm 1cm 6.9cm, clip,width=0.45\textwidth]{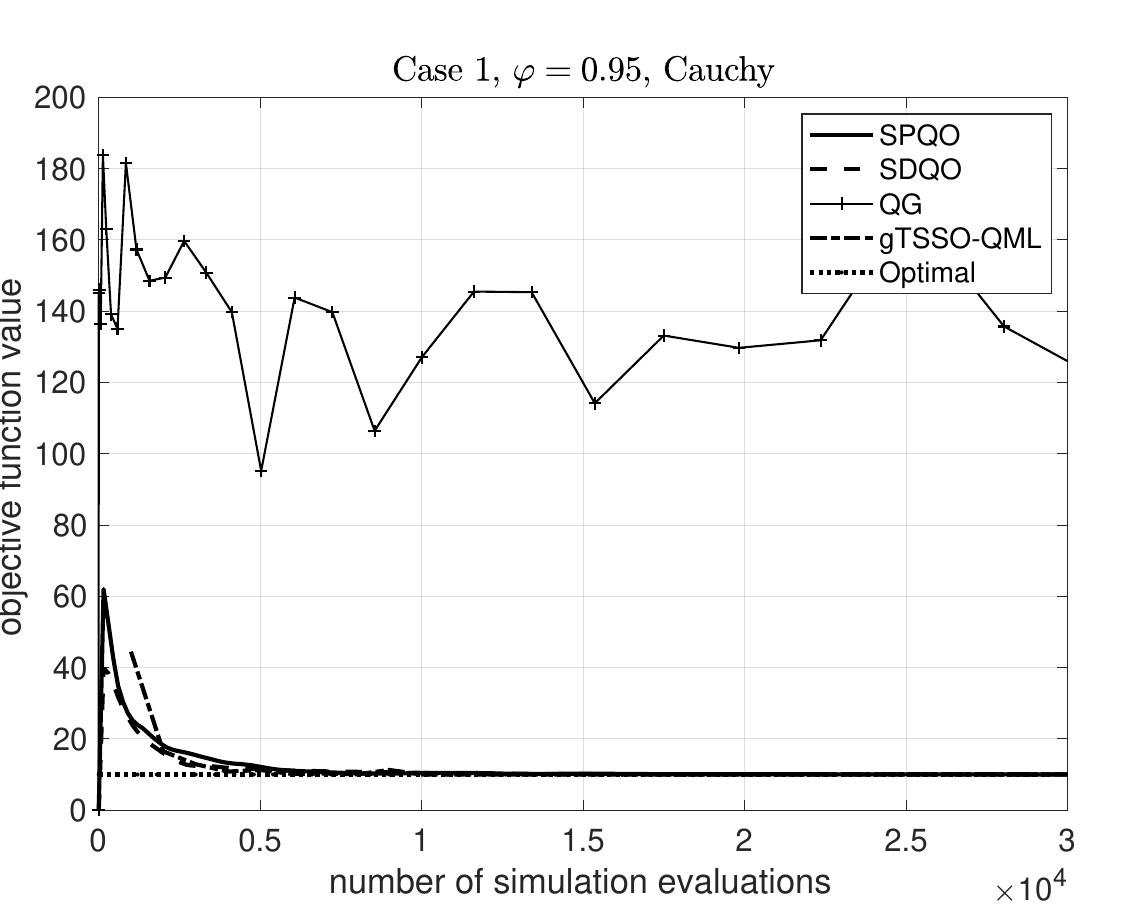}\\
	\vspace{1pt}
	\includegraphics[trim=1cm 7.5cm 1cm 6.9cm, clip,width=0.45\textwidth]{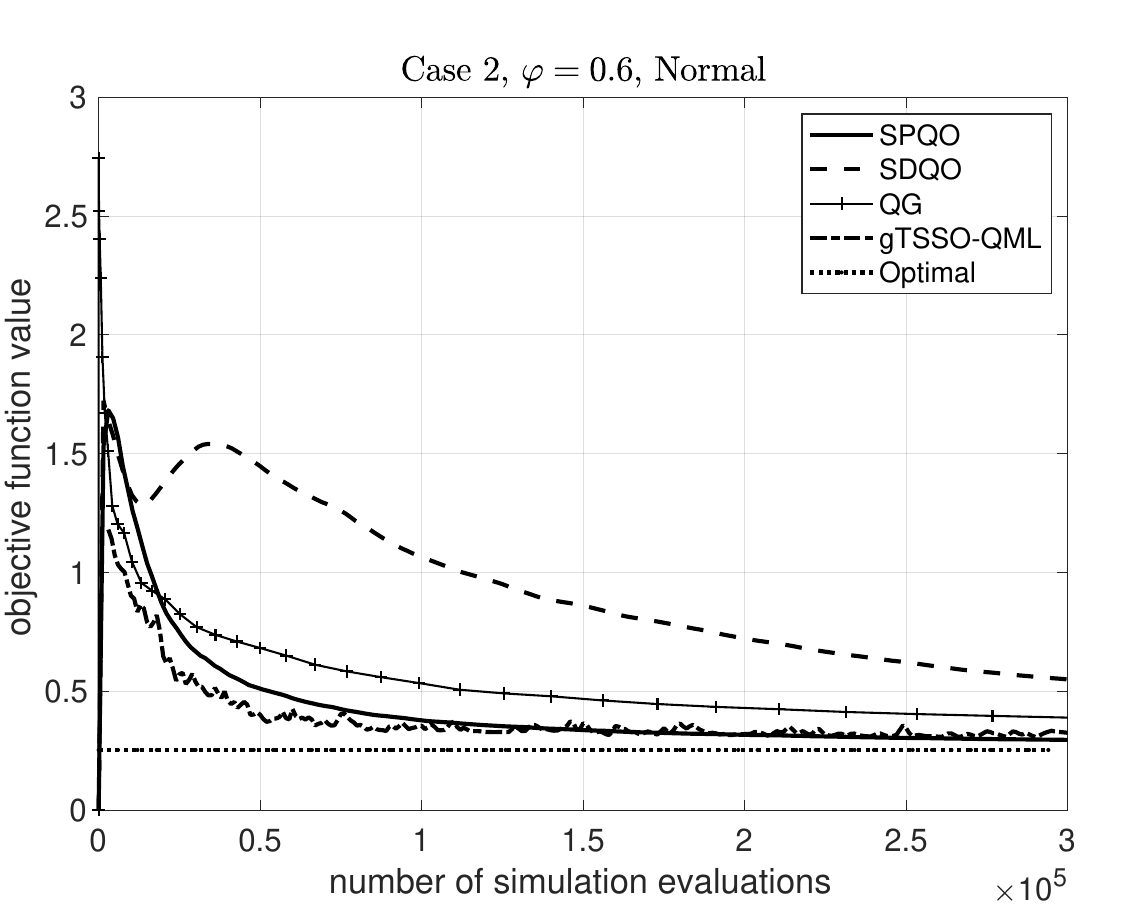}
	\includegraphics[trim=1cm 7.5cm 1cm 6.9cm, clip,width=0.45\textwidth]{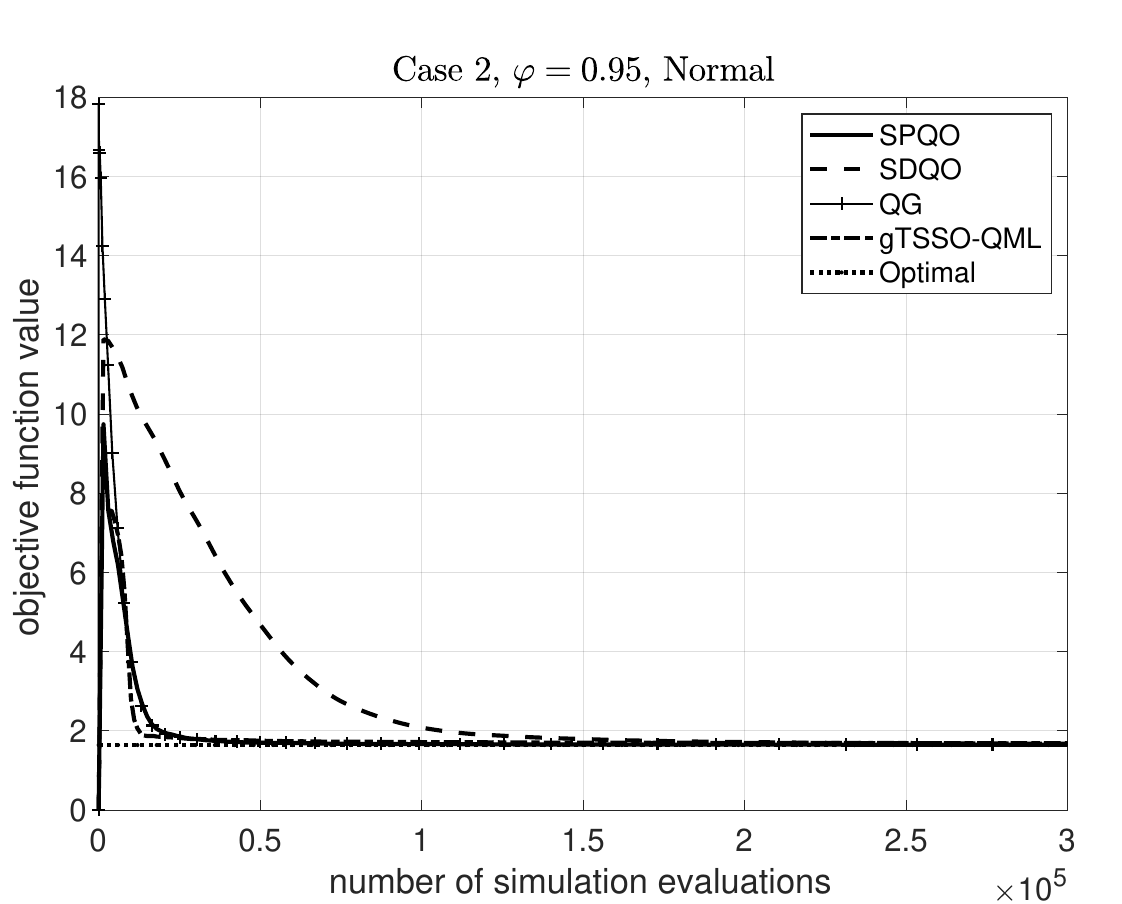}\\

   \includegraphics[trim=1cm 7.5cm 1cm 6.9cm, clip,width=0.45\textwidth]{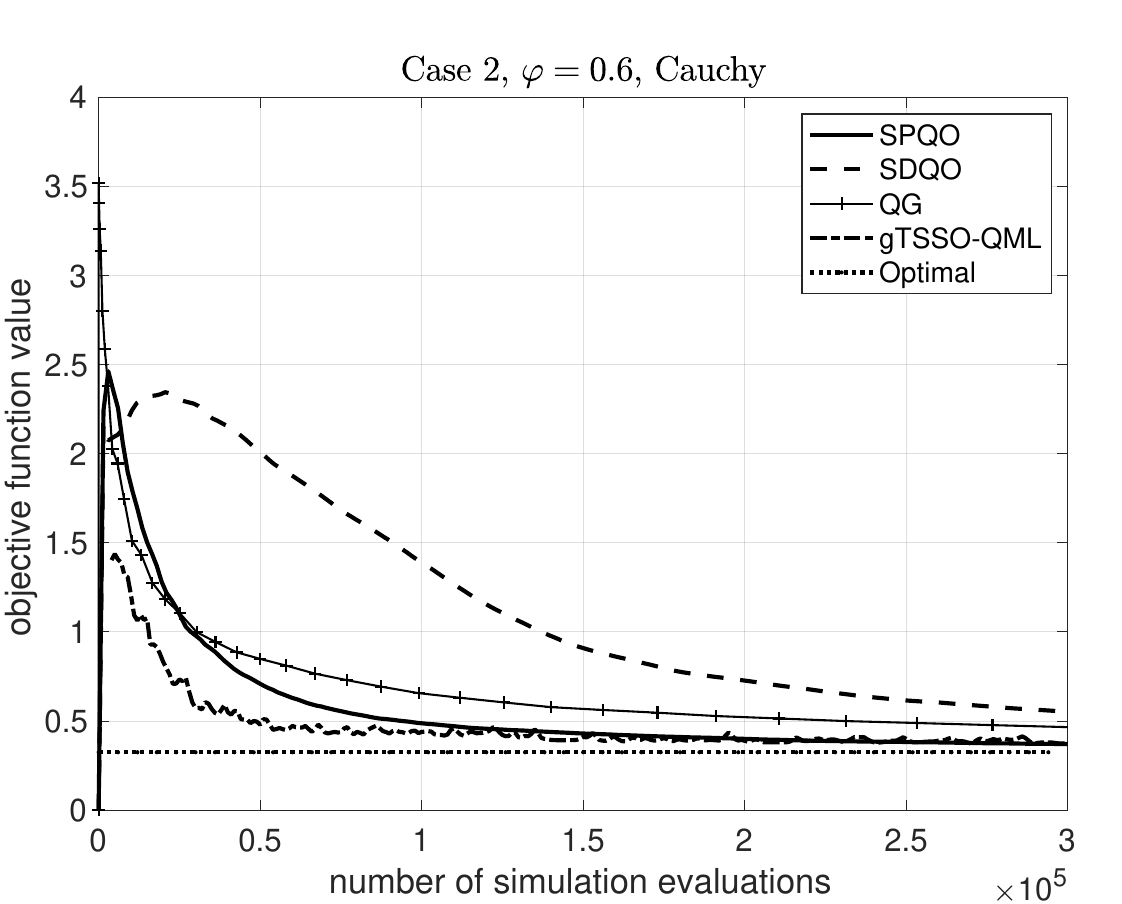}
	\includegraphics[trim=1cm 7.5cm 1cm 6.9cm, clip,width=0.45\textwidth]{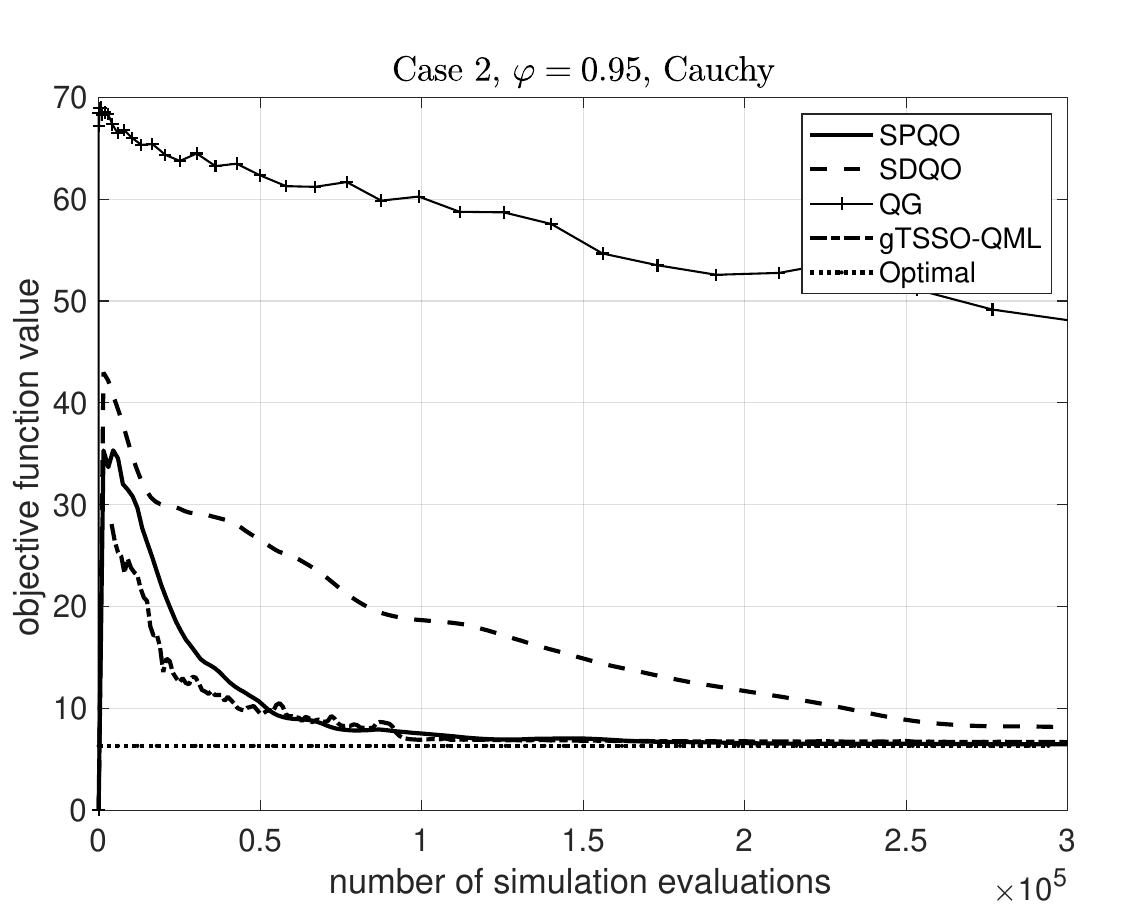}
	
	\caption{Performance of SPQO, SDQO, QG, and gTSSO-QML on test cases 1 and 2.}
	\label{fig1}
\end{figure}

\begin{figure}[!htbp]
	\centering
	\includegraphics[trim=1cm 7.5cm 1cm 6.9cm, clip,width=0.45\textwidth]{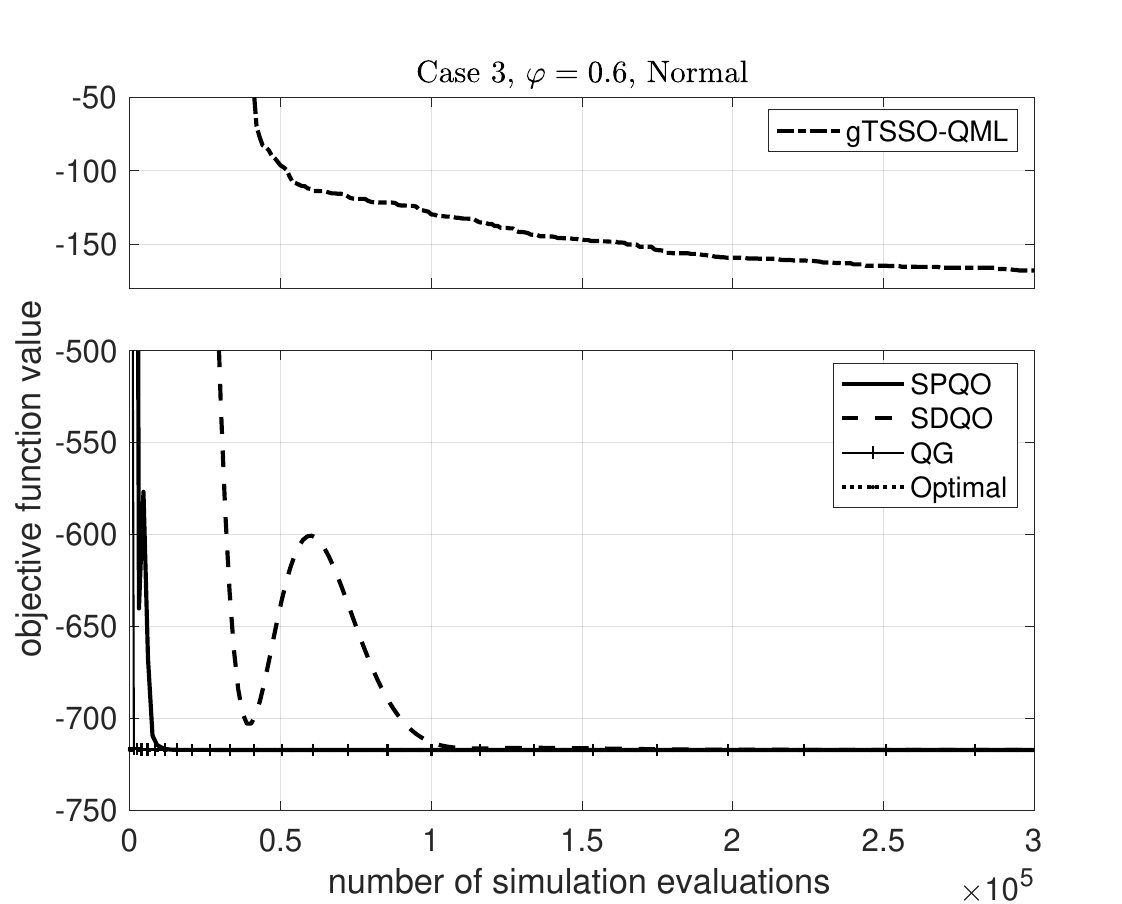}
	\includegraphics[trim=1cm 7.5cm 1cm 6.9cm, clip,width=0.45\textwidth]{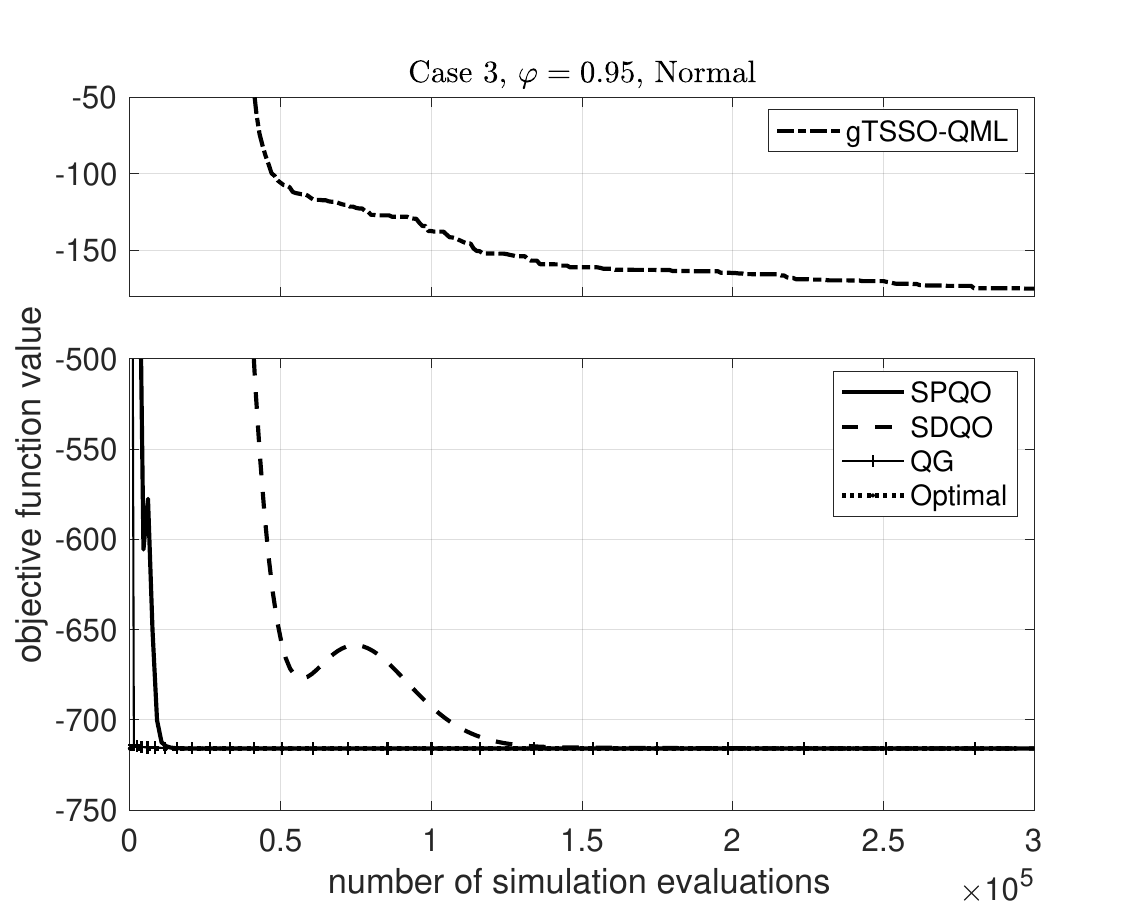}\\
	\vspace{1pt}
	\includegraphics[trim=1cm 7.5cm 1cm 6.9cm, clip,width=0.45\textwidth]{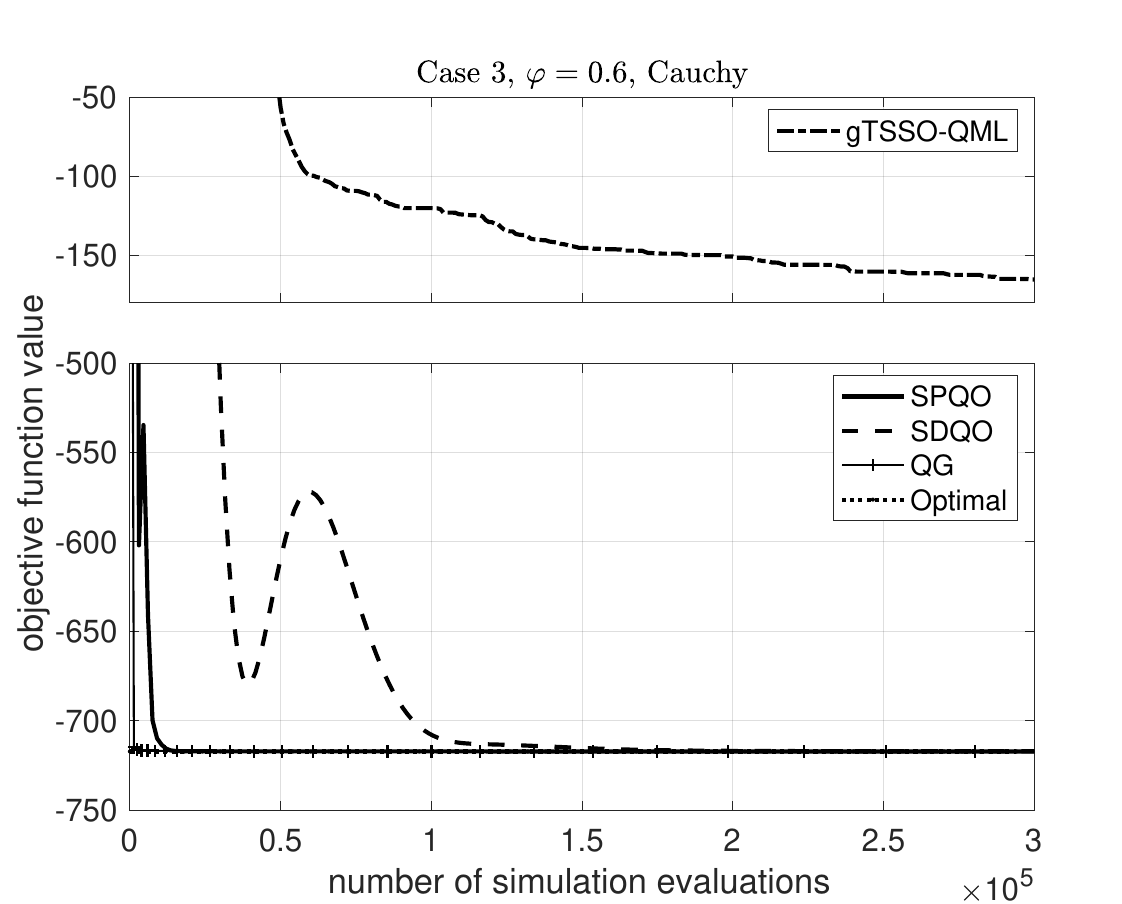}
	\includegraphics[trim=1cm 7.5cm 1cm 6.9cm, clip,width=0.45\textwidth]{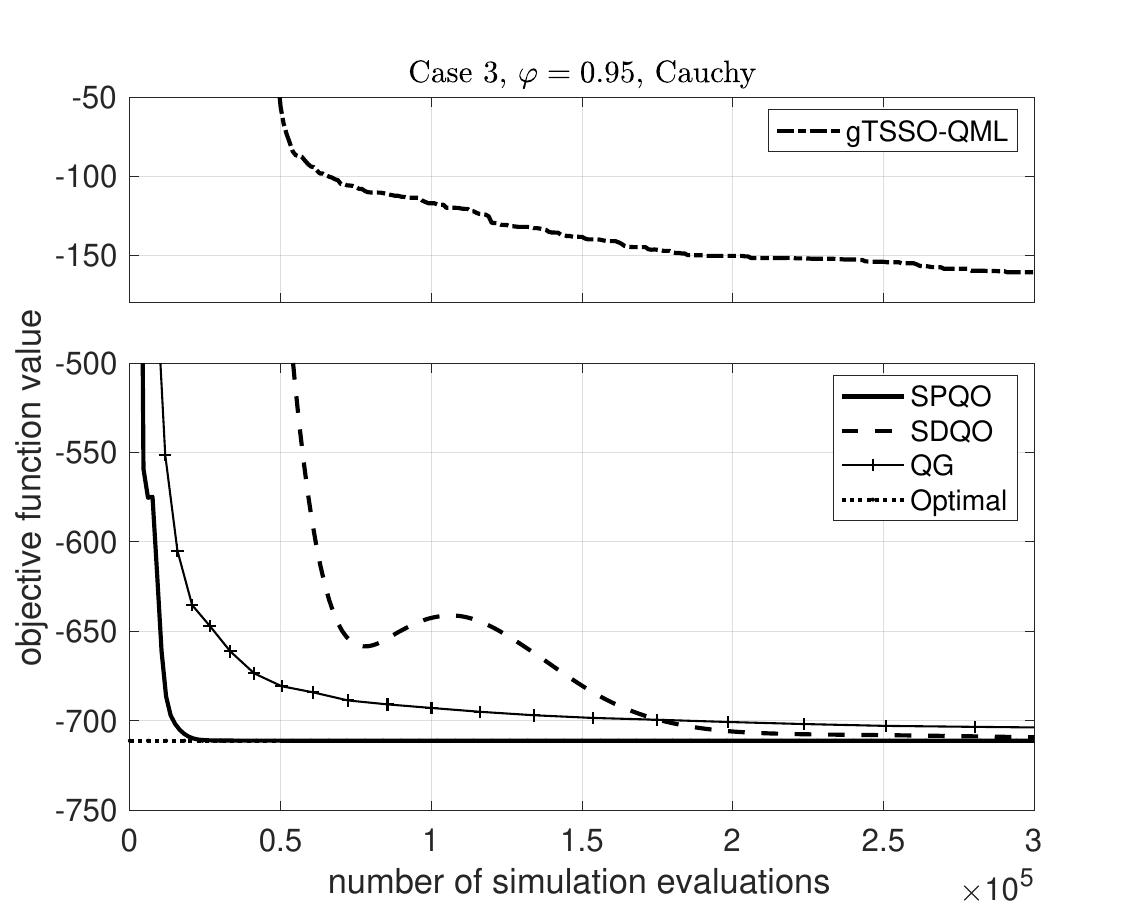}\\
	\vspace{1pt}
	\includegraphics[trim=1cm 7.5cm 1cm 6.9cm, clip,width=0.45\textwidth]{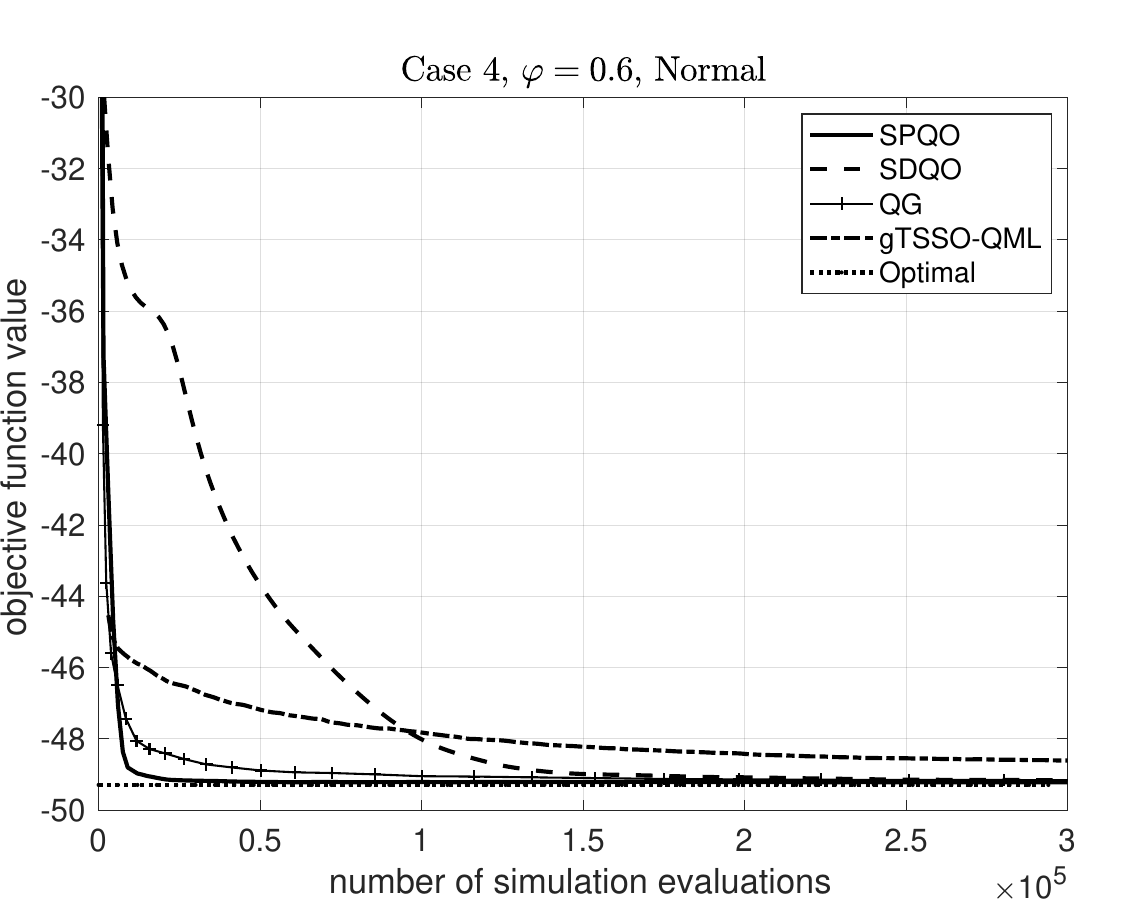}
	\includegraphics[trim=1cm 7.5cm 1cm 6.9cm, clip,width=0.45\textwidth]{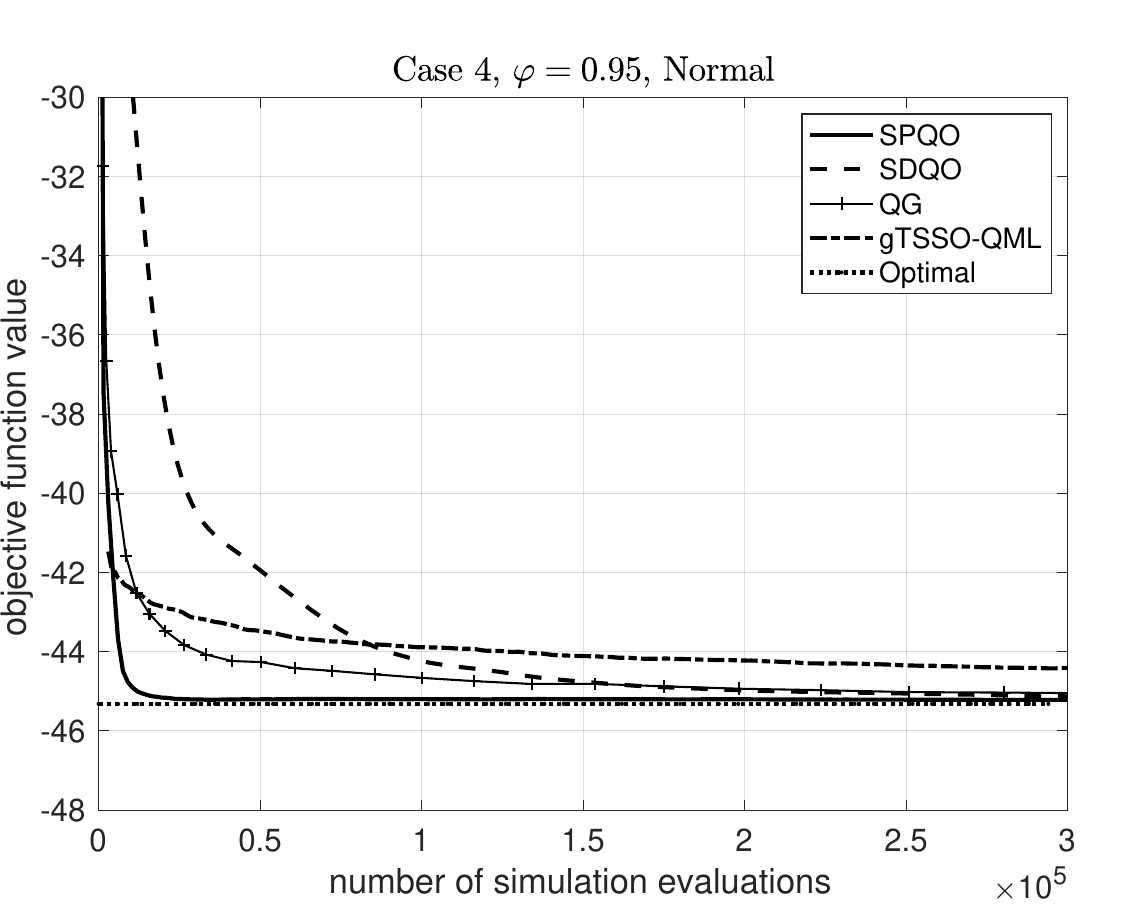}\\

   \includegraphics[trim=1cm 7.5cm 1cm 6.9cm, clip,width=0.45\textwidth]{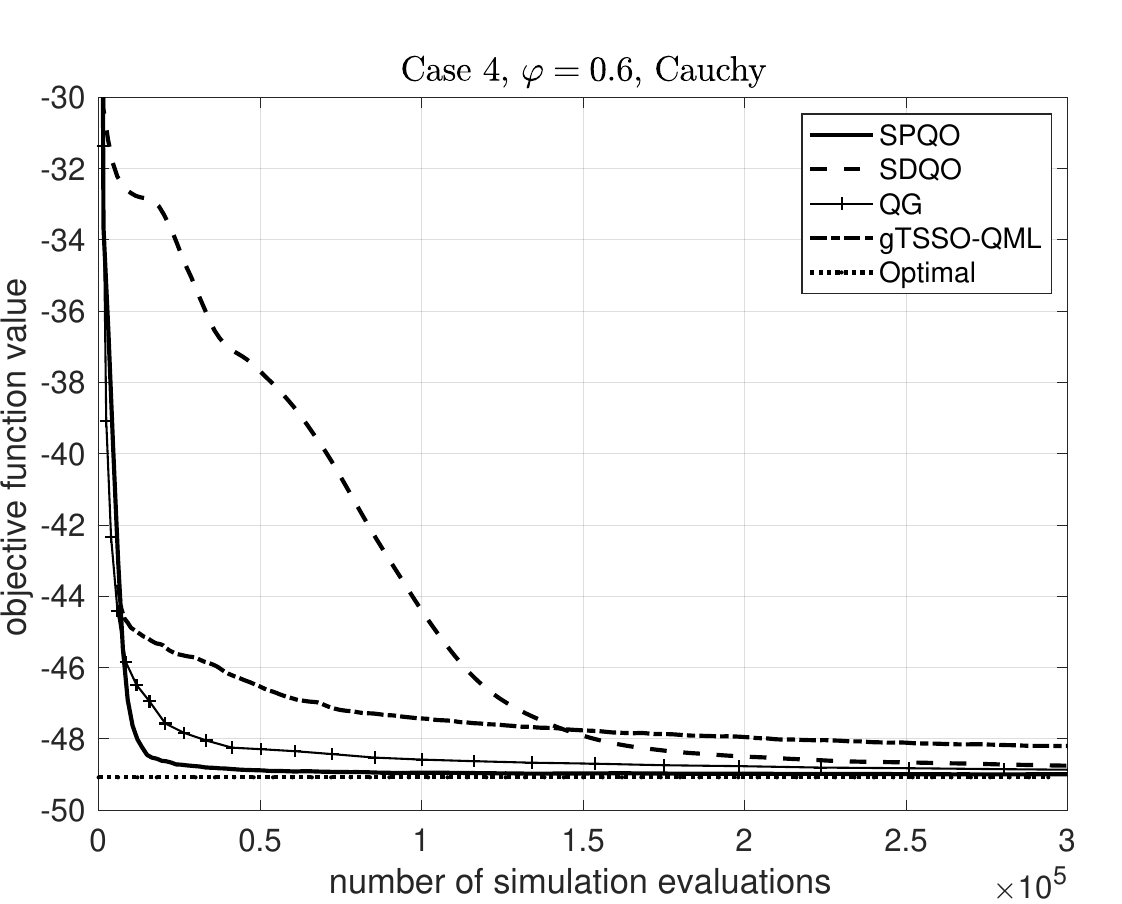}
	\includegraphics[trim=1cm 7.5cm 1cm 6.9cm, clip,width=0.45\textwidth]{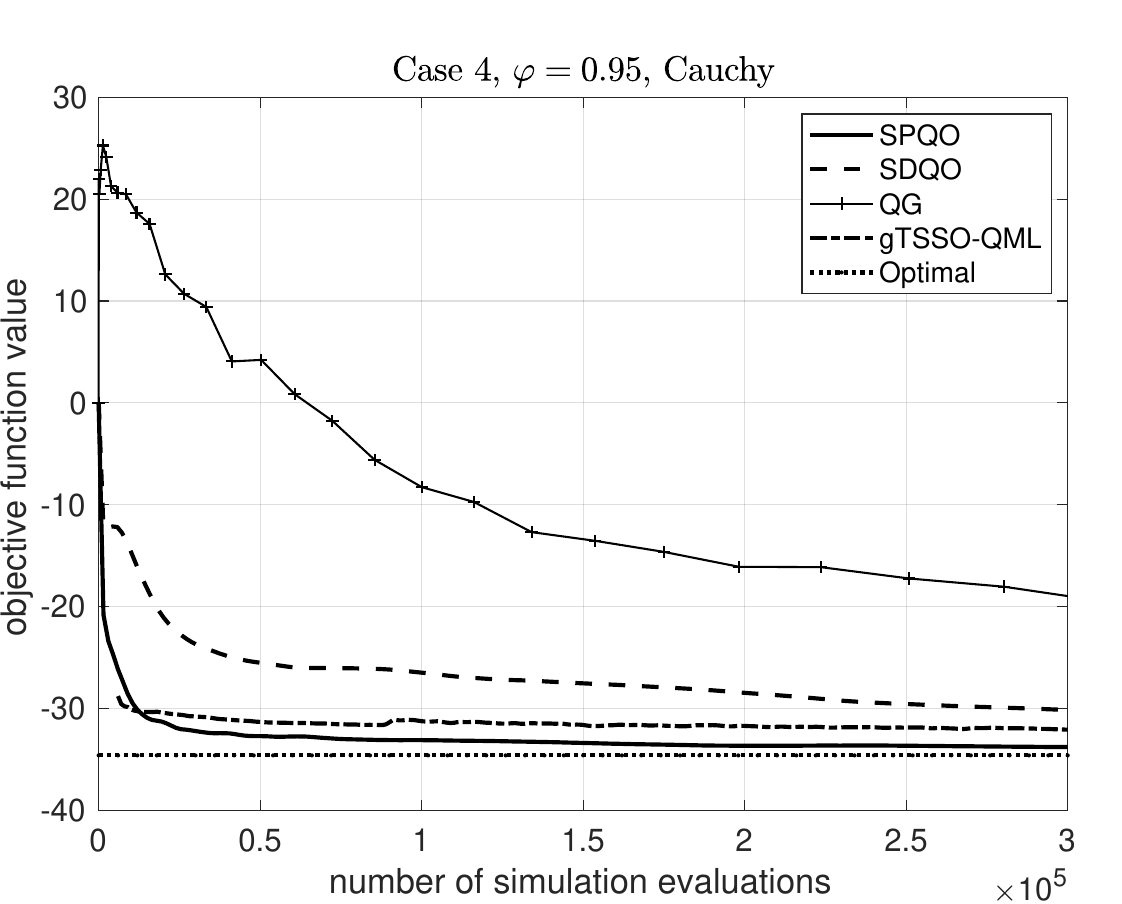}
	
	\caption{Performance of SPQO, SDQO, QG, and gTSSO-QML on test cases 3 and 4.}
	\label{fig2}
\end{figure}

\begin{figure}[!htbp]
	\centering
	\includegraphics[trim=1cm 7.5cm 1cm 6.9cm, clip,width=0.45\textwidth]{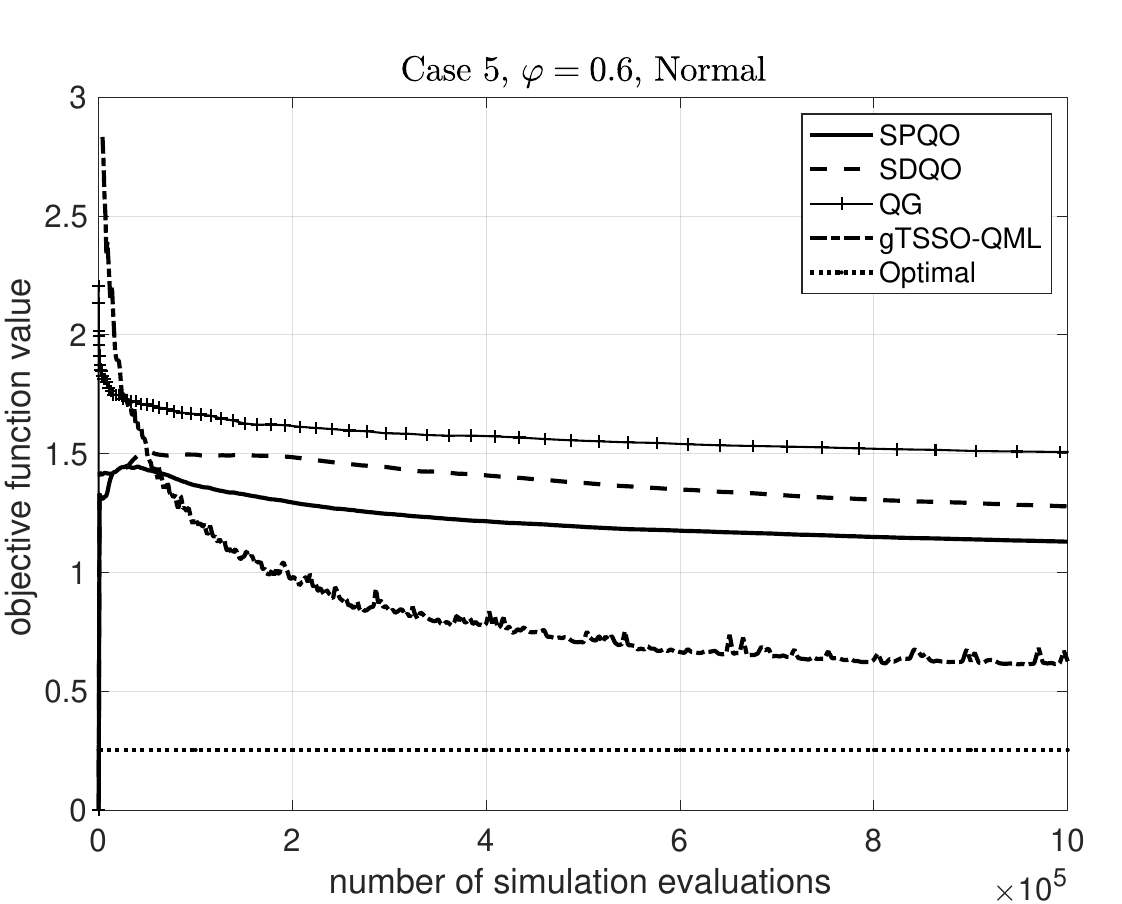}
	\includegraphics[trim=1cm 7.5cm 1cm 6.9cm, clip,width=0.45\textwidth]{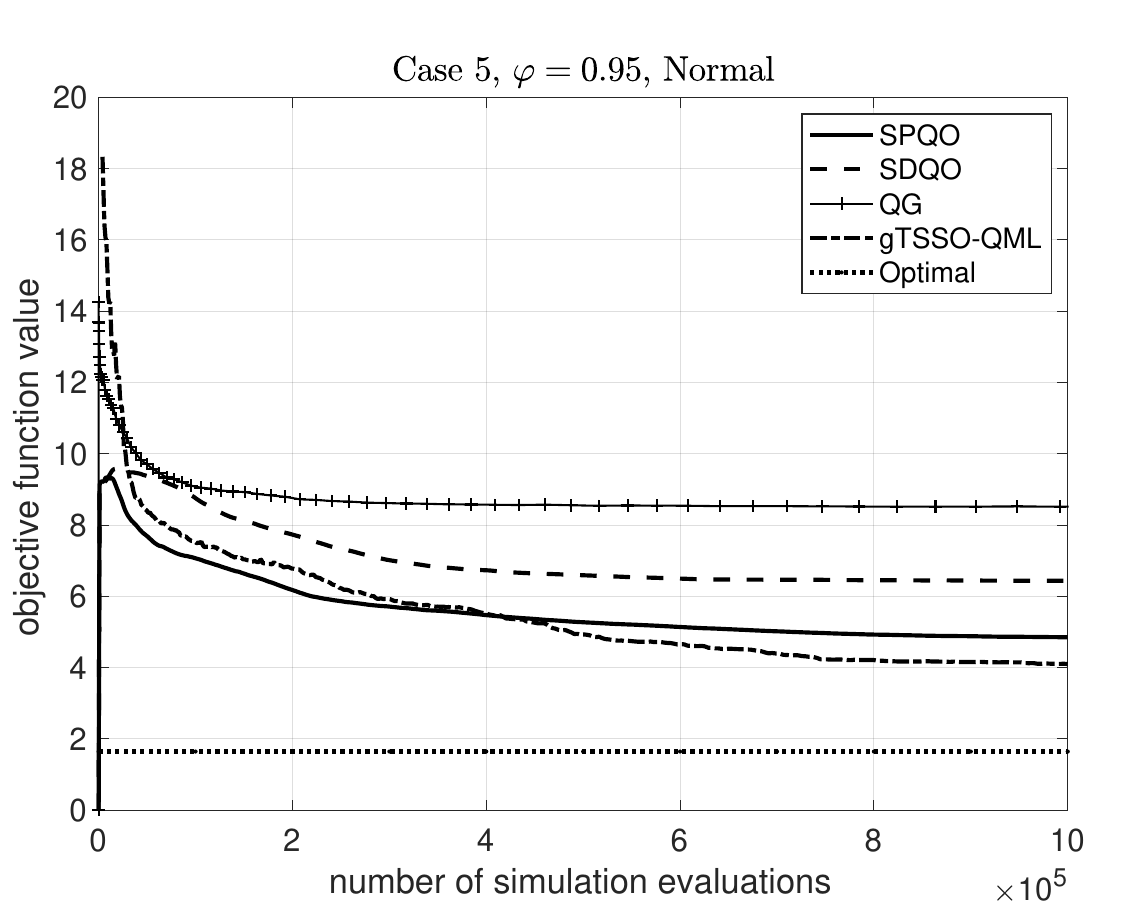}\\
	\vspace{1pt}
	\includegraphics[trim=1cm 7.5cm 1cm 6.9cm, clip,width=0.45\textwidth]{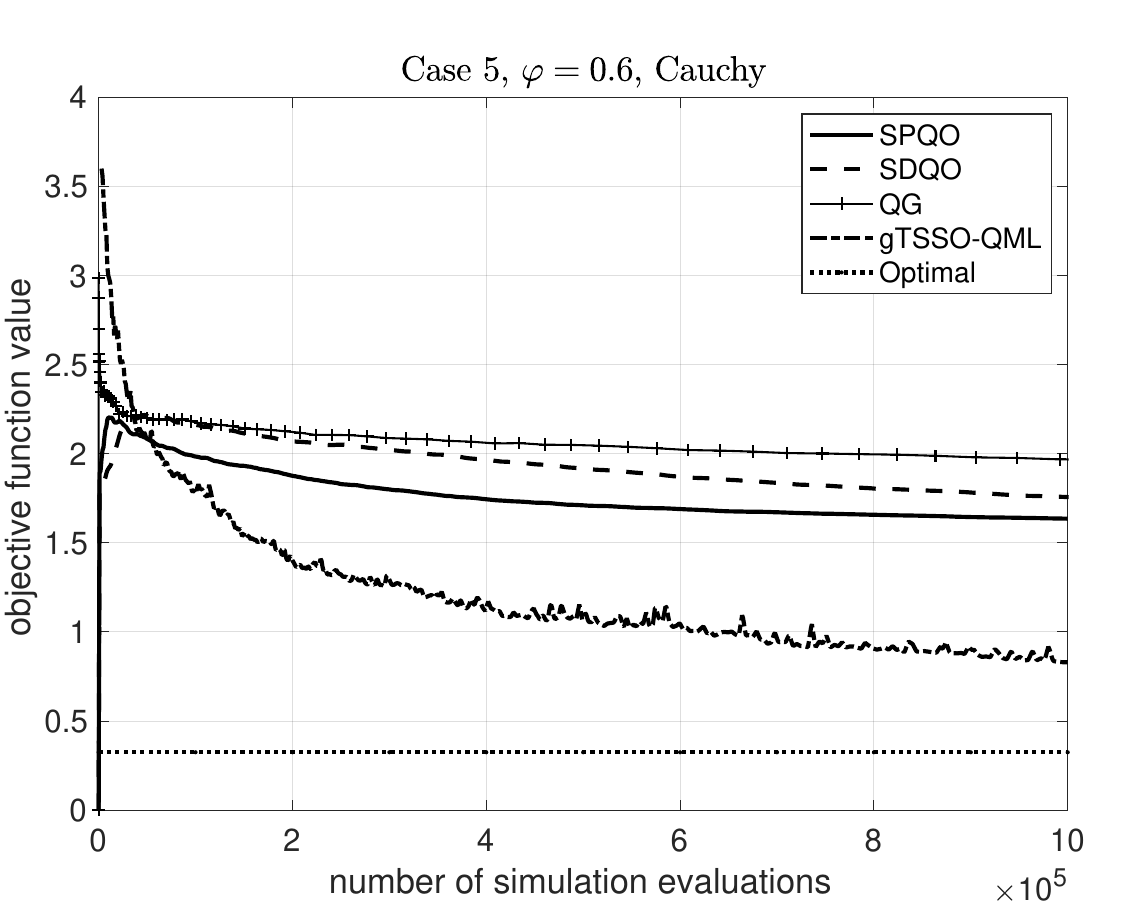}
	\includegraphics[trim=1cm 7.5cm 1cm 6.9cm, clip,width=0.45\textwidth]{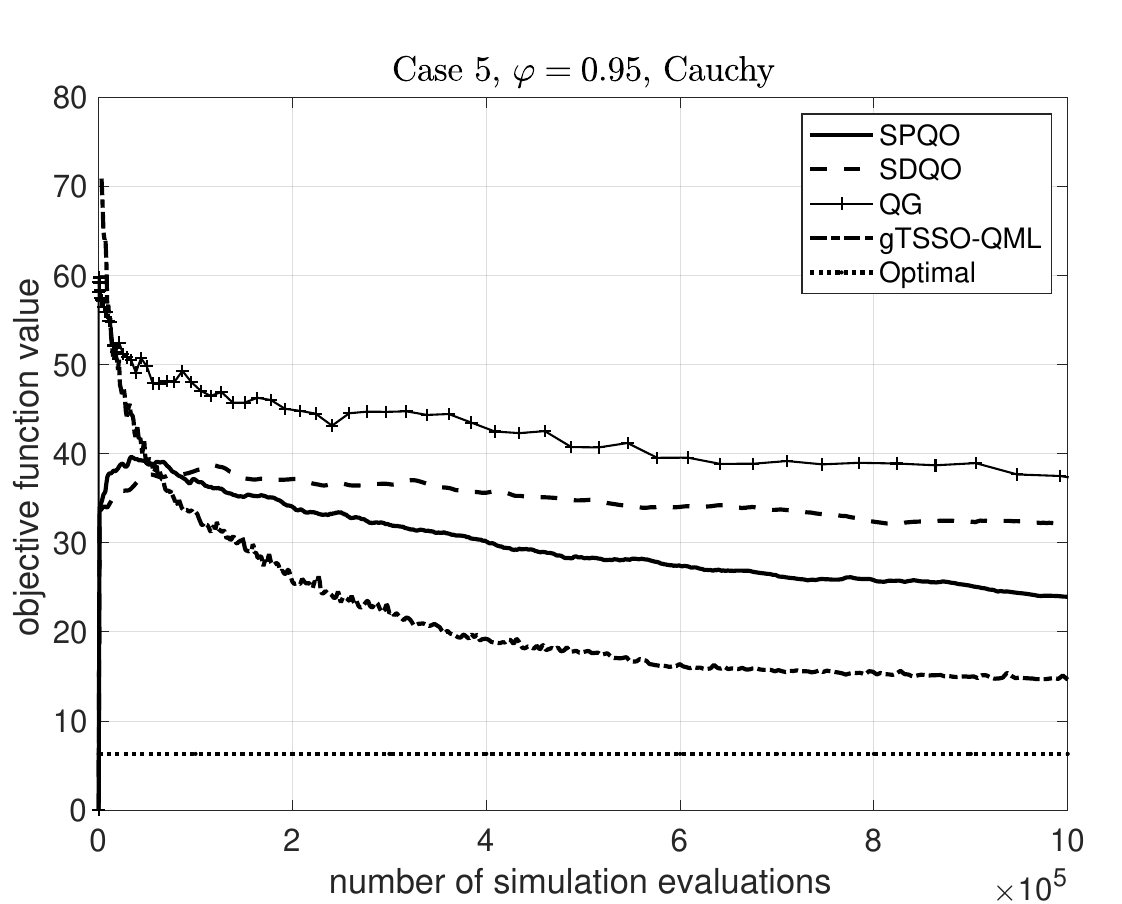}\\
	\vspace{1pt}
	\includegraphics[trim=1cm 7.5cm 1cm 6.9cm, clip,width=0.45\textwidth]{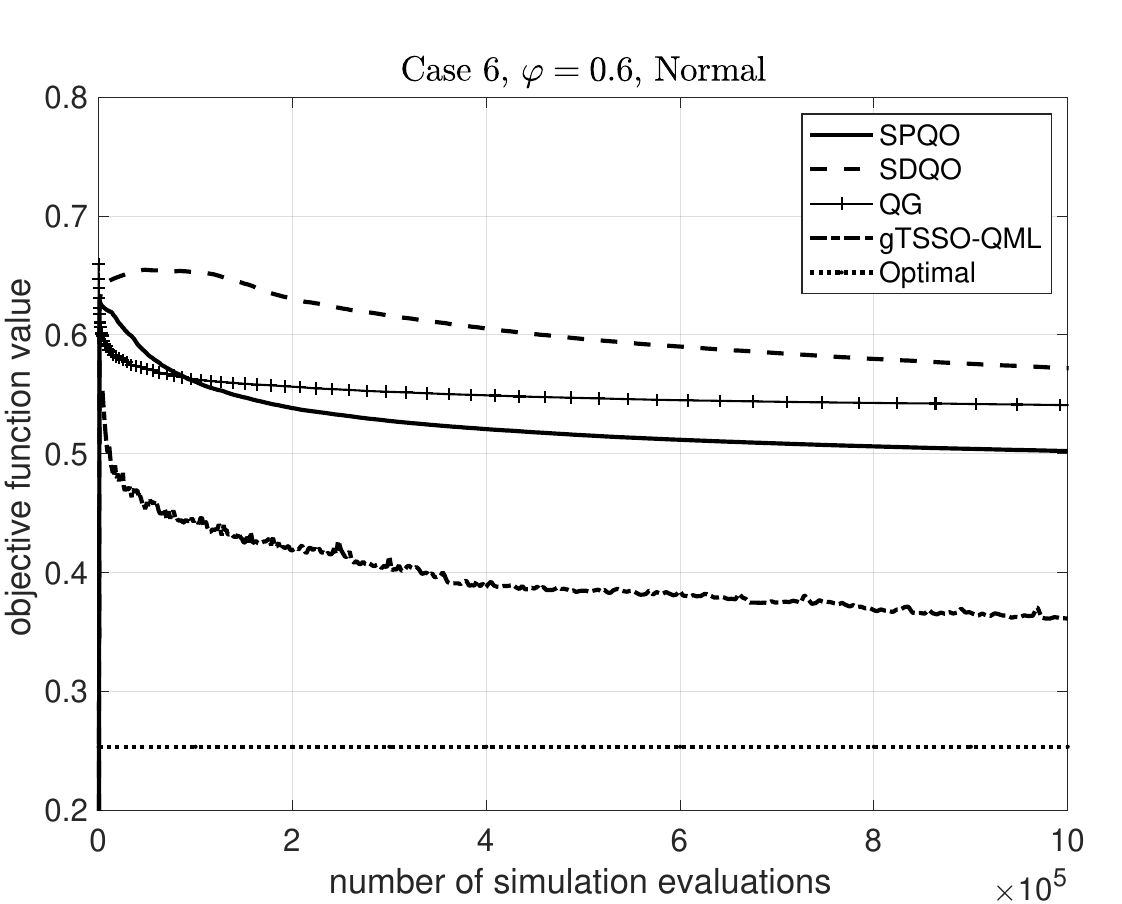}
	\includegraphics[trim=1cm 7.5cm 1cm 6.9cm, clip,width=0.45\textwidth]{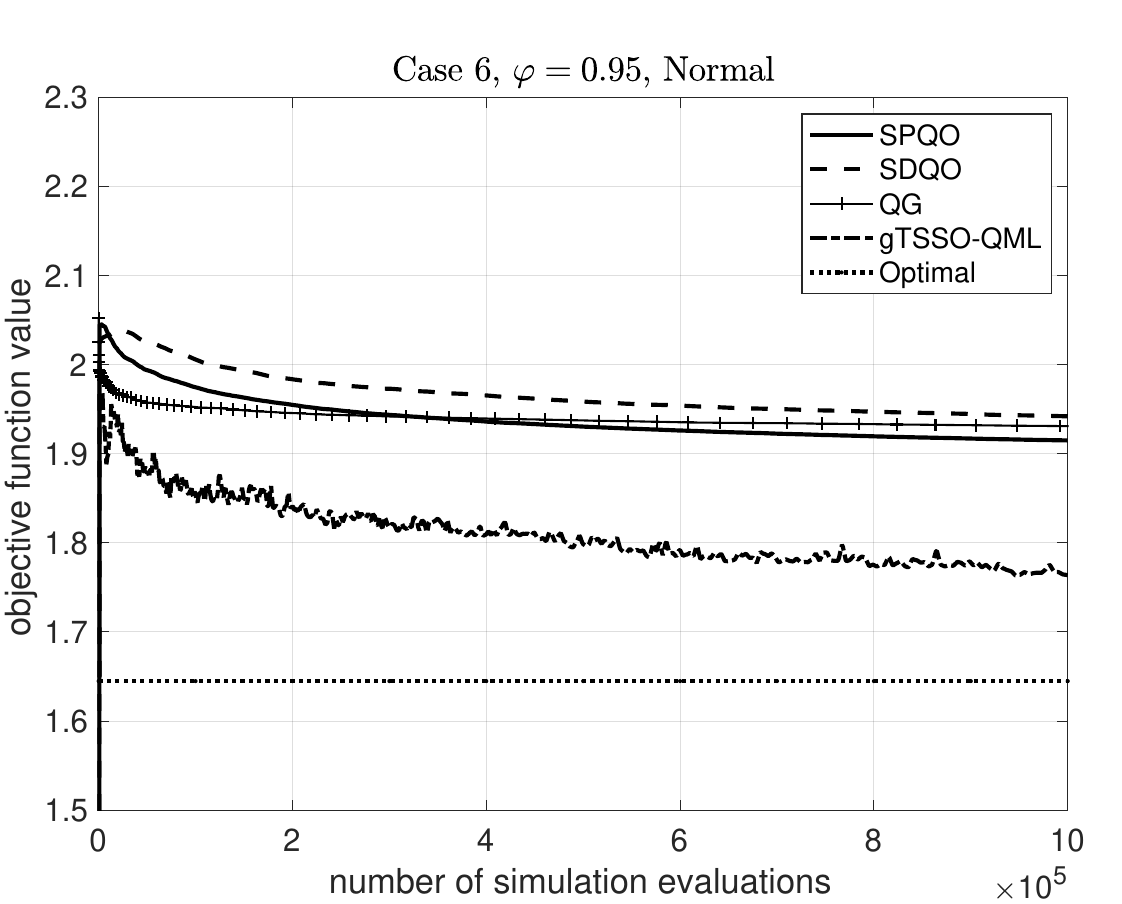}\\

   \includegraphics[trim=1cm 7.5cm 1cm 6.9cm, clip,width=0.45\textwidth]{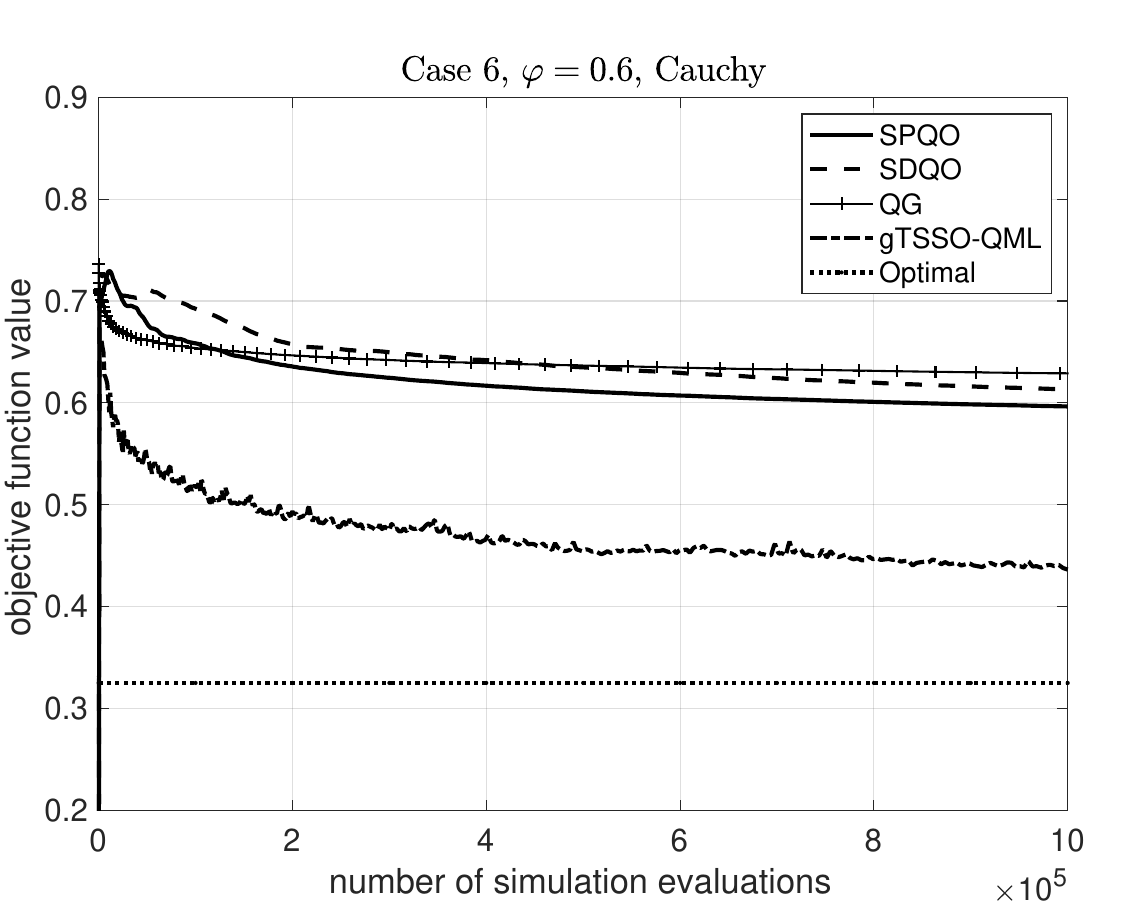}
	\includegraphics[trim=1cm 7.5cm 1cm 6.9cm, clip,width=0.45\textwidth]{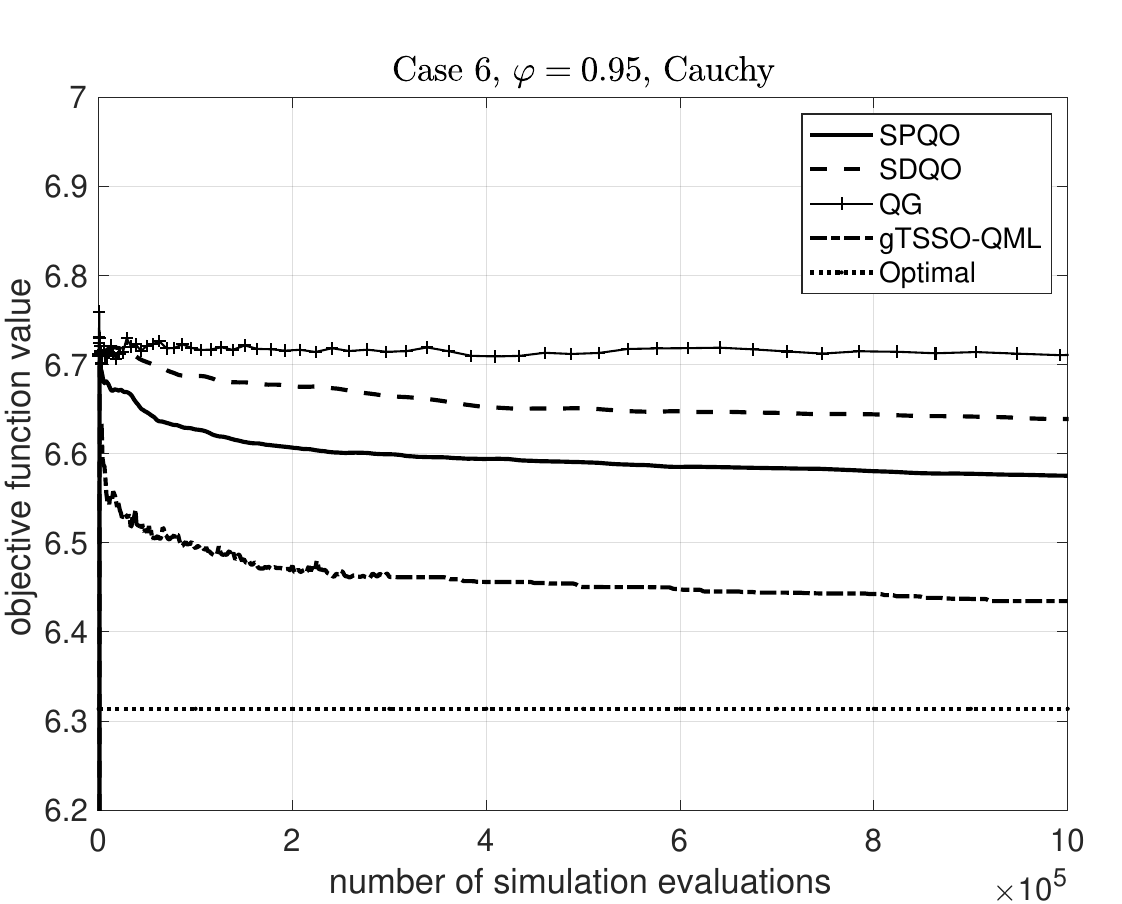}
	
	\caption{Performance of SPQO, SDQO, QG, and gTSSO-QML on test cases 5 and 6.}
	\label{fig3}
\end{figure}

It can be verified that the monotonicity condition assumed in Proposition~\ref{prop:crn} is satisfied for all test functions.
Therefore, in addition to SPQO and SDQO, we have also implemented their CRN versions: SPQO-CRN and SDQO-CRN. For each of the respective test cases, the six algorithms SPQO, SPQO-CRN, SDQO, SDQO-CRN, QG and gTSSO-QML are run using the same computational budget, where the total number of evaluations is set to $3\times 10^4$ for case 1, $3\times 10^5$ for {cases 2-4, and $10^6$ for the multi-modal test cases 5 and 6}. In SPQO (as well as SDQO and their CRN versions), the initial estimates are taken to be $D_0=(0,\ldots,0)^T$ and $q_0=0$. The initial $\theta_0$ is uniformly generated from $\Theta$ for all algorithms. Each algorithm is then independently repeated 40 times, and
the numerical results (averaged over 40 runs) obtained in the respective test
cases are presented in Tables {\ref{tab2}-\ref{tab3}}, which show the means and standard errors of the true quantile function values at the final solutions found by the six comparison algorithms.
In each {row} of the tables, the result that is closest to the true optimal value
is shown in bold (in case of a tie, the one with a smaller standard error is highlighted). The convergence behavior of SPQO, SDQO, QG, and gTSSO-QML is also illustrated in Figures~{\ref{fig1}-\ref{fig3}}, which plot the true quantile values at the current estimated solutions
as functions of the numbers of function evaluations consumed.

Our comparison results indicate that SPQO (SPQO-CRN) has the most consistency {on
test cases 1-4.} The final results obtained by SDQO are close to those of SPQO. However, as the figures clearly show, the convergence behavior of SDQO (in terms of the number of function evaluations) becomes slower as the problem dimension increases.
{For example, in cases 3 and 4, because SDQO uses 41 function evaluations per iteration, its total number of iterations is more than 13 times smaller than that of SPQO, resulting in inferior performance within the allowed budget.}
As expected, it can be seen from the tables that SPQO-CRN and SDQO-CRN generally outperform their original versions and yield smaller standard errors in almost all test cases, indicating their consistent performance over repeated runs.

The performance of QG is comparable to SPQO under the normal noise setting. In particular,
since QG is based on order statistics, the algorithm does not require the specification of an initial estimate and is less susceptible to the magnitude/location of the optimal quantile. This could be beneficial in extreme circumstances such as case 3, especially when
the true quantile happens to be very far away from its initial guessed value in a recursive procedure like SPQO. We see from {Figure~\ref{fig2}} that in case 3, QG quickly identifies the correct quantile range in the first few iterations and shows a very fast initial improvement. However, due to the increasing sample size required at each step, the parameter update in QG is carried out at a frequency that becomes much lower as search progresses, which results in sluggish performance in the long run. From {Table~\ref{tab2}}, we observe that the mean results found by SPQO are closer to the true optimal values than those obtained by QG.

Under the Cauchy noise, QG shows a significant performance degradation. We believe that this is owning to the heavy tail feature of the Cauchy distribution, which makes its high-level quantiles more difficult to estimate than the normal distribution. Therefore, for an order statistic-based estimator, a reliable quantile approximation can only be obtained after a large amount of simulation observations has been collected. This issue is especially manifested on cases 1, 2, and 4, where in each problem the standard Cauchy input distribution is further stretched by a large factor. We observe that at the $\varphi=0.95$ quantile level, QG may fail to locate a near-optimal solution within the prescribed simulation budget. SPQO instead estimates quantiles and gradients by averaging all simulation data collected in past iterations. Thus, it works equally well under the Cauchy noise setting and show significantly faster convergence behavior than QG.

{Functions 5 and 6 are highly multi-modal, which makes it very difficult for a gradient-based algorithm to escape local optima. On these test functions, SPQO, SDQO, and QG may quickly get stuck at a local optimum and stop making improvement even during the early search phase (see Figure~\ref{fig3}). The gTSSO-QML algorithm, instead, shows more robustness in dealing with local optima and yields much superior performance than other algorithms. However, note from Figure~\ref{fig2} that the algorithm is not as efficient on test cases 3 and 4.} We conjecture that this is primarily caused by the high dimensionality of these functions so that a close approximation of the true response curve might require a large amount of data that exceeds the given budget. In addition, because each step of gTSSO-QML involves an expensive optimization procedure, the algorithm {could be very} time consuming to run on high-dimensional problems. For example,
the average running time of gTSSO-QML for solving case 3 on the parallel platform is more than 8 hours. In contrast, the execution time of a single run of SPQO on a Windows PC is under 5 seconds. {The above comparison suggests that our proposed algorithms, particularly SPQO, are best suited to high-dimensional differentiable problems that contain few local optimal solutions, whereas gTSSO-QML is better adapted to the optimization of complex multi-modal objective functions with relatively small numbers of decision variables.}

\subsection{A Queueing Example}\label{sec52}

We consider a first-come, first-served single-server queue with parameterized service rate $\mu(\theta)=1/v^T\theta+\lambda$,
where $v\in \Re^d$ is a fixed positive vector and $\lambda$ is the arrival rate.
Denoting by $Y(\theta)$ the steady-state waiting time in system and $q_\varphi(\theta)$ the corresponding $\varphi$-quantile of $Y$, the objective is to determine an optimal parameter vector $\theta^*$ that minimizes the weighed cost of waiting and service given by
\begin{align}\label{cost}
y(\theta) = c_1 q_\varphi(\theta)+c_2(\theta-\vartheta)^T A (\theta-\vartheta),
\end{align}
where $c_1,\,c_2>0$ are cost coefficients, $\vartheta\in\Re^d$ is a nominal vector, and $A\in \Re^{d\times d}$ is a positive definite matrix; these are all assumed known.
The cost function (\ref{cost}) reflects the tradeoff between decreasing $\theta$ to increase the service rate (and hence reduce the waiting time quantile) and choosing $\theta$ to make the quadratic penalty term small.

Due to the cost-of-service penalty term, optimizing (\ref{cost}) becomes finding the zeros of $\nabla y(\theta)= 0$ rather than $\nabla q_\varphi(\theta)= 0$,
where
$$
\nabla y(\theta)= c_1 \nabla q_\varphi(\theta) + 2c_2A(\theta-\vartheta).
$$
Consequently, the three algorithms SPQO, SDQO, and QG are adjusted accordingly to solve this slightly modified root-finding problem.
All other steps of the algorithms remain intact.
For the simulation experiments, we take i.i.d. exponentially distributed interarrival times and services times, i.e., an M/M/1 queue,
with $\Theta=[1, 20]^4$, $\lambda=1$, $v=(0.1,0.2,0.3,0.4)^T$, $c_1=0.1$, $c_2=0.02$, $\vartheta=(7,8,9,10)^T$,
$$A  =
\begin{pmatrix}
	10~ & 2~ & 1~ & 2\\
	2~ & 9~ & 2~ & 4\\
	1~& 2~ & 8~ & 0 \\
    2~ & 4~ & 0~ & 7
	\end{pmatrix},$$
and consider two cases: $\varphi=0.5$ and $\varphi=0.95$.

\begin{table}[t]
\caption{\baselineskip8pt Performance on the queueing example, based on 40 independent runs (standard error in parentheses).}
  \label{tab:queue}
  \begin{center}
  \begin{tabular}{c c c}
    \hline
      &\multicolumn{1}{c}{$\varphi=0.50$} &\multicolumn{1}{c}{$\varphi=0.95$} \\
      \hline
      \multicolumn{1}{c}{Optimal Cost }& 0.62 & 2.66\\

      \multicolumn{1}{c}{SPQO} & 0.70 (1.2e-2) & 2.78 (1.9e-2) \\
      \multicolumn{1}{c}{SPQO-CRN} & \bf{0.67 (8.5e-3)} &\bf{2.75 (1.5e-2)} \\
      \multicolumn{1}{c}{SDQO} & 0.72 (1.6e-2) &2.80 (2.0e-2)\\
      \multicolumn{1}{c}{SDQO-CRN} & 0.73 (2.2e-2) &  2.78 (1.7e-2)\\
      \multicolumn{1}{c}{QG}& 1.17 (6.6e-2) & 3.57 (1.0e-1)\\
      \multicolumn{1}{c}{gTSSO-QML} & 0.87 (2.8e-2) & 3.19 (3.7e-2)\\
      \hline
  \end{tabular}
  \end{center}
\end{table}

\begin{figure}[thpb]
	\centering
\begin{tabular}{cc}
	\includegraphics[width=0.48\textwidth]{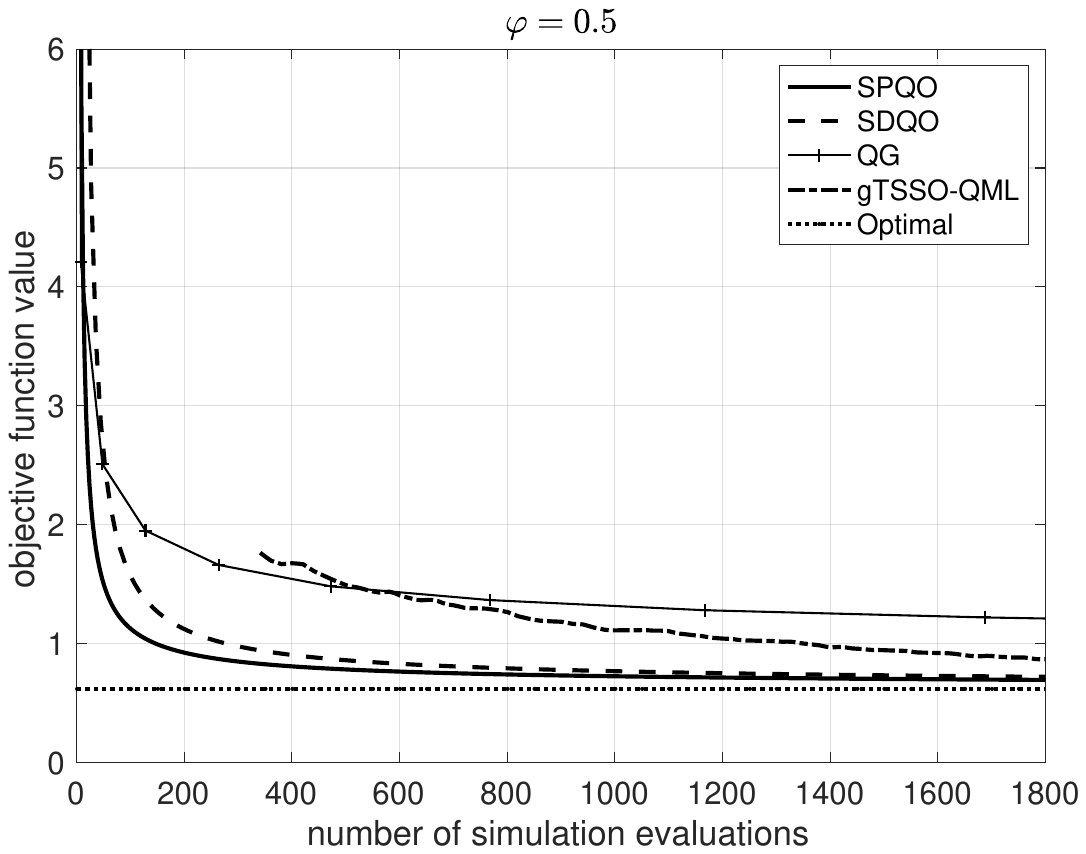} &
	\includegraphics[width=0.48\textwidth]{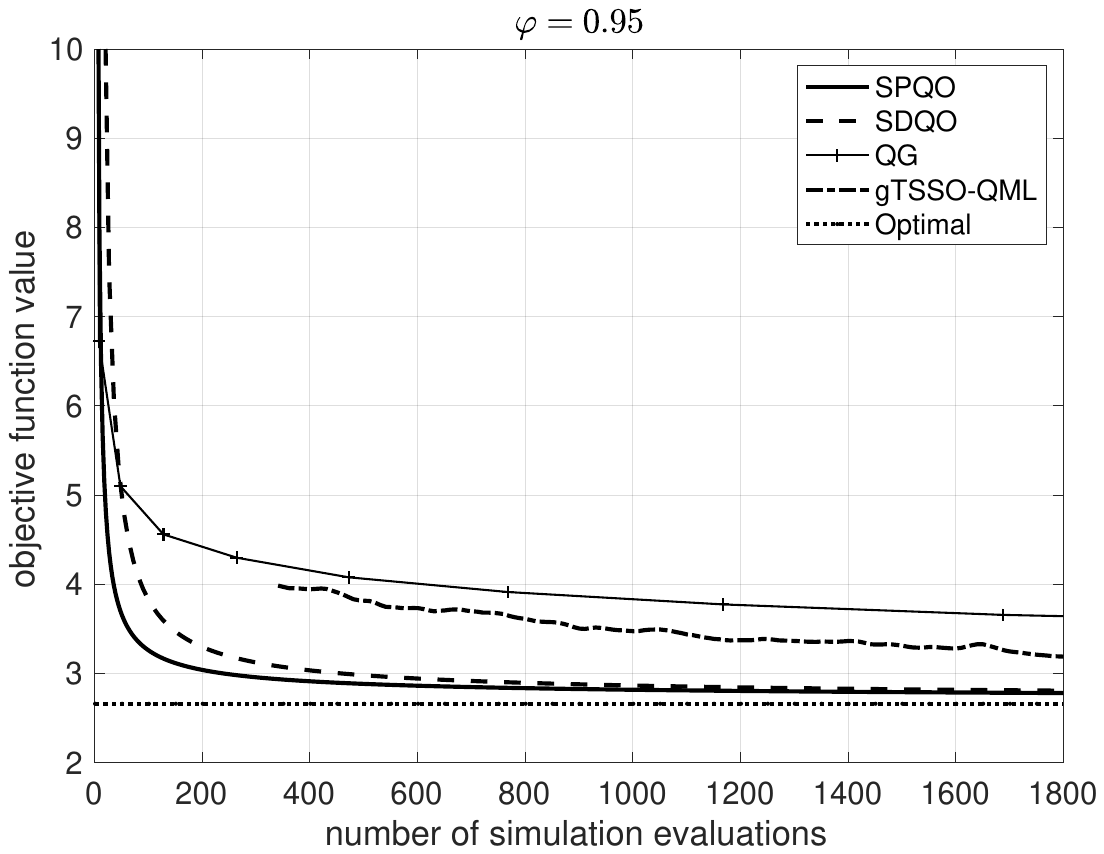} \\
\end{tabular}
	\caption{Performance of SPQO, SDQO, QG, and gTSSO-QML on the queueing example.}
	\label{fig:queue}
\end{figure}

All algorithms except gTSSO-QML are implemented using the same parameter settings
as for the black-box test functions in Section~\ref{sec51}, and the total number of simulation evaluations is set to $1800$.
The gTSSO-QML algorithm requires an initial set of design points and additional simulation samples to initialize the kriging model. In our implementation, we set the number of initial design points to {16} and allocate $20$ observations to each point to compute its sample quantile. {Thus, a total of $320$ simulation evaluations are used during the initialization step.}

The steady-state waiting time is approximated by the waiting time of the $1000$th customer. Each algorithm is independently repeated 40 times.
The simulation results (means and standard errors based on 40 independent runs of each algorithm) obtained on the two test cases are presented in Table~\ref{tab:queue}.
For the purpose of comparing with the true optimum,
we note that by basic queueing theory, the steady-state waiting time of an M/M/1 queue is exponentially distributed with parameter $\mu(\theta)-\lambda$, so the cost function can be calculated in closed form in terms of $\theta$ as
$y(\theta) =c_1\ln({1-\varphi})v^T\theta+c_2(\theta-\vartheta)^TA(\theta-\vartheta)$, which attains its minimum at  $\theta^*=\vartheta+\frac{c_1}{2c_2}\ln(1-\varphi)A^{-1}v$.
In Figure~\ref{fig:queue}, we also plot the true objective function values (averaged over 40 runs) obtained by SPQO, SDQO, QG, and gTSSO-QML as a function of the number of simulation evaluations.

The conclusions are generally consistent with the results for the black-box test functions {1-4} in Section~\ref{sec51}, with SPQO (SPQO-CRN) showing the best performance. SPQO and SDQO outperform QG by a large margin in terms of both mean performance and consistency (standard error). In particular, the large sample size required by QG results in only 8 iterations being carried out under the limited budget, whereas
the numbers of iterations for SPQO (SPQO-CRN) and SDQO (SDQO-CRN) are 600 and 200, respectively.
Moreover, because QG computes new quantile estimates at each iteration independently of past values, the observations collected in previous iterations are discarded, causing inefficient use of simulation data.
In contrast, SPQO and SDQO allow the quantile/gradient estimates to be constructed incrementally based on all historical simulation data, and this in turn offers superior finite-sample performance under a limited simulation budget.
{Note that because this is a uni-modal problem and gTSSO-QML does not exploit gradient information, it is not as efficient as SPQO/SDQO. However, the algorithm still outperforms QG, and its performance may be further improved through a more careful tuning of algorithm parameters.}


\section{Conclusions}\label{sec6}

For solving quantile optimization problems in the setting of noisy black-box functions, we have proposed two new three-timescale gradient-based SA algorithms.
For this quantile BBO setting, there are very few existing algorithms, so these algorithms represent a methodological contribution to the BBO literature.
These algorithms can also be applied to (stochastic) simulation optimization problems where direct gradient estimators based on techniques such as perturbation analysis or the likelihood ratio method, which rely on knowledge of the underlying model, are not readily available or are difficult to implement, so the algorithms also advance the state of the art in simulation optimization.
The SPQO algorithm is especially promising for high-dimensional problems, requiring only three function evaluations per iterative update. 
Compared with methods relying on order statistics, the algorithms proposed here have the potential to achieve substantial computational savings. Variants of the algorithms using CRN offer the opportunity for further reductions in the variance of the gradient estimator and hence faster convergence of the algorithms in the simulation optimization setting.

Under the assumption of differentiability of the quantile function and other appropriate conditions, we have analyzed the bias effect of the proposed gradient estimation scheme and established the local convergence of the resultant algorithms.
More importantly, through a novel fixed-point argument, we have also provided detailed characterizations of the convergence rates of the quantile and quantile gradient estimators. These results extend existing work in the single-timescale setting and indicate that
{an upper bound on the MAEs of the algorithms diminishes at the optimal rate  $O(k^{-1/4})$.}
Simulation experiments indicate that the algorithms perform well, and, in particular, SPQO is very promising for solving high-dimensional problems, in terms of the sample size required to achieve reasonable performance.

Future avenues of potential research building on the results here include (i) developing and analyzing other multi-timescale algorithms for quantile BBO, such as the two-timescale version alluded to in the introduction, for comparison;
(ii) applying the fixed-point approach to study the convergence rate of other multi-timescale SA algorithms;
(iii) investigating a more systematic way of tuning the step-size parameters needed to implement multi-timescale SA algorithms;
and (iv) designing a more comprehensive computational/experimental study to characterize when the multi-timescale approach is most effective.






\bibliographystyle{informs2014}
\bibliography{qospsa}

\clearpage
\begin{APPENDICES}
\section{Proof of Lemma~\ref{lem:dk}}\label{appendix_dk}
We have from (\ref{dk}) that
\begin{align}\label{tmp1}\nonumber
\|D_{k+1}\|^2&=\|D_k\|^2+\frac{\beta_kM_k}{c_k}(-I_k^++I_k^-)\Delta_k^TD_k+\frac{M_k^2}{4}\frac{\beta_k^2}{c_k^2}(-I_k^++I_k^-)^2\Delta_k^T\Delta_k\\
&\leq \|D_k\|^2+\frac{\beta_kM_k}{c_k}(-I_k^++I_k^-)\Delta_k^TD_k+\frac{M_k^2}{4}\frac{\beta_k^2}{c_k^2}d \nonumber\\
&=\|D_k\|^2+\frac{\beta_kM_k}{c_k}(-I_k^++I_k^-)\Delta_k^TD_k+\frac{\max\{d,\|D_k \|^2\}}{4}\frac{\beta_k^2}{c_k^2} \nonumber\\
&\leq \Big(1+\frac{\beta_k^2}{4c^2_k}\Big)\|D_k\|^2+\frac{\beta_kM_k}{c_k}(-I_k^++I_k^-)\Delta_k^TD_k+\frac{d}{4}\frac{\beta_k^2}{c_k^2}
\end{align}
where the first inequality above follows from $(-I_k^++I_k^-)^2\leq 1$ and $\Delta_k^T\Delta_k=d$ by A2.
Taking conditional expectation at both sides and then applying Taylor's theorem, we have
\begin{align*}
E[\|D_{k+1}\|^2|\mathcal{F}_k]&\leq\Big(1+\frac{\beta_k^2}{4c^2_k}\Big)\|D_k\|^2+\frac{\beta_kM_k}{c_k}E\big[(-I_k^++I_k^-)\Delta_k^T\big|\mathcal{F}_k\big]D_k+
\frac{d}{4}\frac{\beta_k^2}{c_k^2}\\
&\hspace{-2.5cm}=\Big(1+\frac{\beta_k^2}{4c^2_k}\Big)\|D_k\|^2+\frac{\beta_kM_k}{c_k}E\Big[E\big[(-I_k^++I_k^-)|\Delta_k,\mathcal{F}_k\big]\Delta_k^T\big|\mathcal{F}_k\Big]D_k+
\frac{d}{4}\frac{\beta_k^2}{c_k^2}\\
&\hspace{-2.5cm}=\Big(1+\frac{\beta_k^2}{4c^2_k}\Big)\|D_k\|^2+\frac{\beta_kM_k}{c_k}E\big[(-F^+_k+F^-_k)\Delta^T_k|\mathcal{F}_k\big]D_k+\frac{d}{4}\frac{\beta_k^2}{c_k^2}\\
&\hspace{-2.5cm}=\Big(1+\frac{\beta_k^2}{4c^2_k}\Big)\|D_k\|^2+\frac{\beta_kM_k}{c_k}E\Big[(-f(\bar q_k^+;\bar{\theta}_k^+)-f(\bar q_k^-;\bar{\theta}_k^-) )c_k\frac{(D_k^T\Delta_k)^2}{M_k} \Big|\mathcal{F}_k\Big]\\
&-\frac{\beta_kM_k}{c_k}E\Big[\big(\nabla^T_{\theta}F(\bar q_k^+;\bar{\theta}_k^+)+\nabla^T_{\theta}F(\bar q_k^-;\bar{\theta}_k^-) \big)\bar c_k {\Delta}_k\Delta_k^TD_k \Big|\mathcal{F}_k\Big]+\frac{d}{4}\frac{\beta_k^2}{c_k^2},
\end{align*}
where $\bar q_k^+$, $\bar q_k^-$ are on the line segments between $q_k$ and $q_k\pm \bar c_k {D}_k^T\Delta_k$, and $\bar{\theta}_k^+$, $\bar{\theta}_k^-$ are on the line segments connecting $\theta_k$ and $\theta_k\pm \bar c_k{\Delta}_k$. Since $c_k\rightarrow 0$ and $\bar c_k\|{D}_k\|\leq c_k\sqrt{d}$, $\bar c_k\|{\Delta}_k\|\leq c_k\|\Delta_k\|= c_k\sqrt{d}$, A1 implies that there exist a constant $\mathcal{B}>0$ and an integer $N>0$ such that $f(\bar q_k^+;\bar{\theta}_k^+),\,f(\bar q_k^-;\bar{\theta}_k^-)\geq \epsilon$ and $\|\nabla_{\theta}F(\bar q_k^+;\bar{\theta}_k^+)\|,\,\|\nabla_{\theta}F(\bar q_k^-;\bar{\theta}_k^-)\|\leq \mathcal{B}$ for all $k\geq N$. Next, by noting that
$E[(D_k^T\Delta_k)^2|\mathcal{F}_k]=D_k^TE[\Delta_k \Delta_k^T|\mathcal{F}_k ]D_k=\|D_k\|^2$ (because $E[\Delta_{k,i}\Delta_{k,j}|\mathcal{F}_k]=0$ $\forall\,i\neq j$) and applying the Cauchy-Schwarz inequality, we have that for all $k\geq N$,
\begin{align*}
E[\|D_{k+1}\|^2|\mathcal{F}_k]&\leq \Big(1+\frac{\beta_k^2}{4c^2_k}\Big)\|D_k\|^2-2\beta_k\epsilon E\big[(D_k^T\Delta_k)^2\big|\mathcal{F}_k\big]\\
&+\beta_kE\Big[\big|\big(\nabla^T_{\theta}F(\bar q_k^+;\bar{\theta}_k^+)+\nabla^T_{\theta}F(\bar q_k^-;\bar{\theta}_k^-) \big){\Delta}_k\Delta_k^TD_k\big| \Big|\mathcal{F}_k\Big]+\frac{d}{4}\frac{\beta_k^2}{c_k^2}\\
&\leq \Big(1+\frac{\beta_k^2}{4c^2_k}-2\beta_k\epsilon\Big)\|D_k\|^2+2\beta_k \mathcal{B}d \|D_k\|+\frac{d}{4}\frac{\beta_k^2}{c_k^2}.
\end{align*}
From A3(a), we know that $\beta_k^2/c_k^2=o(\beta_k)$. Therefore, there exist an $\epsilon'\in(0,\epsilon)$ and an integer $N'>0$ such that $2\beta_k\epsilon'<1$ and $\beta_k^2/{4c^2_k}-2\beta_k\epsilon\leq -2\beta_k\epsilon'$ for all $k\geq N'$. Taking expectations at both sides, we thus obtain that for all
$k\geq \bar N:=\max\{N,N'\}$,
\begin{align*}
E[\|D_{k+1}\|^2]&\leq (1-2\beta_k\epsilon')E[\|D_k\|^2]+2\beta_k \mathcal{B}d\sqrt{E[\|D_k\|^2 ]}+\frac{d}{4}\frac{\beta_k^2}{c_k^2}\\
&\leq \Big(\sqrt{1-2\beta_k\epsilon'}\sqrt{E[\|D_k\|^2]}+\beta_k\epsilon'\frac{\mathcal{B}d}{\epsilon'\sqrt{1-2\beta_k\epsilon'}} \Big)^2+\frac{d}{4}\frac{\beta_k^2}{c_k^2}\\
&\leq \big((1-\beta_k\epsilon')\sqrt{E[\|D_k\|^2]}+\beta_k\epsilon' C_D\big)^2+\frac{d}{4}\frac{\beta_k^2}{c_k^2}\\
&\leq \max\{\sqrt{E[\|D_k\|^2]},C_D\}^2+\frac{d}{4}\frac{\beta_k^2}{c_k^2}\\
&=\max\{E[\|D_k\|^2],C_D^2\}+\frac{d}{4}\frac{\beta_k^2}{c_k^2}\\
&\ldots \\
&\leq \max\{E[\|D_{\bar N}\|^2],C^2_D \}+\frac{d}{4}\sum_{i=\bar N}^k\frac{\beta_i^2}{c_i^2}\\
&<\infty,
\end{align*}
where the third step follows from the inequality $\sqrt{1-x}\leq 1-x/2$ for $x\in [0,1]$, $C_D:=\frac{\mathcal{B}d}{\epsilon'\sqrt{1-2\beta_{\bar N}\epsilon'}}$, and the finiteness of the last bound is a direct consequence of A3(a).
This shows $\sup_kE[\|D_k \|^2 ]<\infty$.

To show part ($ii$) of the lemma, we let
$e^+_k:=F^+_k-I^+_k$, $e^-_k:=I^-_k-F^-_k$, and $e_k:=\|D_k\|^2-(D_k^T\Delta_k)^2$. Using a similar argument as before, we have that for all $k\geq \bar N$, the term on the right-hand side of (\ref{tmp1}) can be bounded by
\begin{align*}
\|D_{k+1}\|^2 &\leq \Big(1+\frac{\beta_k^2}{4c^2_k}\Big)\|D_k\|^2+\frac{\beta_kM_k}{c_k}(-F_k^++F_k^-)\Delta_k^TD_k+\frac{\beta_kM_k}{c_k}(e^+_k+e^-_k)\Delta^T_kD_k +\frac{d}{4}\frac{\beta_k^2}{c_k^2} \\
&=\Big(1+\frac{\beta_k^2}{4c^2_k}\Big)\|D_k\|^2+\beta_k\Big(-f(\bar q_k^+;\bar{\theta}_k^+)-f(\bar q_k^-;\bar{\theta}_k^-) \Big)(D_k^T\Delta_k)^2\\
&\hspace{0.5cm}-\beta_k \Big(\nabla_{\theta}^TF(\bar q_k^+;\bar{\theta}_k^+)+\nabla_{\theta}^TF(\bar q_k^-;\bar{\theta}_k^-) \Big)\Delta_k\Delta_k^TD_k+\frac{\beta_kM_k}{c_k}(e^+_k+e^-_k)\Delta^T_kD_k +\frac{d}{4}\frac{\beta_k^2}{c_k^2} \\
&\leq \Big(1+\frac{\beta_k^2}{4c^2_k}\Big)\|D_k\|^2-2\beta_k\epsilon (D_k^T\Delta_k)^2+\beta_k \Big|\Big(\nabla_{\theta}^TF(\bar q_k^+;\bar{\theta}_k^+)+\nabla_{\theta}^TF(\bar q_k^-;\bar{\theta}_k^-) \Big)\Delta_k\Delta_k^TD_k\Big|\\
&\hspace{0.5cm}+\frac{\beta_kM_k}{c_k}(e^+_k+e^-_k)\Delta^T_kD_k +\frac{d}{4}\frac{\beta_k^2}{c_k^2} \\
&\leq \Big(1+\frac{\beta_k^2}{4c^2_k}\Big)\|D_k\|^2-2\beta_k\epsilon (D_k^T\Delta_k)^2+2\beta_k\mathcal{B}d\|D_k \|+\frac{\beta_kM_k}{c_k}(e^+_k+e^-_k)\Delta^T_kD_k +\frac{d}{4}\frac{\beta_k^2}{c_k^2}\\
&\leq (1-2\beta_k\epsilon')\|D_k\|^2+2\beta_k\epsilon e_k+2\beta_k\mathcal{B}d\|D_k \|+\frac{\beta_kM_k}{c_k}(e^+_k+e^-_k)\Delta^T_kD_k +\frac{d}{4}\frac{\beta_k^2}{c_k^2}\\
&\leq \Big((1-\beta_k\epsilon')\|D_k\|+\beta_k\epsilon' C_D \Big)^2+2\beta_k\epsilon e_k+ \frac{\beta_kM_k}{c_k}(e^+_k+e^-_k)\Delta^T_kD_k +\frac{d}{4}\frac{\beta_k^2}{c_k^2}\\
&\leq \max\{\|D_k\|^2,C_D^2 \}+2\beta_k\epsilon e_k+ \frac{\beta_kM_k}{c_k}(e^+_k+e^-_k)\Delta^T_kD_k +\frac{d}{4}\frac{\beta_k^2}{c_k^2},
\end{align*}
where $\epsilon'$ and $C_D$ are defined as in the proof of part ($i$) of the lemma. Let
$\xi_k:=2\beta_k\epsilon e_k+ \frac{\beta_kM_k}{c_k}(e^+_k+e^-_k)\Delta^T_kD_k$ and $\xi'_k:=\frac{d}{4}\frac{\beta_k^2}{c_k^2}$. Expanding the above inequality recursively, we obtain
\begin{align}\label{tmp2}
\|D_{k+1}\|^2&\leq \max\big\{\|D_k\|^2+\xi_k+\xi'_k,C^2_D+\xi_k+\xi'_k  \big\}\nonumber \\
&\leq \max\Bigg\{\max\Big\{\|D_{k-1}\|^2+\sum_{i=k-1}^k\xi_i+\sum_{i=k-1}^k\xi'_i,C_D^2+ \sum_{i=k-1}^k\xi_i+\sum_{i=k-1}^k\xi'_i\Big\}, C^2_D+\xi_k+\xi'_k \Bigg\}\nonumber\\
&\leq \max\Bigg\{\max\Big\{\|D_{k-1}\|^2+\sum_{i=k-1}^k\xi_i,C_D^2+\sum_{i=k-1}^k\xi_i \Big\},C^2_D+\xi_k \Bigg\}+\sum_{i=k-1}^k\xi'_i \nonumber\\
&=\max\Bigg\{\|D_{k-1}\|^2+\sum_{i=k-1}^k\xi_i, C^2_D+\max\Big\{\sum_{i=k-1}^k\xi_i,\xi_k \Big\} \Bigg\}+\sum_{i=k-1}^k\xi'_i \nonumber\\
&\ldots \nonumber\\
&\leq \max\Bigg\{\|D_{\bar N}\|^2+\sum_{i=\bar N}^k\xi_i, C_D^2+\max\Big\{\xi_k,\ldots,\sum_{i=\bar N}^k\xi_i \Big\} \Bigg\}+\sum_{i=\bar N}^k\xi'_i.
\end{align}

By squaring both sides of (\ref{tmp1}) and applying essentially the same argument as in the proof of part ($i$) (noting that $E[\|D_k\|^3]\leq \sqrt{E[\|D_k\|^4 ]}\sqrt{E[\|D_k \|^2 ]}$), it is tedious but straightforward to show that $\sup_k E[\|D_k\|^4]<\infty$. 
This, together with A3(a) and the fact $E[\xi_k|\mathcal{F}_k]=0$, indicates that $\sum_{i=0}^k\xi_i$ is an $L^2$-bounded martingale and thus converges to a finite random variable for almost all $\omega\in\Omega$ \citep[e.g.,][]{shiry96}. Along each such $\omega$, the sequence $\{\sum_{i=0}^k \xi_i\}$ is Cauchy convergent. Therefore, for a given $C>0$, there is an integer $N''>0$ such that $|\sum_{i=k'}^k\xi_i|=|\sum_{i=0}^k \xi_i-\sum_{i=0}^{k'-1}\xi_i|\leq C$ for all $k\geq k'\geq N''$. Consequently, the absolute values of both terms $\sum_{i=\bar N}^k\xi_i$ and $\max\big\{\xi_k,\ldots,\sum_{i=\bar N}^k\xi_i \big\}$ in (\ref{tmp2}) can be made smaller than $C$ by taking $\bar N\geq N''$. Hence, the proof is completed by noticing that $\sum_{i=\bar N}^k\xi'_i=\frac{d}{4}\sum_{i=\bar N}^k\beta_i^2/c_i^2<\infty$ by A3(a). \Halmos

\section{Proof of Theorem ~\ref{thm:con}}\label{appendix_thm1}

Lemma~\ref{lem:bias} allows us to put (\ref{dk}) in the form of a recursion in \cite{huetal2021}:
\begin{equation}\label{tmp3}
D_{k+1}=D_k+\beta_k\Big(G_{k,1}(Y^+_k,Y^-_k;q_k,\theta_k)-G_{k,2}(Y^+_k,Y^-_k;q_k,\theta_k)D_k\Big),
\end{equation}
where $G_{k,2}:=f(q_k;\theta_k)$ and $G_{k,1}:=-\nabla_{\theta}F(q_k;\theta_k)+\big(\frac{-I_k^++I_k^-}{2\bar c_k{\Delta}_k}-E[\frac{-I_k^++I_k^-}{2\bar c_k{\Delta}_k} |\mathcal{F}_k] \big)+b_k(q_k,D_k,\theta_k)$. Thus, the results can be proved by verifying that the conditions used in \cite{huetal2021} hold. In particular, it has been shown that for coupled recursions of the forms (\ref{qk}), (\ref{tmp3}), and (\ref{projection}), the stated convergence properties of $\{\theta_k\}$ hold if $i)$ there are constants $\epsilon\in(0,1)$ and $C_2>0$ such that $\inf_kE[G_{k,2}|\mathcal{F}_k]\geq \epsilon$ and $\sup_k E[G_{k,2}^2|\mathcal{F}_k]\leq C_2$ w.p.1; $ii)$ there exists a constant $C_1>0$ such that $\sup_k E[\|G_{k,1}\|^2|\mathcal{F}_k]\leq C_1$ w.p.1; $iii)$ $E[G_{k,2}|\mathcal{F}_k]-f(q_k;\theta_k)\rightarrow 0$ and  $E[G_{k,1}|\mathcal{F}_k]+\nabla_{\theta}F(q_k;\theta_k)\rightarrow 0$ as $k\rightarrow \infty$ w.p.1.

Evidently, $i)$ follows trivially from A1 and our definition of $G_{k,2}$, whereas $iii)$ is a direct consequence of part ($i$) of Lemma~\ref{lem:bias}, A3(a) ($c_k\rightarrow 0$), and the fact that \begin{equation}\label{tmp4}
E\Big[\Big(\frac{-I_k^++I_k^-}{2\bar c_k{\Delta}_k}-E\big[\frac{-I_k^++I_k^-}{2\bar c_k{\Delta}_k} \big|\mathcal{F}_k\big]\Big)\Big|\mathcal{F}_k \Big]=0.
\end{equation}
Regarding $ii)$, we note that this assumption is in fact stronger than what is needed and can instead be relaxed to the weaker condition $\sum_{k=0}^{\infty}\beta_k^2E[\|G_{k,1}\|^2]<\infty$. It is then easy to check that all steps in the convergence proof given in \cite{huetal2021} would still remain valid without any modification. In our case, since $\nabla_{\theta}F(q_k;\theta_k)$ is bounded in magnitude uniformly in $k$ and $\omega$ (by A1(a)), we have from (\ref{tmp4}) that
\begin{align*}
E[\|G_{k,1}\|^2]&=E[\|\nabla_{\theta}F(q_k;\theta_k) \|^2 ]+E\Big[\Big\|\frac{-I_k^++I_k^-}{2\bar c_k{\Delta}_k}-E\big[\frac{-I_k^++I_k^-}{2\bar c_k{\Delta}_k} |\mathcal{F}_k\big]\Big\|^2 \Big]+O(c^2_k)\\
&\leq E[\|\nabla_{\theta}F(q_k;\theta_k) \|^2 ]+E\Big[\Big\|\frac{-I_k^++I_k^-}{2\bar c_k{\Delta}_k}\Big\|^2\Big]+O(c_k^2)\\
&\leq E[\|\nabla_{\theta}F(q_k;\theta_k) \|^2 ]+\frac{1}{4c_k^2}E \big[\max\big\{d,{\|D_k\|^2}\big\}\big]+O(c_k^2),
\end{align*}
where the $O(c_k^2)$ term is due to part ($ii$) of Lemma~\ref{lem:bias}.
Then from A1(a), A3(a), and part ($i$) of Lemma~\ref{lem:dk}, it immediately follows that $\sum_{k=0}^{\infty}\beta_k^2E[\|G_{k,1}\|^2]<\infty$, which completes the proof. \Halmos
\endproof

\section{Proof of Lemma~\ref{lem:order}}\label{appendix_order}

Since $w(i)=O(1/i^s)$ and $u(i)=a/i^{p}$, we can find an integer $N>0$ and a constant $C>0$ such that $w(i)\leq Ci^{-s}$ and $1/i\leq u(i)< 1$ for all $i\geq N$. Thus for all $k\geq N$, we have that
\begin{align}\label{tmp5}
&\sum_{i=1}^k\Big[ \prod_{j=i+1}^k \big(1-u(j) \big)\Big]u(i)w(i) \nonumber\\
&=\sum_{i=1}^{N-1}\Big[ \prod_{j=i+1}^k \big(1-u(j) \big)\Big]u(i)w(i)+\sum_{i=N}^k\Big[ \prod_{j=i+1}^k \big(1-u(j) \big)\Big]u(i)w(i) \nonumber \\
&=\prod_{j=N}^k\big(1-u(j) \big)\sum_{i=1}^{N-1}\Big[\prod_{j=i+1}^{N-1}\big(1-u(j) \big) \Big]u(i)w(i)+\sum_{i=N}^k\Big[ \prod_{j=i+1}^k \big(1-u(j) \big)\Big]u(i)w(i) \nonumber\\
&\leq e^{-\sum_{j=N}^k u(j)}\Big|\sum_{i=1}^{N-1}\prod_{j=i+1}^{N-1}\big(1-u(j) \big)u(i)w(i)\Big|+\sum_{i=N}^k\Big[ \prod_{j=i+1}^k \big(1-u(j) \big)\Big]u(i)w(i) \nonumber\\
&=o(k^{-1})+\sum_{i=N}^k\Big[ \prod_{j=i+1}^k \big(1-u(j) \big)\Big]u(i)w(i)~~~\mbox{since $p\in(0,1)$} \nonumber \\
&\leq o(k^{-1})+C\sum_{i=N}^k\Big[ \prod_{j=i+1}^k \big(1-u(j) \big)\Big]u(i)i^{-s}.
\end{align}

Let $s_k:=\sum_{i=N}^k\Big[ \prod_{j=i+1}^k \big(1-u(j) \big)\Big]u(i)i^{-s}$ and define a mapping
$T_k(x):=\big(1-u(k)\big)x+u(k)k^{-s}$. It is clear that
\begin{align*}
s_k&=\sum_{i=N}^{k-1}\Big[ \prod_{j=i+1}^k \big(1-u(j) \big)\Big]u(i)i^{-s}+u(k)k^{-s}\\
&=\big(1-u(k)\big)s_{k-1}+ u(k)k^{-s}\\
&=T_k(s_{k-1}).
\end{align*}
In addition, since $|T_k(x)-T_k(y) |=(1-u(k))|x-y|$, $T_k(\cdot)$ is a contraction mapping and has a unique fixed point given by $s^*_k=k^{-s}$. Next, by noting that $|s^*_k-s^*_{k+1}|\leq s k^{-(s+1)}$, we obtain \begin{align*}
|s_k-s_k^*|&=|T_k(s_{k-1})-T_k(s_k^*)  |\\
&= (1-u(k))|s_{k-1}-s_k^*|\\
&\leq (1-u(k))|s_{k-1}-s_{k-1}^* |+(1-u(k))|s_k^*-s_{k-1}^* |\\
& \ldots\\
&\leq \prod_{i=N+1}^k \big(1-u(i) \big)|s_N-s_N^*|+\sum_{i=N+1}^k \Big[\prod_{j=i}^k\big(1-u(j) \big) \Big]|s_i^*-s^*_{i-1}|\\
&\leq e^{-\sum_{i=N+1}^k u(i)}|s_N-s^*_N|+\sum_{i=N+1}^k \Big[\prod_{j=i}^k\big(1-u(j) \big) \Big]\frac{s}{(i-1)^{s+1}}\\
&\leq o(k^{-1})+\sum_{i=N+1}^k \Big[\prod_{j=i}^k\big(1-\frac{1}{j} \big) \Big]\frac{s}{(i-1)^{s+1}}~~~\mbox{because $u(j)\geq 1/j$}\\
&\leq o(k^{-1})+\frac{s}{k}\sum_{i=N}^{k-1}\frac{1}{i^s}\\
&\leq o(k^{-1})+\frac{s}{1-s}k^{-s}.
\end{align*}
In view of (\ref{tmp5}), the result is hence proved by noting that $0<s_k \leq s_k^*+o(k^{-1})+\frac{s}{1-s}k^{-s}=o(k^{-1})+\frac{1}{1-s} k^{-s}=O(k^{-s})$.
\Halmos

\section{Proof of Lemma~\ref{lem:qkrate}}\label{appendix_qkrate}
Since $q(\theta)$ is twice continuously differentiable and $\Theta$ is compact, $q(\cdot)$ is Lipschitz continuous on $\Theta$. Let us denote its associated Lipschitz constant by $L_q$. Thus, we have from (\ref{projection}) that
\begin{equation}\label{tmp6}
|q(\theta_k)-q(\theta_{k+1})|\leq L_q\|\theta_k-\theta_{k+1}\|\leq \alpha_kL_q(\|D_k\|+\|Z_k\|)\leq 2\alpha_kL_q\|D_k\|,
\end{equation}
where we have used the fact that $\|Z_k\|\leq \|D_k\|$ (see the proof of Lemma 5 in \cite{huetal2021}).

For brevity, define $\zeta_k=q_k-q(\theta_k)$ and write iteration (\ref{qk}) in terms of $\zeta_k$ as
\begin{align*}
\zeta_{k+1}=\zeta_k+\gamma_k(\varphi-I\{Y_k\leq q_k \})+q(\theta_k)-q(\theta_{k+1}).
\end{align*}
Squaring both sides of the above equation and using (\ref{tmp6}) yield
\begin{align*}
\zeta_{k+1}^2&=\zeta_k^2+\gamma_k^2(\varphi-I\{Y_k\leq q_k \})^2+\big(q(\theta_k)-q(\theta_{k+1}) \big)^2+2\gamma_k\zeta_k(\varphi-I\{Y_k\leq q_k \})\\
&\hspace{0.5cm}+2\zeta_k\big(q(\theta_k)-q(\theta_{k+1}) \big)+2\gamma_k(\varphi-I\{Y_k\leq q_k \})\big(q(\theta_k)-q(\theta_{k+1})\big)\\
&\leq \zeta_k^2+\gamma_k^2+4\alpha_k^2L_q^2\|D_k\|^2+2\gamma_k\zeta_k(\varphi-I\{Y_k\leq q_k \})+4|\zeta_k|\alpha_kL_q\|D_k\|\\
&\hspace{1cm}+2\gamma_k(\varphi-I\{Y_k\leq q_k \})\big(q(\theta_k)-q(\theta_{k+1})\big).
\end{align*}
Then, by taking conditional expectation at both sides, we get
\begin{align*}
E[\zeta_{k+1}^2 |\mathcal{F}_k]&\leq \zeta_k^2+\gamma_k^2+4\alpha_k^2L_q^2\|D_k\|^2+2\gamma_k\zeta_k\Big(F\big(q(\theta_k);\theta_k\big)-F(q_k;\theta_k) \Big)+4|\zeta_k|\alpha_kL_q\|D_k\|\\
&\hspace{1cm}+2\gamma_k\Big(F\big(q(\theta_k);\theta_k\big)-F(q_k;\theta_k) \Big)\big(q(\theta_k)-q(\theta_{k+1})\big)\\
&=\zeta_k^2+\gamma_k^2+4\alpha_k^2L_q^2\|D_k\|^2-2\gamma_k f(\tilde q_k;\theta_k)\zeta_k^2+4|\zeta_k|\alpha_kL_q\|D_k\|\\
&\hspace{1cm}-2\gamma_kf(\tilde q_k;\theta_k)\zeta_k\big(q(\theta_k)-q(\theta_{k+1})\big)\\
&\leq (1-2\gamma_k\varepsilon)\zeta_k^2+\gamma_k^2+4\alpha_k^2L_q^2\|D_k\|^2+4\alpha_kL_q|\zeta_k|\|D_k\|+4\gamma_k\alpha_kC_fL_q|\zeta_k|\|D_k\|,
\end{align*}
where $\tilde q_k$ lies on the line segment between $q_k$ and $q(\theta_k)$, and the last step follows from B1 and (\ref{tmp6}).

We assume without loss of generality that $\varepsilon$ is sufficiently small such that $2\gamma_k\varepsilon<1$ for all $k$. Next by unconditioning on $\mathcal{F}_k$ and applying the Cauchy Schwarz inequality, we further obtain
\begin{align*}
E[\zeta_{k+1}^2]&\leq (1-2\gamma_k\varepsilon)E[\zeta_k^2]+4\alpha_kL_q(1+\gamma_kC_f)\sqrt{E[\zeta_k^2]E[\|D_k\|^2 ]}+\gamma_k^2+4\alpha_k^2L_q^2E[\|D_k\|^2]\\
&\leq \bigg(\sqrt{1-2\gamma_k\varepsilon}\sqrt{E[\zeta_k^2]}+\frac{2\alpha_kL_q(1+\gamma_kC_f)\sqrt{E[\|D_k\|^2]}}{\sqrt{1-2\gamma_k\varepsilon}} \bigg)^2+\gamma_k^2\\
&\leq \Big((1-\gamma_k\varepsilon)\sqrt{E[\zeta_k^2]}+\alpha_k\mathcal{M} \Big)^2+\gamma_k^2,
\end{align*}
where $\mathcal{M}:=\frac{2L_q(1+\gamma_0C_f)\sqrt{\sup_kE[\|D_k\|^2]}}{\sqrt{1-2\gamma_0\varepsilon}}<\infty$ (by part ($i$) of Lemma~\ref{lem:dk}) and we have used the inequality $\sqrt{1-x}\leq 1-x/2$ for $x\in[0,1]$.

Now define a mapping
\begin{align*}
\mathcal{T}_k(x):=\sqrt{\big((1-\gamma_k\varepsilon)x+\alpha_k\mathcal{M} \big)^2+\gamma_k^2}
\end{align*}
and consider the sequence of real numbers $\{x_k\}$ generated by
$x_{k+1}=\mathcal{T}_k(x_k)$ for all $k$ with $x_0:=\sqrt{E[\zeta_0^2]}$. A simple inductive argument shows that $\sqrt{E[\zeta_k^2]}\leq x_k$ for all $k$. In addition, notice that $\mathcal{T}_k(x)$ is a function of the form $h(x)=\sqrt{x^2+b^2},~x>0$ for some constant $b\neq 0$. Since $dh(x)/dx=x/\sqrt{x^2+b^2}<1$, we have $|h(x)-h(y)|\leq |x-y|$ for all $x,\,y>0$. This implies that $\mathcal{T}_k(\cdot)$ is a contraction mapping satisfying
\begin{equation*}
|\mathcal{T}_k(x)-\mathcal{T}_k(y)|\leq (1-\gamma_k\varepsilon)|x-y|.
\end{equation*}

The unique fixed point $x^*_k$ of $\mathcal{T}_k$ can be obtained by solving the quadratic equation
$x^2=((1-\gamma_k\varepsilon)x+\alpha_k \mathcal{M} )^2+\gamma_k^2$ and written in the form
\begin{equation}\label{fixed}
x_k^*=C_1\frac{\alpha_k}{\gamma_k}+C_2\gamma_k^{\frac{1}{2}}+h.o.t.,
\end{equation}where $C_1,\,C_2\geq0$ are two constants and $h.o.t.$ stands for higher-order terms.
For the specific forms of $\alpha_k$ and $\gamma_k$ given in Section~\ref{sec4}, it is not hard to verify that $|x_{k+1}^*-x_k^*|=O(k^{\gamma-\alpha-1})+O(k^{-\gamma/2-1})$. Consequently, we have by the contraction property of $\mathcal{T}_k$ and Lemma~\ref{lem:order},
\begin{align*}
|x_{k+1}-x_{k+1}^*|&\leq |x_{k+1}-x_k^*|+|x_k^*-x_{k+1}^*|\\
&=|\mathcal{T}_k(x_k)-\mathcal{T}_k(x_k^*) |+ |x_k^*-x_{k+1}^*|\\
&\leq (1-\gamma_k\varepsilon)|x_k-x_k^* |+ |x_k^*-x_{k+1}^*|\\
&\leq (1-\gamma_k\varepsilon)|x_k-x_{k-1}^* |+(1-\gamma_k\varepsilon)|x_k^*-x_{k-1}^* |+ |x_k^*-x_{k+1}^*|\\
&\ldots \\
&\leq \prod_{i=0}^k(1-\gamma_i\varepsilon)|x_0-x_0^* |+\sum_{i=0}^k\Big[\prod_{j=i+1}^k(1-\gamma_j\varepsilon)  \Big]\gamma_i\varepsilon \frac{|x_{i+1}^*-x_i^*|}{\gamma_i\varepsilon}\\
&=O(k^{2\gamma-\alpha-1})+O(k^{\gamma/2-1})\\
&=o(x_{k+1}^*).
 \end{align*}
Finally, combining this with (\ref{fixed}) shows that $\sqrt{E[\zeta_k^2]}\leq x_k\leq x_k^*+o(x_k^*)=O(\alpha_k/\gamma_k)+O(\gamma_k^{{1}/{2}})$, which completes the proof. \Halmos

\section{Proof of Lemma~\ref{lem:dkrate}}\label{appendix_dkrate}
Since $q(\theta)$ is twice continuously differentiable on $\Theta$, there exists a Lipschitz constant $L_g$ such that $\|\nabla_{\theta}q(\theta_k)-\nabla_{\theta}q(\theta_{k+1}) \|\leq L_g\|\theta_k-\theta_{k+1}\|\leq 2\alpha_k L_g\|D_k\|$; see, (\ref{tmp6}). Note that here for notational simplicity, we write $\nabla_{\theta}q(\theta_k)$ for $\nabla_{\theta}q(\theta)|_{\theta=\theta_k}$ and $b_k$ for $b_k(q_k,D_k,\theta_k)$.

Let $\eta_k=D_k-\nabla_{\theta}q(\theta_k)$ and write (\ref{dk}) as
\begin{align*}
\eta_{k+1}=\eta_k+\beta_k\frac{-I_k^++I^-_k}{2\bar c_k{\Delta}_k}+\nabla_{\theta}q(\theta_k)-\nabla_{\theta}q(\theta_{k+1}).
\end{align*}
Using $1/\Delta_k=\Delta_k$, $\Delta_k^T\Delta_k=d$, and the Lipschitz continuity of $\nabla_\theta q(\theta)$, it follows that
\begin{align*}
\|\eta_{k+1}\|^2&\leq \|\eta_k\|^2+\frac{\beta_k^2 M_k^2 d}{4c_k^2}+4\alpha_k^2L_g^2\|D_k\|^2+2\beta_k\eta_k^T\frac{(-I^+_k+I^-_k)}{2\bar c_k{\Delta}_k}\\
&\hspace{0.3cm}+2\big(\nabla_{\theta}q(\theta_k)-\nabla_{\theta}q(\theta_{k+1})\big)^T\eta_k+2\beta_k\big(\nabla_{\theta}q(\theta_k)-\nabla_{\theta}q(\theta_{k+1})\big)^T \frac{-I^+_k+I^-_k}{2\bar c_k{\Delta}_k}.
\end{align*}
Then, taking conditional expectations at both sides and noticing that $\theta_{k+1}$ is $\mathcal{F}_k$-measurable, we have by Lemma~\ref{lem:bias} that
\begin{align}\label{tmp7}
E[\|\eta_k\|^2 |\mathcal{F}_k]&\leq \|\eta_k\|^2+\frac{\beta_k^2 M_k^2 d}{4c_k^2}+4\alpha_k^2L_g^2\|D_k\|^2-2\beta_k\big(f(q_k;\theta_k)D_k^T\eta_k+\nabla_{\theta}^TF(q_k;\theta_k)\eta_k \big)\nonumber\\
&\hspace{0.2cm}+2\beta_kb_k^T\eta_k+2\big(\nabla_{\theta}q(\theta_k)-\nabla_{\theta}q(\theta_{k+1})\big)^T\eta_k \nonumber\\
&\hspace{0.2cm}-2\beta_k\big(\nabla_{\theta}q(\theta_k)-\nabla_{\theta}q(\theta_{k+1})\big)^T\big(f(q_k;\theta_k)D_k+\nabla_{\theta}F(q_k;\theta_k) \big)\nonumber\\
&\hspace{1cm}+2\beta_kb_k^T\big(\nabla_{\theta}q(\theta_k)-\nabla_{\theta}q(\theta_{k+1})\big)\nonumber\\
&\hspace{-2cm}=\|\eta_k\|^2+\frac{\beta_k^2 M_k^2 d}{4c_k^2}+4\alpha_k^2L_g^2\|D_k\|^2-2\beta_kf(q_k;\theta_k)\big(D_k^T-\nabla_{\theta}^Tq(\theta_k) \big)\eta_k \nonumber\\
&\hspace{-1.5cm}+2\beta_k f(q_k;\theta_k)\Big(-\nabla_{\theta}^Tq(\theta_k)-\frac{\nabla_{\theta}^TF(q_k;\theta_k)}{f(q_k;\theta_k)} \Big)\eta_k+2\beta_kb_k^T\eta_k \nonumber\\
&\hspace{-1.5cm}+2\big(\nabla_{\theta}q(\theta_k)-\nabla_{\theta}q(\theta_{k+1})\big)^T\eta_k-2\beta_kf(q_k;\theta_k)(D_k^T-\nabla_{\theta}^Tq(\theta_k))\big(\nabla_{\theta}q(\theta_k)-\nabla_{\theta}q(\theta_{k+1}) \big)\nonumber \\
&\hspace{-1.5cm}+2\beta_kf(q_k;\theta_k)\Big(-\nabla_{\theta}^Tq(\theta_k)-\frac{\nabla_{\theta}^TF(q_k;\theta_k)}{f(q_k;\theta_k)} \Big)\big(\nabla_{\theta}q(\theta_k)-\nabla_{\theta}q(\theta_{k+1}) \big) \nonumber \\
&\hspace{-1.5cm}+2\beta_kb_k^T\big(\nabla_{\theta}q(\theta_k)-\nabla_{\theta}q(\theta_{k+1})\big)\nonumber\\
&\hspace{-2cm}\leq(1-2\beta_k\epsilon)\|\eta_k\|^2+\frac{\beta_k^2 M_k^2 d}{4c_k^2}+4\alpha_k^2L_g^2\|D_k\|^2-2\beta_kf(q_k;\theta_k)\mu_{k,1}^T\eta_k+2\beta_kb_k^T\eta_k \nonumber\\
&\hspace{-1.5cm}+2\mu_{k,2}^T\eta_k-2\beta_kf(q_k;\theta_k)\mu_{k,2}^T\eta_k-2\beta_kf(q_k;\theta_k)\mu_{k,1}^T\mu_{k,2}+2\beta_kb_k^T\mu_{k,2},
\end{align}
where the last inequality follows from A1(b) and we have defined $\mu_{k,1}=\nabla_{\theta}q(\theta_k)+\frac{\nabla_{\theta}F(q_k;\theta_k)}{f(q_k;\theta_k)}$ and $\mu_{k,2}=\nabla_{\theta}q(\theta_k)-\nabla_{\theta}q(\theta_{k+1})$.

Note that by (\ref{gd0}),  $\nabla_{\theta}q(\theta_k)=-\frac{\nabla_{\theta}F(y;\theta_k)|_{y=q(\theta_k)}}{f(q(\theta_k);\theta_k)}$.
Due to the compactness of $\Theta$, the continuity of $q(\theta)$ and $f(y;\theta)$ (Assumption B2) implies that $f(q(\theta_k);\theta_k)$ is bounded from below uniformly in $\theta_k$. This, together with B1 and the Lipschitz conditions in B2, indicates that there must be a constant $\mathcal{L}$ such that $\|\mu_{k,1}\|\leq \mathcal{L}|q_k-q(\theta_k) |=\mathcal{L}|\zeta_k|$, where recall that $\zeta_k=q_k-q(\theta_k)$.
In addition, $\|\mu_{k,2} \|\leq 2 \alpha_k L_g\|D_k\|$.
On the other hand, from the proof of Lemma~\ref{lem:bias}, it can be seen that
$\sqrt{E[\|b_k\|^2]}=O(c_k^2)$.

Next, taking expectations at both sides of (\ref{tmp7}) and applying the Cauchy–Schwarz inequality, we further obtain
\begin{align*}
E[|\eta_{k+1}\|^2]&\leq (1-2\beta_k\epsilon)E[\|\eta_k\|^2]+\frac{\beta_k^2 M_k^2 d}{4c_k^2}+4\alpha_k^2L_g^2E[\|D_k\|^2]+2\beta_kC_f\mathcal{L}\sqrt{E[\zeta_k^2 ]}\sqrt{E[\|\eta_k \|^2]}\\
&+2\beta_k\sqrt{E[\|b_k \|^2 ]}\sqrt{E[\|\eta_k \|^2 ]}
+4(1+\beta_kC_f)\alpha_kL_g\sqrt{E[\|D_k\|^2 ]}\sqrt{E[\|\eta_k\|^2]}\\
&+4\alpha_k\beta_kC_fL_g\mathcal{L}\sqrt{E[\zeta_k^2]}\sqrt{E[\|D_k \|^2]}+O(\alpha_k\beta_kc_k^2)\\
&\hspace{-2cm}=(1-2\beta_k\epsilon)E[\|\eta_k\|^2]+O \Big(\frac{\beta_k^2}{c_k^2}\Big)\\
&\hspace{-1.5cm}+2\Big(\beta_kC_f\mathcal{L}\sqrt{E[\zeta_k^2]}+
\beta_k\sqrt{E[\|b_k \|^2 ]}
+2\alpha_k(1+\beta_kC_f)L_g\sqrt{E[\|D_k \|^2 ]} \Big)\sqrt{E[\|\eta_k \|^2 ]},
\end{align*}
where the last step above follows from the observations that: (1) $4\alpha_k^2L_g^2E[\|D_k\|^2]=o(\beta_k^2/c_k^2)$; (2) $\alpha_k\beta_kc_k^2=o(\beta_k^2/c_k^2)$; (3) $4\alpha_k\beta_kC_fL_g\mathcal{L}\sqrt{E[\| \zeta_k\|^2]}\sqrt{E[\|D_k \|^2]}=O(\alpha_k^2\beta_k/\gamma_k )+O(\alpha_k\beta_k\gamma_k^{{1}/{2}})$ (by Lemmas~\ref{lem:dk} and \ref{lem:qkrate}), which is of a higher order than $O(\beta_k^2/c_k^2)$.

Let $\mathcal{A}_k=\big(\beta_kC_f\mathcal{L}\sqrt{E[\zeta_k^2]}+\beta_k\sqrt{E[ \| b_k\|^2]}+2\alpha_k(1+\beta_kC_f)L_g\sqrt{E[\|D_k \|^2 ]}\big)/\sqrt{1-2\beta_k\epsilon}$. By Lemmas~\ref{lem:dk} and \ref{lem:qkrate}, $\mathcal{A}_k=O(\alpha_k\beta_k/\gamma_k)+O(\beta_k\gamma_k^{1/2})+O(\beta_kc_k^2)+O(\alpha_k)=O(\alpha_k\beta_k/\gamma_k)+O(\beta_k\gamma_k^{1/2})+O(\beta_kc_k^2)$ (due to A3(d)) and $\mathcal{A}_k^2=o(\beta_k^2/c_k^2)$. Thus, we arrive at the following bound:
\begin{align*}
E[\|\eta_{k+1}\|^2 ]&\leq \Big(\sqrt{1-2\beta_k\epsilon}\sqrt{E[\|\eta_k \|^2 ]}+\mathcal{A}_k \Big)^2+O\Big(\frac{\beta_k^2}{c_k^2}\Big)\\
&\leq \big((1-\beta_k\epsilon)\sqrt{E[\|\eta_k \|^2 ]}+\mathcal{A}_k \big)^2+O\Big(\frac{\beta_k^2}{c_k^2}\Big).
\end{align*}

Let $\mathcal{N}>0$ be an integer such that $\mathcal{A}_k\leq \mathcal{C}(\alpha_k\beta_k/\gamma_k+\beta_k\gamma_k^{1/2}+\beta_k c_k^2)$ for all $k\geq\mathcal{N}$ for some constant $\mathcal{C}>0$.
We then proceed as in the proof of Lemma~\ref{lem:qkrate} and define a mapping
$$\mathcal{G}_k(x)=\sqrt{\big((1-\beta_k\epsilon)x+\bar{\mathcal{A}}_k\big)^2+O\big(\beta_k^2/c_k^2 \big)},$$ where $\bar{\mathcal{A}}_k=\mathcal{C}(\alpha_k\beta_k/\gamma_k+\beta_k\gamma_k^{1/2}+\beta_k c_k^2)$.
Clearly, if we take $x_{\mathcal{N}}=\sqrt{E[\|\eta_{\mathcal{N}} \|^2 ]}$, then $\sqrt{E[\|\eta_k \|^2 ]}\leq x_k$ for all $x_k$ generated by $x_{k+1}=\mathcal{G}_k(x_k)$, $k\geq \mathcal{N}$. In addition, $\mathcal{G}_k(\cdot)$ is a contraction and its unique fixed point is of order {$x_k^*=O(\bar{\mathcal{A}}_k/\beta_k)+O(\beta_k^{1/2}/c_k)=O(\alpha_k/\gamma_k)+O(\gamma_k^{1/2})+O(c_k^2)+O(\beta_k^{1/2}/c_k)
=O(\alpha_k/\gamma_k)+O(c_k^2)+O(\beta_k^{1/2}/c_k)$, where the last step follows because $\gamma_k^{1/2}=o(\beta_k^{1/2})=o(\beta_k^{1/2}/c_k)$ by A3(d).} The rest of the proof amounts to showing $|x_k-x_k^*|=o(x_k^*)$ so that $x_k$ has the same order as $x_k^*$. This follows straightforwardly from the same contraction argument used in the proof of Lemma~\ref{lem:qkrate}. \Halmos
\end{APPENDICES}

\end{document}